\documentclass{article}
\usepackage[utf8]{inputenc}
\usepackage[english]{babel}
\usepackage{mathtools,amsmath,amssymb,amsbsy}
\usepackage[normalem]{ulem} 
\usepackage{enumitem}
\usepackage{fancyhdr}
\usepackage{csquotes}
\usepackage{geometry}
\usepackage{graphicx}
\usepackage{microtype}
\usepackage{graphicx}
\usepackage{subfigure}
\usepackage{booktabs} 

\usepackage{fancyhdr}
\usepackage{ntheorem}
\usepackage{times}  
\usepackage{helvet}  
\usepackage{courier}  
\usepackage{url}  
\usepackage{booktabs}       
\usepackage{amsfonts}       
\usepackage{nicefrac}       
\usepackage{microtype}      
\newtheorem{theorem}{Theorem}[section]
\newtheorem{remark}{Remark}[section]
\newtheorem{proposition}{Proposition}[section]
\newtheorem{definition}{Definition}[section]
\newtheorem{corollary}{Corollary}[section]
\newtheorem{lemma}{Lemma}[section]

\newtheorem{example}[theorem]{Example}

\theoremstyle{empty}

\newtheorem{refproof}{Proof}[section]

\usepackage[
backend=bibtex,
style=alphabetic,
sorting=ynt
]{biblatex}

\title{Optimal Capital Injections with the Risk of Ruin: A Stochastic Differential Game of Impulse Control and Stopping Approach
  \newline \newline}

\author{David Mguni\footnote{ Quantitative and Applied Spatial Economic Research Laboratory, University College London, Gower Street, London, WC1E 6BT, UK.}\hspace{1.5 mm}\footnote{Centre for Doctoral Training in Financial Computing \& Analytics, University College London, Gower Street, London, WC1E 6BT, UK. davidmguni@hotmail.com} }

\date{}

\begin{document}

\maketitle
\vspace{0.05 cm}
\hspace{-6.6 mm} \begin{abstract}
We consider an investment problem in which an investor performs capital injections to increase the liquidity of a firm for it to maximise profit from market operations. Each time the investor performs an injection, the investor incurs a fixed transaction cost. In addition to maximising their terminal reward, the investor seeks to minimise risk of loss of their investment (from a possible firm ruin) by exiting the market at some point in time. We show that the problem can be reformulated in terms of a new  stochastic differential game of control and stopping in which one of the players modifies a (jump-)diffusion process using impulse controls and an adversary chooses a stopping time to end the game. We show that the value of this game can be computed by solving a double obstacle problem described by a quasi-variational inequality. We then characterise the value of the game via a set of HJBI equations, considering both games with zero-sum and non-zero-sum payoff structures. Our last result demonstrates that the solution to the investment problem is recoverable from the Nash equilibrium strategies of the game.
\end{abstract}

\maketitle
\section{Introduction}

There are numerous environments in which financial agents incur fixed or minimal costs when adjusting their investment positions; trading environments with transaction costs, real options pricing and real estate and large-scale infrastructure investing are a few important examples. The study of optimal investments by an economic agent who seeks to minimise the probability that they go bankrupt within their lifetime is known as the probability of lifetime ruin problem. The problem was introduced by \cite{[1]} and studied in depth by \cite{[2]}, for a complete treatment of the background and origins of this problem, we refer the reader to \cite{[15],[16],[17]} and references therein.

Despite the breadth of the literature concerning the lifetime ruin problem and  the significant effect of transaction costs on investment behaviour, current models within the literature have yet to include those in which the investor faces financial transaction costs.\footnote{One exception is \cite{[10]} in which a single-controller problem is analysed --- however, in \cite{[10]} the controller's action space is limited to two actions wherein the model in \cite{[10]} can thus be viewed as \textit{regime switching model} with switching costs.} Moreover, the absence of transaction costs within the theoretical analysis limits the scope of application of the lifetime ruin model within financial systems. The objective of this paper is therefore to generalise the lifetime ruin problem to problems the investor now faces transaction costs when modifying their position.\

In order to tackle this problem, we introduce a new stochastic differential game of control and stopping in which the controller uses impulse controls to modify the dynamics of a (jump-)diffusion process. In contrast to existing controller and stopper games, the game we study involves a controller that now faces control costs that are bounded from beneath. The inclusion of fixed minimal adjustment costs results in a game that is an appropriate modelling framework for analysing financial applications in which we are required to find both an optimal (market) exit criterion and, an optimal investment strategy (over a time horizon) when the investor faces transaction costs. In this setting, interdependencies between the two objectives, that is, finding an optimal exit time and finding an optimal investment control strategy produce a two-player game structure in which each control strategy takes consideration of the other.  

We perform our analysis in sympathy with an optimal firm liquidity control problem with lifetime ruin in which an investor performs capital injections to increase the liquidity of a firm. In this instance, the investor seeks to maximise his capital injections to buoy the firm's liquidity process whilst seeking to minimise the probability of loss of investment by exiting the market (selling all firm holdings). The results found in the paper are however general to the problem of impulse control and discretionary stopping in strategic settings. 

We analyse the stochastic game and provide a full characterisation of the value of the game as a solution to a PDE, namely a Hamilton-Jacobi-Bellman-Isaacs (HJBI) equation. We then extend the results to cover the game with a non-zero-sum payoff structure hence providing a characterisation of the Nash equilibrium of the non-zero-sum game. To complete our analysis, we lastly show that the solution to the optimal liquidity control and lifetime ruin with transaction costs can be recovered from the equilibrium controls of the non-zero-sum stochastic differential game. This enables the solution to be extracted from solutions to a joint set of PDEs.
\subsection*{Theoretical Background}

The mathematical framework that underpins the investment problem is a stochastic differential game which combines impulse controls with optimal stopping. Impulse control problems are stochastic control models in which the cost of control is bounded below by some fixed positive constant which prohibits continuous control. This augments the problem to one of finding both an optimal sequence of times to apply the control policy, in addition to determining optimal control magnitudes. We refer the reader to \cite{bensoussan1982controle} as a general reference to impulse control theory and to \cite{korn1999some,ma2008impulse,[8]} for articles on applications. Additionally, matters relating to the application of impulse control models have been surveyed extensively in \cite{[9]}. Impulse control frameworks therefore underpin the description of financial environments with transaction costs and liquidity risks and more generally, applications of optimal control theory in which the system dynamics are modified by a sequence of discrete actions.

Stochastic differential games with impulse control (in which two players modify the system dynamics) have recently appeared in the stochastic impulse control literature. Deterministic versions of a game in were first studied by \cite{[11],[12]} --- in the model presented in \cite{[11]}, impulse controls are restricted to use by one player and the other uses continuous control. Similarly, in \cite{[13]} stochastic differential games in which one player uses impulse control and the other uses continuous controls were studied. Using a verification argument, the conditions under which the value of game is a solution to a HJBI equation is also shown in \cite{[13]}. In \cite{[14]}, a stochastic differential game in which both players use impulse control is analysed using viscosity theory.

Problems that combine both discretionary stopping and stochastic optimal control have attracted much attention over recent years; in particular there is a notable amount of literature on models of this kind in which a single controller uses continuous controls to modify the system dynamics. Discretionary stopping and stochastic optimal control problems in which the controller exercises modifications through the drift component of the state process (using continuous controls) have been studied by \cite{[22],[23],[24],[25]}. Another version of these problems which has attracted significant interest is problems in which the controller acts to modify the system dynamics by finite variations of the state process --- such problems have been studied by \cite{[26],[27]}.

A related family of models has recently emerged in which the task of controlling the system dynamics and exit time is divided between two players who act according to separated interests \cite{[15],[16]}. Controller-Stopper games were introduced by Maitra \& Sudderth in \cite{[36]}; recent papers on the matter include \cite{[37]} in which a game with underlying system dynamics are given by a one-dimensional diffusion within a given interval in $\mathbb{R}$ is studied. Other papers on the topic include \cite{[37]} and \cite{[39]}; in the latter, a multidimensional model is studied wherein the state process is controlled on a diffusion in a multidimensional Euclidean space.  A game-theoretic approach to stochastic optimal control problems with discretionary stopping has been used to analyse the lifetime ruin problem in \cite{[3],[15],[16]} amongst others.\

Within these models, the task of controlling the investment process and selecting the market exit time is assigned to two individual players who each seek to maximise some form of the same objective payoff functional. Game-theoretic formulations of the optimal stochastic control with discretionary stopping model can be viewed as  generalised versions of the single controller models wherein the investor is now allowed to seek multiple objectives which are each defined over multiple payoff functions.\

Within the body of literature concerning stochastic differential games of control and stopping however, the set of controller is restricted to a continuous class of controls (e.g. \cite{[4],[5],[36],[37],[38]}). This renders the aforementioned models unsuitable for prescribing solutions for investment problems with fixed minimal costs as continuous adjustments would result in immediate ruin. Our goal in this paper is to address the absence of impulse controls within stochastic differential games of control and stopping.

\subsection*{Organisation}

The paper is organised as follows: in Sec. \ref{section_investment_problem}, we give a complete description of  the optimal liquidity control problem and construct the main investment model of the paper. In Sec. \ref{section_lit_review}, we perform a literature review covering the related literature within stochastic differential game theory and the study of optimal liquidity control. In Sec. \ref{section_game_description}, we give a technical description of the game and introduce some of the underlying concepts for the analysis of the game. In Sec. \ref{section_main_results}, we give a statement of the main results of the paper. In Sec. \ref{section_zero_sum_game}, we perform the main study of game within a zero-sum setup starting with a  progressive development of a set of arguments that describe the equilibrium conditions. Following that, we then give  a full characterisation of the value function for the game in terms of a verification argument (Theorem \ref{Verification_theorem_for_Zero-Sum_Stochastic_Differential_Games_of_Control_and_Stopping_with_Impulse_Controls}). In Sec. \ref{section_non_zero_sum_game}, we extend the analysis to non zero-sum games, firstly constructing a definition of the appropriate equilibrium concept in this setting. This is followed by some background arguments leading to a corresponding verification theorem for non zero-sum stochastic differential games of control and stopping (Theorem \ref{Verification theorem for Non-Zero-Sum Stochastic Differential Games of Control and Stopping with Impulse Controls}). In Section \ref{section_examples}, we apply the theory to solve examples drawn from financial settings and extract the optimal investment strategy for the optimal liquidity and lifetime ruin problem of Section \ref{section_investment_problem} by way of characterising the equilibrium policies of the game.

\section{Investment Problem} \label{section_investment_problem}


Consider an investor that can choose a set of times to inject capital into a firm which increases the firm's market capabilities. Each time the investor performs a capital injection, the investor faces a transaction cost. The investor may exit the market by selling all holdings in the firm, moreover at any point, the firm may go into ruin at which point the investor faces loss of their firm investments. Thus the investor seeks to maximise their terminal returns by performing the maximal sequence of capital injections at selected times that their wealth process can tolerate. However, the investor also seeks an optimal time to exit the market by selling all firm holdings before firm bankruptcy.

If the investor seeks to maximise their reward and exit at a time that minimises some risk criterion, what is the sequence of capital injections and at what time should the investor exit the market in order to minimise the risk of loss of investment? 

In order to address this question, we seek both an  optimal control process of the investor's capital injections and, an optimal time to exit the market in advance of possible firm ruin (\textit{under a risk-minimising policy}). In the problem we study, the investor faces transaction costs so that at each capital injection incurs some fixed minimal cost. The problem combines two distinct problems; the first problem addresses the optimal sequence of capital injections and times to perform such injections. The second problem is concerned with finding an optimal exit criterion in order to minimise a notion of risk. 
We now provide a formal description of the problem facing the investor. Some of the ideas for the following description of the problem are loosely adapted from the (continuous control) descriptions of problem presented in \cite{[3],[18]}.

\subsubsection*{Description of The Problem}

The firm's liquidity at time $t_0\leq s\leq T$ is described by a stochastic process $X^{t_0,x_0}_s=X(s,\omega ):[0,T]\times \Omega \to \mathbb{R}$ over some time horizon $T\in ]0,+\infty [$ where $t_0\in[0,T]$ and $x_0\in\mathbb{R}_{>0}$ are given parameters that describe the start time of the problem and the firm's initial surplus respectively.

When there are no capital injections, the firm's liquidity process evolves according to the following expression:
\begin{equation}
dX^{t_0,x_0}_s=erX^{t_0,x_0}_{s}ds+\sigma_f X^{t_0,x_0}_{s}dB_f(s)+S_f(X^{t_0,x_0}_s,s),\qquad X_{t_0}^{t_0,x_0}:=x_0; \qquad  \mathbb{P}-{\rm a.s}., \label{ch2firmpassiveliqprocess}
\end{equation}
where without loss of generality we assume that $X^{t_0,x_0,\cdot}_s=x_0$ for any $s\leq t_0$. The constant $e\in\mathbb{R}_{>0}$ that describes the firm's rate of expenditure and $r\in ]0,1[$ is the firm's return on capital. The term $S_f$ captures the exogenous shocks in the firm's liquidity process and is given by $S_f(X_s^{t_0,x_0},s):=\int\gamma_f(s,z)X_{s^-}^{t_0,x_0} \tilde{N}_f(ds,dz)$  where $\tilde{N}_f(ds,dz)\equiv N_f(ds,dz)-\nu(dz)ds$ is a compensated Poisson random measure, $N_f(ds,dz)$ is a jump measure and  $\nu(\cdot):= \mathbb{E}[N(]0,1],V)]$ for $V\subset\mathbb{R}\backslash\{0\}$, and $B_f(s)$ is a $1-$dimensional standard Brownian motion. The constant $\sigma_f>0$ and the function $\gamma_f:[0,T]\times\mathbb{R}\to\mathbb{R}$  describe the volatility and the jump-amplitude of the firm's liquidity process  respectively. The firm's liquidity process therefore follows a \textit{geometric L\'evy process}. Geometric L\'evy processes are widely used to model financial processes due to their close empirical fit \cite{chanlevy2004}.
    
    Each time $\tau\leq T$ the investor performs a capital injection, the investor incurs a cost given by $c (\tau ,z ):=\exp^{-\delta \tau}(\kappa_I+(1+\lambda)z)$, where $\kappa_I\in\mathbb{R}_{>0}$ is a fixed transaction cost and the parameter $\lambda\in\mathbb{R}_{>0}$  determines the proportional cost for an injection of size $z\in\mathbb{R}_{>0}$ and $\delta\in ]0,1]$ is the investor's discount rate. Since performing continuous actions would result in immediate bankruptcy, the investor's capital injections must be performed over a discrete sequence of investments.

The investor therefore  performs a sequence of capital injections $\{z_k \}_{k\in \mathbb{N}}$ over the horizon of the problem which are performed over a sequence of \textit{intervention times} $\{\tau_k (\omega)\}_{k\in \mathbb{N}}$. We denote the investor's control by the double sequence $(\tau,Z )\equiv \sum_{j\in\mathbb{N}}z_j\cdot1_{\{t_0<\tau_j\leq T\}}\in \Phi$ where $\mathcal{Z}$ is a feasible set of investor capital injections and $\mathcal{T}$ is the set of intervention times and  $\Phi\subseteq\mathcal{T}\times\mathcal{Z}$. 

Denote by $T_s^{(\tau,Z)}:=\sum_{m\geq 1}z_m\cdot 1_{\{t_0<\tau_m\leq s\}}$ the investor's capital injections process. Since the investor's capital injections are transferred to the firm, the firm's liquidity with capital injections at time $t_0< s\leq T$ is given by the following expression:
\begin{align}\nonumber
X_s^{t_0,x_0,(\tau, Z) }=x_0+\int_{t_0}^{s\wedge \rho}erX_r^{t_0,x_0,(\tau, Z) }dr&+\int_{t_0}^{s\wedge \rho}\sigma_f X_r^{t_0,x_0,(\tau, Z) } dB_f(r) 
\\&\begin{aligned}+T_{s\wedge\rho}^{(\tau,Z)}
+\int\int_{t_0}^{s\wedge \rho}\gamma_f(r,z)X_{r-}^{t_0,x_0,(\tau, Z) } \tilde{N}_f(dr,dz),&
\\ \mathbb{P}-{\rm a.s.},&	\label{ch2firmliquidityprocess} 
\end{aligned}
\end{align}
where $\rho\in\mathcal{T}$ is an $\mathcal{F}-$measurable stopping time which will be defined shortly.

In order to complete the description of the investor's problem, we construct the notion of risk of ruin facing the investor. As in \cite{baghery2013optimal} and in the sense given by \cite{artzner1999coherent,follmer2002convex}, let $\theta$ be a convex risk measure acting on a stochastic process $X$, then we can write the risk measure associated to the problem as (see Theorem 2.9 in \cite{Roccioletti2016}):
\begin{equation}
	\theta (X)=\sup_{\mathbb{Q}\in \mathcal{M}_a }\mathbb{E}_\mathbb{Q} [-X]-\chi (\mathbb{Q}), \label{convex_risk_measure} \end{equation}
where $\mathcal{M}_a$ is some family of equivalent measures\footnote{The measure $\mathbb{Q}$ is said to be equivalent or \textit{absolutely continuous} w.r.t. the measure $\mathbb{P}$ iff the null set of $\mathbb{Q}$ is a proper subset of the null set of $\mathbb{P}$. We denote the equivalence of $\mathbb{Q}$ w.r.t. the measure $\mathbb{P}$ by $\mathbb{Q}\ll\mathbb{P}$.} i.e. $\mathbb{Q}\ll\mathbb{P}$ and where $\mathbb{E}_\mathbb{Q}$ denotes the expectation w.r.t. the measure $\mathbb{Q}\in \mathcal{M}_a$ and $\chi :\mathcal{M}_a\to \mathbb{R}$ is some convex (penalty) function. 

Since the investor seeks to minimise the risk of null returns, the investor seeks to exit the market by selling all holdings at a point $\rho(\Omega  )\in\mathcal{T}$ that minimises the risk $\theta (X) $ of the investor's returns falling below $0$ (after firm ruin) before $T$ where $\mathcal{T}$ is a set of $\mathcal{F}-$measurable stopping times.

We now observe that since the investor seeks to exit the market in advance of firm ruin, the investor's optimal stopping problem admits the following representation:
\begin{equation}
\inf_{\rho\in\mathcal{T}}\left[\sup_{\mathbb{Q}\in \mathcal{M}_a} \mathbb{E}_{\mathbb{Q}}[-X(\rho)]-\chi (\mathbb{Q})\right], \label{ch2divprobgame1}
\end{equation}
where $X(\rho)$ denotes the process (\ref{ch2firmliquidityprocess}) stopped at time $\rho\in\mathcal{T}$ \footnote{We observe that the problem in (\ref{ch2divprobgame1}) can be viewed as a zero-sum game between two players; namely a player that controls the measure $\mathbb{Q}$ which may be viewed as an adverse market and the investor who selects the stopping time $\rho\in\mathcal{T}$. Games of this type are explored in \cite{bayraktar2013multidimensional} and \cite{mataramvura2008risk}.}$^{,}$\footnote{We shall hereon specialise to the case $\chi\equiv0$ in which case the risk measure $\theta$ is called \textit{coherent}.}. The convex risk measure that appears in (\ref{convex_risk_measure}) and the subsequent minimax structure appearing in (\ref{ch2divprobgame1}) are closely related to notions of \textit{robust Bayesian control} and \textit{entropy maximisation control}  (see for example, \cite{grunwald2004game,Roccioletti2016} for exhaustive discussions).

The firm's liquidity process is therefore raised by capital injections performed by the investor, however in performing capital injections, the investor's wealth is reduced since liquidity is transferred from investor to firm. The investor however, receives a return on capital through some running stream and some terminal reward after liquidating all holdings in the firm.

The investor's wealth at time $s\leq T$ is described by a stochastic process $Y_s=Y(s,\omega ):[0,T]\times\Omega\to\mathbb{R}$. Denote by  $\pi\in [0,1]$, the portion of the investor's wealth invested in risky assets and by  $\bar{T}^{(\tau,Z)}_s:=\sum_{m\geq 1}\exp^{-\delta \tau_m}[(1+\lambda)z_m +\kappa_I]1_{\{t_0<\tau_m\leq\rho \wedge s\}}$, the total deductions from the investor's wealth process due to the injections; the process $Y_s$ is hence given by the following:
\begin{align}
Y_s^{t_0,x_0,(\tau, Z) }=y_0+\int_{t_0}^{s\wedge \rho}\Gamma Y_r^{t_0,x_0,(\tau, Z) }dr-\bar{T}_{s\wedge\rho}^{(\tau,Z)}&+\int_{t_0}^{s\wedge \rho}\pi\sigma_I Y_r^{t_0,x_0,(\tau, Z) } dB_I(r)\nonumber
\\&\begin{aligned}+\int\int_{t_0}^{s\wedge \rho}\pi\gamma_I(r,z)Y_{r-}^{t_0,x_0,(\tau, Z) } \tilde{N}_I(dr,dz),&
\\Y_{t_0}^{t_0,y_0,(\tau,Z)}\equiv y_0,\; \label{ch2invliquidityprocess}\; \mathbb{P}-{\rm a.s}.,&
\end{aligned}
\end{align}
where $y_0\in\mathbb{R}_{>0}$ is the investor's initial wealth; we assume that $Y^{t_0,y_0,\cdot}_s=y_0$ for any $s\leq t_0$. The constant $\Gamma$ is given by $\Gamma:=(1-\pi)r_0+\pi\mu_R$ where $r_0, \mu_R\in\mathbb{R}$ are constants that describe the interest rate and the return on the risky assets respectively. The term $\tilde{N}_I(dr,dz)\equiv N_I(dr,dz)-\nu(dz)dr$ is a compensated Poisson random measure where $N_I(dr,dz)$ is a jump measure and $B_I(r)$ is a  $1-$dimensional standard Brownian motion. The constant $\sigma_I>0$ and the function $\gamma_I:[0,T]\times\mathbb{R}\to\mathbb{R}$ describe the volatility and the jump-amplitude of the investor's wealth process.

If we now interpret optimality of the stopping time $\rho(\Omega  )$ in a sense of risk-minimal w.r.t. the risk measure $\theta$, we can reformulate the problem in (\ref{ch2divprobgame1}) and the investor's maximisation problem in terms of a decoupled pair of objective functions. Focusing firstly on the investor's capital injection problem, we can write the problem as:

Find an admissible strategy $(\hat{\tau},\hat{Z})\in \Phi$ s.th.
\begin{align}
(\hat{\tau},\hat{Z})\in \underset{(\tau, Z) \in \Phi}{\arg \sup}\quad J_I^{(1)} (t_0,x_0,y_0,(\tau, Z),\rho),
\end{align}
where
\begin{align}
J_I^{(1)} (t_0,x_0,y_0,(\tau, Z),\rho)=  \mathbb{E}\left[\sum_{m\geq 1}e^{-\delta \tau_m}z_m \cdot 1_{\{t_0< \tau_m\leq\tau_S \wedge \rho \}}+g_2e^{-\delta (\tau_S\wedge \rho)}{Y}_{\tau_S \wedge \rho }^{t_0,y_0,(\tau, Z) }\right],&\label{ch2divinvestorproblemcontrol}	
\\\forall(t_0,x_0,y_0)\in [0,T]\times \mathbb{R}_{>0}\times \mathbb{R}_{>0}, \forall\rho\in \mathcal{T},&\nonumber
\end{align}
where $g_2\in [0,1]$ is a constant that represents the fraction of the firm's wealth returned to the investor upon exit and $\tau_S:= \inf\{s\in [0,T]: X_s,Y_s\leq 0\}\wedge T$.

Therefore, the above set of expressions describe the investor's objective to maximise the injections to the firm over the horizon of the problem when the investor's capital injections are performed over a sequence of discrete wealth transfers.  

We now turn to describing the component of the investor's objective to exit the market at an optimal time. The following expression represents the investor's optimal stopping problem:

Find an admissible strategy $\hat{\rho}\in \mathcal{T}$ s.th.
\begin{align}
\hat{\rho}\in \underset{\rho \in \mathcal{T}}{\inf}\underset{\mathbb{Q}\in\mathcal{M}_a}{\sup}J_I^{(2)} (t_0,x_0,y_0,(\tau, Z),{\rho}),
\end{align}
where
\begin{align} 
J_I^{(2)} (t_0,x_0,y_0,(\tau, Z),{\rho})=-\mathbb{E}_{\mathbb{Q}}\left[e^{-\delta (\tau_S\wedge\rho)}\left(g_1X_{\tau_S \wedge \rho }^{t_0,x_0,(\tau, Z) }+\lambda_T\right)\right]&,\label{ch2divinvestorproblemtime} 
\\
\forall(t_0,x_0,y_0)\in [0,T]\times \mathbb{R}_{>0}\times \mathbb{R}_{>0}, &\forall(\tau,Z)\in\Phi,	\nonumber
\end{align}
where $g_1\in ]0,1]$ and $\lambda_T\geq0$ represent the fraction of the firm's liquidity process and some fixed amount each received by the investor upon exit respectively.

The expressions (\ref{ch2divinvestorproblemcontrol}) and (\ref{ch2divinvestorproblemtime}) fully express the investor's set of objectives. We can combine the expressions (\ref{ch2divinvestorproblemcontrol}) and (\ref{ch2divinvestorproblemtime}) to construct a single problem with an objective function $\Pi$  given by the following:

Find an admissible strategy $(\hat{\rho},(\hat{\tau},\hat{Z}))\in \mathcal{T}\times \Phi$ s.th.
 \begin{align}
 \hat{\rho}&\in \underset{\rho \in \mathcal{T}}{\arg\inf}\hspace{1 mm}\Pi (t_0,x_0,y_0,(\hat{\tau}, \hat{Z}) ,\rho), \label{investor_stopping_problem}\\
  (\hat{\tau},\hat{Z})&\in \underset{(\tau,Z)\in\Phi}{\arg\sup}\hspace{1 mm}\Pi (t_0,x_0,y_0,(\tau, Z) ,\hat{\rho}), \qquad \forall(t_0,x_0,y_0)\in [0,T]\times \mathbb{R}_{>0}\times \mathbb{R}_{>0}, \label{investor_control_problem}
 \end{align}
 where
 \begin{align}\nonumber
 \Pi (t_0,x_0,y_0,(\tau, Z) ,\rho)=\mathbb{E}\Bigg[\underset{\mathbb{Q}\in\mathcal{M}_a}{\sup}&\mathbb{E}_{\mathbb{Q}}\Big[-e^{-\delta (\tau_S\wedge\rho)}\left(g_1X_{\tau_S \wedge \rho }^{t_0,x_0,(\tau, Z) }+\lambda_T\right)\Big] 
 \\
 &+\sum_{m\geq 1}e^{-\delta \tau_m} z_m \cdot 1_{\{t_0<\tau_m\leq\tau_S \wedge \rho \}}+g_2e^{-\delta (\tau_S\wedge \rho)}{Y}_{\tau_S \wedge \rho }^{t_0,y_0,(\tau, Z) }\Bigg].\label{ch2divinvestorpayoff}	
\end{align}

It can now be seen that the problem is to find the interdependent set of controls $(\hat{\rho},(\hat{\tau},\hat{Z}))\in \mathcal{T}\times \Phi$. If we think of the two objectives (\ref{ch2divinvestorproblemcontrol}) and (\ref{ch2divinvestorproblemtime}) as being assigned to two individual players, we recognise the pair of problems (\ref{ch2divinvestorproblemcontrol}) and (\ref{ch2divinvestorproblemtime}) as jointly representing a stochastic differential game of control and stopping in which the controller modifies the system dynamics using impulse controls. Each of the investor's objectives is delegated to an individual player who plays so as to maximise their own objective whilst seeking an optimal response to the other player.

We develop the general underlying structure of the investment problem namely, the stochastic differential game of control and stopping. We, for the first time, solve the game providing a complete characterisation of the value of the game via a verification theorem.

\section{Current Literature}\label{section_lit_review}
Since its introduction to the literature, a considerable amount of work has been dedicated to the study of the lifetime ruin problem in addition to a number of variants of the problem. Variations of the original problem include models with stochastic consumption \cite{[3]}, stochastic volatility \cite{[4]}, ambiguity aversion \cite{[5]} amongst many other works. The probability of lifetime ruin model can be extended to address an analogous problem within the context of an investor who holds some portfolio of risky assets who seeks to both maximise their return whilst finding the optimal time to exit the market.

Typically, the lifetime ruin problem is analysed using a optimal stochastic control formulation in which the controller seeks both an optimal investment strategy (modelled using continuous controls) and an optimal time to sell all market holdings. Thus, in general the lifetime ruin problem in which the investor also seeks to maximise their returns can be formulated as an optimal stochastic control problem with discretionary stopping. In general, lifetime ruin problems in which the investor also seeks to maximise some performance criterion can be reformulated as stochastic differential games. The intuition behind this is that given a sufficient player aversion to lifetime ruin, nature can be viewed as a second player with the first player responding to nature's actions so as to avoid the occurrence of lifetime ruin. 

In \cite{[3]}, it is shown that the single investor portfolio problem in a Black-Scholes market in which an investor seeks to both maximise a running reward and minimise the probability of lifetime bankruptcy exhibits duality with controller-stopper games. Indeed, in \cite{[3]} it is shown that the value function of the investment problem is the convex dual of the value of a controller-stopper game. Similarly, in \cite{[15]} an investor portfolio problem with discretionary stopping is analysed by studying an optimal stopping-stochastic control differential game and proving an equivalence.

In \cite{[15]}, the value for a game in which the stopper seeks to minimise a convex risk measure defined over a common (zero-sum) payoff objective is characterised in terms of a Hamilton-Jacobi-Bellman Variational Inequality (HJBVI) to which it is proven that the value is a viscosity solution. The inclusion of a convex risk measure, as outlined in \cite{artzner1999coherent,follmer2002convex}, serves as a means by which risk attitudes of the investor are incorporated into the model. Moreover, the zero-sum payoff structure of the model implies that the strategies are appropriate for the extraction of optimal strategies in worst-case scenario analyses.

The problem of when capital injections should be performed (and when dividends should be paid) by the firm is an area of active research within theoretical actuarial science to which a great deal of attention has been focused. In general, the optimal capital injections and dividends model is represented as a single-player impulse control problem in which the controller seeks the optimal sequence of capital injections. In \cite{[9]} a model in which the firm can seek to raise capital (by issuing new equity) to be injected so as to allow the firm to remain solvent is considered. We refer the reader to \cite{[9]} and \cite{[10]} and references therein for exhaustive discussions.

\subsection*{Contributions}

This paper introduces a controller-stopper game in which the controller uses impulse controls; the results cover a general setting in which the underlying state process is a jump diffusion process. We extend existing game-theoretic impulse control results to now cover games in which i) the underlying state process is a jump-diffusion process and in contrast to \cite{[14]} ii) the payoff is no longer restricted to a zero-sum structure. 

To our knowledge, this paper is the first to deal with a jump-diffusion process within a stochastic differential game in which the players use impulse controls to modify the state process. Lastly, also to the best of our knowledge, this paper is the first to provide results pertaining to a non-zero-sum payoff structure within stochastic differential games for a controller-stopper game in which impulse controls used.

In the following section we describe the details of a general version of the controller-stopper game thereafter, we prove two key results: we firstly prove a set of verification theorems that characterise the conditions for a HJBI equation in non-zero-sum and zero-sum games. As in the Dynkin game case and controller-controller case, the HJBI equation is an obstacle problem in particular, the HJBI equation is an obstacle quasi-variational inequality. 

We begin by giving a canonical description of the game dynamics, starting with the zero sum game.

\section{Description of The Game}\label{section_game_description}
\subsection*{The Dynamics: Canonical Description}
The uncontrolled passive state evolves according to a stochastic process $X:[0,T]\times \Omega \to S\subset  \mathbb{R}^p,\; (p\in \mathbb{N})$, which is a jump-diffusion on $(\mathcal{C}([0,T]; \mathbb{R}^p ),(\mathcal{F}_{{(0,s)}_{s\in [0,T] }}, \mathcal{F}, \mathbb{P}_{0})$ that is, the state process obeys the following SDE:
\begin{align} 
dX_s^{t_0,x_{0}}=\mu( s,X_s^{t_0,x_{0}} )ds+\sigma(s,X_s^{t_0,x_{0}} )dB_s+\int  \gamma(X_{s^-}^{t_0,x_{0}},z) \tilde{N}(ds,dz)  ,\; X_{t_0}^{t_0,x_{0}}:= x_{0},&
\\\mathbb{P}-{\rm a.s.},\label{uncontrolledstateprocess}
\forall s\in [0,T],\; \forall (t_0,x_{0})\in [0,T]\times S,& 
\end{align}
where $B_s$ is an $m-$dimensional standard Brownian motion, $\tilde{N}(ds,dz)=N(ds,dz)-\nu(dz)ds$ is a compensated Poisson random measure where $N(ds,dz)$ is a jump measure and $\nu(\cdot):= \mathbb{E}[N(1,\cdot)]$ is a L\'{e}vy measure. Both $\tilde{N}$ and $B$ are supported by the filtered probability space and $\mathcal{F}$ is the filtration of the probability space $(\Omega ,\mathbb{P},\mathcal{F}=\{\mathcal{F}_s\}_{s\in [0,T] } )$. We assume that $N$ and $B$ are independent. 

We assume that the functions $\mu:[0,T]\times S\to S,\; \sigma:[0,T]\times S\to \mathbb{R} ^{p\times m}$ and $\gamma:\mathbb{R}^p\times \mathbb{R}^l\to \mathbb{R} ^{p\times l}$  are deterministic, measurable functions that are Lipschitz continuous and satisfy a (polynomial) growth condition so as to ensure the existence of (\ref{uncontrolledstateprocess}) \cite{ikeda2014stochastic}.

We note that the above specification of the filtration ensures stochastic integration and hence, the controlled jump-diffusion is well defined (this is proven in  \cite{stroock2007multidimensional}).

The generator of $X$ (of the uncontrolled process) acting on some function $\phi\in \mathcal{C}^{1,2} (\mathbb{R}^l,\mathbb{R}^p)$ is given by: 
\begin{equation} \mathcal{L}\phi(\cdot,x)=\sum_{i=1}^p   \mu_i (x)    \frac{\partial \phi}{\partial x_i}(\cdot,x)+\frac{1}{2} \sum_{i,j=1}^p  (\sigma \sigma^T )_{ij} (x)    \frac{ \partial^2\phi}{\partial x_i\partial x_{j} }+I\phi(\cdot,x)	
\label{generator}
\end{equation}
where $I$ is the integro-differential operator defined by:
\begin{equation} 
I\phi(\cdot,x):= \sum_{j=1}^{l} \int_{\mathbb{R}^p}  \{\phi(\cdot,x+ \gamma^{j}(x,z_j))-\phi(\cdot,x)-\nabla\phi(\cdot,x)  \gamma^{j} (x,z_{j} )\}  \nu_{j} (dz_{j}),\;{\forall  x\in \mathbb{R}^p.}\label{integro_differential_operator} 
\end{equation}

In this game there are two players, player I and player II.  Throughout the horizon of the game, each player incurs a cost which is a function of the value of the state process. Let the set  $\mathcal{T}$ be a given family of $\mathcal{F}-$measurable stopping times; at any point in the game $\rho \in \mathcal{T}$, player II can choose to terminate the game at which point the state process is stopped and both players receive a terminal cost. Player I influences the state process using impulse controls $u\in \mathcal{U}$ where $u(s)=\sum_{j\geq 1}\xi_j \cdot 1_{\{\tau_j\leq T \}}  (s)$ for all $0\leq t_0< s\leq T$ where $\xi_1,\xi_2,\ldots\in \mathcal{Z}\subset S$ are impulses that are executed at $\mathcal{F}$-measurable stopping times $\{\tau_i\}_{i\in\mathbb{N}}$ where $0\leq t_0\leq \tau_1< \tau_2< \dots <$ and the set $\mathcal{U}$ is a convex cone that defines the set of player I admissible controls. We assume that the impulses $\xi_j \in \mathcal{Z}$ are $\mathcal{F}-$measurable for all $j \in \mathbb{N}$. Hence, let us suppose that an impulse $\zeta \in \mathcal{Z}$ determined by some admissible policy $u\in\mathcal{U}$ is applied at some $\mathcal{F}-$measurable stopping time $\tau:\Omega \to [0,T]$ when the state is $x'=X^{t_0,x_0,\cdot} (\tau^-)$, then the state immediately jumps from $x'=X^{t_0,x_0,\cdot} (\tau^-)$ to $X^{t_0,x_0,u} (\tau)=\Gamma (x',\zeta)$ where $\Gamma :S\times \mathcal{Z}\to S$ is called the impulse response function. 

For notational convenience, as in \cite{chen2013impulse}, we use $u=[\tau_j,\xi_j ]_{j\geq 1} $ to denote the player I control policy $u=\sum_{j\geq 1}\xi_j  \cdot1_{\{\tau_j\leq T \}}  (s)\in \mathcal{U}$.

The evolution of the state process with actions is given by the following:
\begin{align}\nonumber
\hspace{-1.8 mm}X_r^{t_0,x_0,u}=x_0+\int_{t_0}^{r\wedge\rho}\hspace{-1.5 mm}\mu(s,X^{t_0,x_0,u}_s)ds+\int_{t_0}^{r\wedge\rho}\hspace{-1.5 mm}\sigma(s,X^{t_0,x_0,u}_s)dB_s
+\sum_{j\geq 1}\xi_j  \cdot 1_{\{\tau_j\leq r\wedge \rho \}}  (r)&\nonumber 
\\\nonumber+\int_{t_0}^{r}\int\gamma (X^{t_0,x_0,u} (s-),z) \tilde{N}(ds,dz)&
\\
\forall r\in [0,T];\;\forall (t_0,x_0)\in [0,T]\times S,\; \mathbb{P}-{\rm a.s.}&
\label{controlledstateprocessch1player}
\end{align}

Without loss of generality we assume that $X^{t_0,x_0,\cdot}_s=x_0$ for any $s\leq t_0$.

Player I has a cost function which is also the player II gain (or profit) function. The corresponding payoff function is given by the following expression which player I (resp., player II) minimises (resp., maximises):  
\begin{align}\nonumber 
J [t_0, x_0  ;u,\rho ]= \mathbb{E}\Bigg[\int_{t_0}^{\tau_s  \wedge\rho} f (s, X_s^{t_0, x_0  ,u } ) ds  &+ \sum_{m\geq 1}  c (\tau_m  , \xi _m  )  \cdot 1_{\{\tau_m  \leq  \tau_S  \wedge\rho \}}
\\&\begin{aligned}+G (\tau_S  \wedge\rho, X_{\tau_S  \wedge\rho }^{t_0, x_0  ,u } )1_{\{\tau_S  \wedge\rho<\infty\}}\Bigg],&\\ 
\forall (t_0,x_0)\in[0,T]\times S,& 
\end{aligned}
\end{align}
where $\tau_S:\Omega\to [0,T] $ is some random exit time, i.e. $\tau_S(\omega):=\inf\{s\in [0,T]|X_s^{t_0, x_0  ,\cdot }\in S\backslash A;\;\omega\in\Omega\}$ where $A$ is some measurable subset of $S$, at which point $\tau_S$ the game is terminated. The functions $f:[0,T]\times S\to \mathbb{R} , G:[0,T]\times S\to \mathbb{R}$ are the running cost function and the bequest function respectively and the function $c:[0,T]\times\mathcal{Z}\to\mathbb{R}$ is the intervention cost function. 

We assume that the function $G$ satisfies the condition $\underset{s\to\infty}{\lim}G(s,x)=0$ for any $x\in S$. Functions of the form $G(s,x)\equiv e^{-\delta s}\bar{G}(x)$ for some $\delta>0$ with $\bar{G}:|\bar{G}(x)|<\infty$ satisfy this condition among others.

Note that when $\mathcal{U}$ is a singleton the game is degenerate and collapses into a classical optimal stopping problem for player II with a value function and solution as that in ch.3 in \cite{oksendal2005applied}. Similarly, when $\mathcal{T}$ is a singleton the game collapses into a classical impulse control problem for player I with a value function and solution as that in ch.6 in \cite{oksendal2005applied}.

We now give some definitions which we will need to describe the system dynamics modified by impulse controls:

\begin{definition}\label{Definition 1.2.3.}
\

\noindent Denote by $\mathcal{T}_{(t,\tau')}$  the set of all $\mathcal{F}-$measurable stopping times in the interval $[t,\tau']$, where $\tau'$ is some stopping time s.th. $\tau' \leq T$. If $\tau'=T$ then we denote by $\mathcal{T}\equiv \mathcal{T}_{(0,T)}$. Let $u=[\tau_j,\xi_j]_{j\in\mathbb{N}}$ be a control policy where $\{\tau_j\}_{j\in\mathbb{N}}$ and $\{\xi_j\}_{j\in\mathbb{N}}$ are $\mathcal{F}_{\tau_j}-$ measurable stopping times and interventions respectively, then we denote by $ \mu_{[t,\tau]} (u)$ the number of impulses the controller (player I) executes within the interval $[t,\tau]$ under the control policy $u\in\mathcal{U}$ for some $\tau \in \mathcal{T}$.\end{definition} 

\begin{definition}\label{Definition 1.2.4.}
Let $u\in\mathcal{U}$ be a player I impulse control. We say that an impulse control is admissible for player I on $[0,T]$ if either the number of impulse interventions is finite on average that is to say we have that:
\begin{equation}\mathbb{E}[\mu_{[0,T]} (u)]<\infty	 \end{equation}
or if $\mu_{[0,T]} (u)=\infty\implies\lim_{j\to \infty }\tau_j=\infty$.
\end{definition}
We shall hereon use the symbol $\mathcal{U}$ to denote the set of admissible player I controls. For controls $u\in \mathcal{U}$ and $u'\in \mathcal{U}$, we interpret the notion $u\equiv u'$ on $[0,T]$ iff $\mathbb{P}(u=u'$ a.e. on $[0,T])=1$. 

\begin{definition}\label{Definition 1.2.5.}
\

\noindent Let $u(s)=\sum_{j\geq 1}\xi_j  \cdot 1_{\{\tau_j\leq T \}}  (s) \in \mathcal{U}$ be a player I impulse control defined over $[0,T]$, further suppose that $\tau:\Omega \to [0,T]$ and $\tau':\Omega \to [0,T]$ are two $\mathcal{F}-$measurable stopping times with $\tau\geq s>\tau'$, then we define the restriction $u_{[\tau',\tau ]} \in \mathcal{U}$ of the impulse control $u(s)$ to be $u(s)=\sum_{j\geq 1}\xi_{\mu_{]t_0,\tau)}(u)+j}  \cdot 1_{\{\tau_{\mu_{]t_0,\tau)}(u)+j} \geq s\geq \tau' \}}  (s)$. 
\end{definition}

Markov controls are those in which the player uses only information about the current state and duration of the game rather than explicitly incorporating information about the other player's decisions or utilising information on the history of the game. In general, in two player games the player who performs an action first employs the use of a \textit{strategy}, that is, a protocol that involves different responses according to different future actions that their opponent performs. This is in contrast to performing a fixed action or using controls of a Markovian type described above. In general, the use of strategies affords the acting player the ability to increase its rewards since their action is now a function of the other player's latter decisions.

In light of the above remarks, limiting the analysis to Markov controls may incur too strong of a restriction on the abilities of the players to perform optimally. It is however, well known that under mild conditions, for games of the type discussed in this paper that involve a diffusive state process, using Markov controls gives as good performance as an arbitrary $\mathcal{F}$-adapted control (see for example Theorem 11.2.3 in \cite{oksendalapplied2007}). Consequently, in the following analysis we restrict ourselves to Markov controls and hence for player I, the control policy takes the form $u=u(s,\omega)\in\mathcal{U};\; s\in [0,T],\omega\in\Omega$. The control policy $u=u(s,\omega)\in\mathcal{U}$ is hence a stochastic process that modifies the state process directly. Moreover, the player I control can be written in the form $u=\tilde{f}_1(s,X_s)$ for any $s\in[0,T]$ where $\tilde{f}_1:[0,T]\times S\to U$ and $U\subset \mathbb{R}^p$ and $\tilde{f}_1$ is some measurable map  w.r.t. $\mathcal{F}$.

The following definition is a key object in the analysis of impulse control models:
\begin{definition}\label{Definition 2.1.1.}
\

\noindent Let $\tau\in\mathcal{T}$, we define the [non-local] intervention operator $\mathcal{M}:\mathcal{H}\to \mathcal{H}$ acting at a state $X(\tau)$ by the following expression:
\begin{equation}
\mathcal{M} \phi(\tau,X(\tau)):=\inf_{z\in \mathcal{Z}}[\phi(\tau,\Gamma (X(\tau^-),z))+c(\tau, z)\cdot 1_{\{\tau\leq T \}}  ],\label{intervention_operator_definition_equation} \end{equation}
for some function $\phi :[0,T]\times S\to \mathbb{R}$ and $\Gamma : S\times \mathcal{Z}\to  S$ is the impulse response function.
\end{definition}

Of particular interest is the case when the intervention operator is applied to the value function $\mathcal{M}V(\cdot,x)$ --- a quantity which represents the value of a strategy when the controller performs an optimal intervention then behaves optimally thereafter given an immediate optimal intervention taken at a state $x\in S$. The intuition behind (\ref{intervention_operator_definition_equation}) is as follows: suppose at time $\tau^-$ the system is at a state $X(\tau^-)$ and an intervention $z\in\mathcal{Z}$ is applied to the process, then a cost of $c(\tau,z)$ is incurred and the state then jumps from $X(\tau^-)$ to $\Gamma(X(\tau^-),z)$. If the controller acts optimally thereafter, the cost of this strategy, starting at state $\Gamma(X(\tau^-),z)$ is $V(\tau,\Gamma(X(\tau^-),z)+c(\tau,z)$. Lastly, choosing the action that minimises costs leads to $\mathcal{M}V$.

\begin{remark}\label{remark 1.2.7.}
\

\noindent We note that whenever it is optimal for the controller to intervene, $\mathcal{M}V=V$ since the value function describes the player payoff under optimal behaviour. However, at any given point an immediate intervention may not be optimal, hence the following inequality holds pointwise:
\begin{align}
\mathcal{M}V(s,x)&\geq V(s,x), 	\label{intervention_operator_inequality_one_player} \qquad \hfill \forall  (s,x)\in [0,T]\times S.     
\end{align} 
\end{remark}

The results contained in this paper are built exclusively under the following set of assumptions unless otherwise stated:\\
\textbf{Standing Assumptions}\

\noindent A.1.1. Lipschitz Continuity\

We assume there exist real-valued constants $c_{\mu},c_{\sigma}>0$ and $c_{\gamma} (\cdot)\in L^1\cap L^2 ( \mathbb{R}^l,\nu)$ s.th. $\forall  s\in [0,T], \forall  x,y\in S$ and $\forall  z\in \mathbb{R}^l$ we have:
\begin{align*}
|\mu(s,x)-\mu(s,y)|&\leq c_{\mu} |x-y|\\
|\sigma(s,x)-\sigma(s,y)|&\leq c_{\sigma} |x-y|
\\\int_{|z|\geq 1} |\gamma(x,z)-\gamma(y,z)| &\leq c_{\gamma} (z)|x-y|.
\end{align*}
A.1.2. Lipschitz Continuity
\

We also assume the Lipschitzianity of the functions $f$ and $G$ that is, we assume the existence of real-valued constants $c_f,c_G>0$ s.th.  $\forall  s\in [0,T], \forall  (x,y) \in S$ we have for $R\in \{f,G\}$:
\begin{equation*}
|R(s,x)+R(s,y)|\leq c_R |x-y|.
\end{equation*}
A.2. Growth Conditions
\

We assume the existence of a real-valued constants $d_{\mu},d_{\sigma}>0$ and $d_{\gamma} (\cdot)\in L^1\cap L^2 ( \mathbb{R}^l,\nu), \rho\in [0,1)$ s.th. $\forall  (s,x)\in [0,T]\times S$ and $\forall  z\in \mathbb{R}^l$  we have:
\begin{align*}
|\mu(s,x)|\leq d_{\mu} (1+|x|^\rho)
\\
|\sigma(s,x)|\leq d_{\sigma} (1+|x|^\rho)
\\
\int_{|z|\geq 1} |\gamma(x,z)| \leq d_{\gamma} (1+|x|^\rho).
\end{align*}
We also make the following assumptions on the cost function $c:[0,T]\times S\to \mathbb{R}$:\\
A.3.\

Let $\tau,\tau' \in \mathcal{T}$ be $\mathcal{F}-$measurable stopping times s.th. $0\leq \tau<\tau'\leq T$ and let $\xi,\xi' \in \mathcal{Z}$ be measurable impulse interventions. Then we assume that the following statements hold:
	\begin{align} 
	c(\tau,\xi+\xi' )&\leq c(\tau,\xi)+c(\tau,\xi' ),\\
	  c(\tau,\xi)&\geq c(\tau',\xi).
\end{align}
A.4.\ \label{AA.4.}

We also assume that there exists a constant $\lambda_c>0$ s.th.  $\inf_{\xi \in \mathcal{Z})}c(s,\xi)\geq \lambda_c\forall  s \in [0,T]$ where $\xi \in \mathcal{Z}$ is a measurable impulse intervention.\

Assumptions A.1.1 and A.2 ensure the existence and uniqueness of a solution to (\ref{uncontrolledstateprocess}) (c.f. \cite{[17]}). Assumption A.3 (i) (subadditivity) is required in the proof of the uniqueness of the value function. Assumption A.3 (ii) (the player cost function is a decreasing function in time) and may be interpreted as a discounting effect on the cost of interventions. Assumption A.1.2 is required to prove the regularity of the value function (see for example \cite{[11]} and for the single-player case, see for example \cite{[18]}). Assumption A.3 (ii) was introduced (for the two-player case) in \cite{[19]} though is common in the treatment of single-player case problems (e.g. \cite{[18],[15]}). Assumption A.4 is integral to the definition of the impulse control problem.
\

Throughout the script we adopt the following standard notation (e.g. \cite{[12],[15],[16]}):

\subsection*{Notation}

Let $\Omega$  be a bounded open set on  $\mathbb{R}^{p+1}$. Then we denote by:
$\bar{\Omega}$  --- The closure of the set $\Omega$.\\
$Q(s,x;R)={{(s',x' ) \in \mathbb{R} ^{p+1}:\max |s'-s|^{\frac{1}{2}}  ,|x'-x|  }<R,s'<s}$. \\
$\partial \Omega$  --- The parabolic boundary $\Omega$  i.e. the set of points $(s,x) \in \bar{\mathcal{S}}$ s.th. $R>0, Q(s,x;R)\not\subset\bar{\Omega}$.\\
$\mathcal{C}^{\{1,2\}} ([0, T],\Omega )=\{h \in C^{\{1,2\}} (\Omega ): \partial_s h, \partial_{x_i,x_{j} } h \in C(\Omega )\}$, where  $\partial_s$ and  $\partial_{x_{i}, x_{j}}$ denote the temporal differential operator and second spatial differential operator respectively.\\
$\nabla\phi=(\frac{\partial \phi}{\partial x_1 },\ldots,\frac{\partial \phi}{\partial x_p})$ --- The gradient operator acting on some function $\phi \in C^1 ([0,T]\times \mathbb{R}^p)$.\\
$|\cdot|$   --- The Euclidean norm to which $\langle x,y \rangle$  is the associated scalar product acting between two vectors belonging to some finite dimensional space.

\section{Statement of Main Results}\label{section_main_results}
We prove two key results --- in order to characterise the value of the game, we prove a \textit{verification theorem} for the game; this then allows us to characterise equilibrium conditions in terms of some non-linear PDE. As we show, this approach produces a practical scheme by which the solution to the optimal liquidity and lifetime ruin problem is computed. We firstly prove the theorem for the zero-sum game case, we then extend the analysis to stochastic differential games of control and stopping with non zero-sum payoff structures.

In particular, for the non zero-sum case, we have the following result (Theorem \ref{Verification theorem for Non-Zero-Sum Stochastic Differential Games of Control and Stopping with Impulse Controls}) which we prove:

Let $X: [0,T]\times\Omega \to S$ be a stochastic process that evolves according to \eqref{uncontrolledstateprocess} and where $S\subset\mathbb{R}^p$ is some \textit{solvency region}. Let $\phi_i$ be smooth test functions so that we can take first order temporal derivatives and second order spatial derivatives within the interior of $S$ for $i\in\{1,2\}$; then if $\phi_i$ satisfy the following quasi-variational inequalities: 
\begin{align}
 \begin{cases}\max\{\partial_s \phi_i(t,x)+\mathcal{L}\phi_i(t,x)+f_i(t,x),\phi_i(t,x)-\mathcal{M}_i \phi_i(t,x)\}&=0  \\
  \phi_i(\cdot,y)=G_i(\cdot,y) \qquad\qquad \forall y \in  S, \forall (t,x)\in [0,T]\times S,
  \end{cases}
  \end{align}
where $f_i$ and $G_i$ are the player $i$ running cost functions and terminal cost functions respectively and where $\mathcal{L}$ is the stochastic generator of the diffusion process $X$ (c.f. (\ref{generator})), then $\phi_i$ is the player $i$ value function for the non zero-sum game.  

Having proven these results, we then apply the analysis to prove the following set of results relating to the optimal liquidity control and lifetime ruin investment problem stated in Section \ref{section_investment_problem}.

To our knowledge, this paper is the first to deal with a jump-diffusion process within a stochastic differential game in which the players use impulse controls to modify the state process. Additionally, to our knowledge, this is the first game that involves impulse controls in which the role of one of the players is to stop the game at a desirable point.

\section{Stochastic Differential Games of Impulse Control and Stopping}\label{section_zero_sum_game}

In this section, we analyse the zero-sum game and prove a verification theorem (Theorem \ref{Verification_theorem_for_Zero-Sum_Stochastic_Differential_Games_of_Control_and_Stopping_with_Impulse_Controls}). The theorem provides the conditions under which, if a sufficiently smooth solution to a HJBI equation can be found then we have a candidate for the value function of the game. The following verification theorem additionally characterises the conditions under which the value of the game satisfies a HJBI equation and characterises the equilibrium controls for the game. Later, we use the conditions of Theorem \ref{Verification_theorem_for_Zero-Sum_Stochastic_Differential_Games_of_Control_and_Stopping_with_Impulse_Controls} to derive the optimal investment strategy for the optimal liquidity control and lifetime ruin model.

\begin{theorem}[Verification theorem for Zero-Sum Games of Control and Stopping]\label{Verification_theorem_for_Zero-Sum_Stochastic_Differential_Games_of_Control_and_Stopping_with_Impulse_Controls}
\

\noindent Suppose the problem is to find $\phi\in\mathcal{H}$ and $ (\hat{u},\hat{\rho} ) \in \mathcal{U}\times \mathcal{T} $ s.th:
\begin{align}
\phi(t,x)=\sup_{\rho \in \mathcal{T}}\left( \inf_{u \in \mathcal{U}}J^{(u,\rho)} [t,x] \right)= \inf_{u \in \mathcal{U}} \left(\sup_{\rho \in \mathcal{T}}J^{(u,\rho)} [t,x]\right) =J^{(\hat{u},\hat{\rho})} [t,x],
\qquad \forall (t,x) \in [0,T]\times S,
\end{align}
where if $(\hat{u},\hat{\rho}) \in \mathcal{U}\times \mathcal{T} $ exists, it is an optimal pair consisting of the optimal control for player I and the optimal stopping time for player II (resp.).

Let $\tau$ be some $\mathcal{F}-$measurable stopping time and denote by $\hat{X}(\tau)=X(\tau^{-} )+\Delta_N X(\tau) $, where $\Delta_N X(\tau) $ denotes a jump at time $\tau$ due to $\tilde{N}$. Suppose that the value of the game exists. 

Suppose also that there exists a function $\phi \in \mathcal{C}^{1,2} ([0,T],S)\cap\mathcal{C}([0,T],\bar{S})$ that satisfies technical conditions (T1) - (T4) and the following conditions:
\renewcommand{\theenumi}{\roman{enumi}}
 \begin{enumerate}[leftmargin= 6 mm]
	\item \label{verif_theorem_c-s_zero_item_ii_intervention_stopping_inequalities} $\phi\leq \mathcal{M}\phi$ on $S$ and $\phi\geq G$ on $S$ and the regions  $D_1$ and $D_2$ are defined by:\\
$D_1=\{x \in S;\phi(\cdot,x)<\mathcal{M}\phi(\cdot,x)\}$ and  $D_2=\{x \in S;\phi(\cdot,x)>G(\cdot,x)\}$\\
where we refer to $D_1$ (resp., $D_2$) as the player I (resp., player II) continuation region.

	\item \label{verif_theorem_c-s_zero_item_iii_HJBI_stopping_inequality}$ \frac{\partial \phi}{\partial s}+\mathcal{L}\phi(s,X^{\cdot,u } (s))+f(s,X^{\cdot,u } (s))\geq 0,\hspace{1 mm}   \forall  u \in \mathcal{U}$ on $S\backslash{\partial D_1}.$
	
	\item \label{verif_theorem_c-s_zero_item_iv_HJBI_equation}$\frac{\partial \phi}{\partial s}+\mathcal{L}\phi(s,X^{\cdot,\hat{u}} (s))+f(s,X^{\cdot,\hat{u}} (s))=0 \hspace{1 mm}$ in $D_1\cap D_2. $\ \label{verification_theorem_controller_stopper_game_pde}

	\item \label{verif_theorem_c-s_zero_item_v_stopping_time_criterion} For $u \in \mathcal{U}$, define $\rho_D=\rho_D^u= \inf\{s>t_0, X^{\cdot,u} (s)\notin D_2\}$ and specifically, $\hat{\rho}_D=\hat{\rho}= \inf\{s>t_0, X^{\cdot,u} (s)\notin D_2\}$.

	\item \label{verif_theorem_c-s_zero_item_vi_continuity_at_termination} {$X^{\cdot,u} (\tau_S ) \in \partial S, \mathbb{P}-{\rm a.s}$}. for ${\tau_S< \infty}$ and $\phi(s,X^{\cdot,u} (s))\to G(\tau_S,X^{\cdot,u} (\tau_S\wedge\rho)) $  as $s\to \tau_S^-\wedge\rho^-  $, $\mathbb{P}-{\rm a.s}.,\forall  x \in S,\; \forall u \in \mathcal{U}$.
\end{enumerate}

Put $\hat{\tau}_0\equiv t_0$ and define $\hat{u}:=[\hat{\tau}_j,\hat{\xi}_j ]_{j \in \mathbb{N}}$ inductively by:\\ $\hat{\tau}_{j+1}= \inf\{s>\tau_j;X^{\cdot,\hat{u}_{[t_0,s]}} (s)\notin D_1 \}\wedge\tau_S\wedge\rho$,
then $ (\hat{u},\hat{\rho} ) \in \mathcal{U}\times \mathcal{T}$ are an optimal pair for the game, that is to say that we have:\begin{equation}
\phi(t,x)= \inf_{u \in \mathcal{U}}\left(  \sup_{\rho \in \mathcal{T}}  J^{(u,\rho)} [t,x]\right)=  \sup_{\rho \in \mathcal{T}}    \left(\inf_{u \in \mathcal{U}}   J^{(u,\rho)} [t,x]\right)=  J^{(\hat{u},\hat{\rho})} [t,x],	
\qquad\forall (t,x) \in [0,T]\times S.
\end{equation}
\end{theorem}

Theorem \ref{Verification_theorem_for_Zero-Sum_Stochastic_Differential_Games_of_Control_and_Stopping_with_Impulse_Controls} provides a characterisation of the value of the game in terms of a dynamic programming equation (which is simply the non-linear PDE in (\ref{verification_theorem_controller_stopper_game_pde})). In particular, Theorem \ref{Verification_theorem_for_Zero-Sum_Stochastic_Differential_Games_of_Control_and_Stopping_with_Impulse_Controls} says that given some solution to the non-linear PDE in (\ref{verification_theorem_controller_stopper_game_pde}), then this solution coincides with the value of the game from which we can calculate the optimal controls for each player. Moreover, Theorem \ref{Verification_theorem_for_Zero-Sum_Stochastic_Differential_Games_of_Control_and_Stopping_with_Impulse_Controls} provides insight as to the structure of the solution, in particular, provided player II has not terminated the game, player I exercises an impulse control whenever the state process exits the continuation region $D_1$. Similarly, player II does nothing when the state process remains within the region $D_2$ and terminates the game at the first hitting time on $S\backslash{\partial D_2}$.

Before stating the proof of Theorem \ref{Verification_theorem_for_Zero-Sum_Stochastic_Differential_Games_of_Control_and_Stopping_with_Impulse_Controls}, we make the following set of remarks which also applies to Theorem \ref{Verification theorem for Non-Zero-Sum Stochastic Differential Games of Control and Stopping with Impulse Controls}.:

\begin{remark}\label{Remark 5.2.} 
\

\noindent For the jump-diffusion process considered here,  we can automatically conclude that \\$\hat{\xi}_k \in  \arg\hspace{-0.55 mm}\min _{z \in \mathcal{Z}}{\phi(\Gamma(x,z))
+c(\tau_k,z)},\; \forall  k \in \mathbb{N},\; x \in S$ where $\tau_k\in\mathcal{T}$ is an $\mathcal{F}-$measurable stopping time that exists.
\end{remark}

The proof of the lemma is straightforward since we need only prove that the  infimum is in fact a minimum. This follows directly from the fact that the cost function is minimally bounded (c.f.  A.\ref{AA.4.}) and that the boundedness of the value functions.

To prove Theorem \ref{Verification_theorem_for_Zero-Sum_Stochastic_Differential_Games_of_Control_and_Stopping_with_Impulse_Controls}, we firstly require the following results. 

The first result enables us to perform limiting procedures close to the boundary of the player II continuation region:
\begin{theorem}[(Approximation Theorem) (Theorem 3.1 in \cite{oksendalapplied2007})]\label{Approximation_Theorem}
\

\noindent Let $\hat{D} \subset S$ be an open set and let us assume that $X(\tau_S ) \in \partial S$ and suppose that $\partial \hat{D}$ is a Lipschitz surface. Let $\psi: \bar{S}\to \mathbb{R}$ be a function s.th. $\psi \in \mathcal{C}^1 (S)\cap\mathcal{C}(\bar{S})$ and $\psi \in \mathcal{C}^2 (S\backslash\partial \hat{D})$ and suppose the second order derivatives of $\psi$ are locally bounded near $\partial \hat{D}$; then there exists a sequence of functions $\{\psi_m\}_{m=1}^{ \infty} \in \mathcal{C}^2 (S)\cap\mathcal{C}(\bar{S})$ s.th.\begin{align*}
	\qquad \qquad &\lim_{m\to  \infty}\psi_m \to \psi\text{ pointwise dominatedly in }\bar{S}.
	\\&
	\lim_{m\to  \infty}\frac{\partial \psi_m}{\partial x_i}\to \frac{\partial \psi}{\partial x_i}\text{ pointwise-dominatedly in }S.
	\\&
	\lim_{m\to  \infty}\frac{\partial^2 \psi_m}{\partial x_i \partial x_j }\to \frac{\partial^2 \psi}{\partial x_i \partial x_j}\text{ and } \lim_{m\to  \infty}\mathcal{L}\psi_m\to \mathcal{L}\psi\text{ pointwise dominatedly in }S\backslash\partial \hat{D}.\end{align*}
\end{theorem}

We give a statement of the following result without proof:

\begin{lemma}[Lemma 3.10 in \cite{[18]}]
\label{Lemma 4.9.}
\

\noindent The non-local intervention operator $\mathcal{M}$ is continuous wherein we can deduce the existence of a constants $c_1,c_2>0$ s.th. $\forall  x,y \in S$ and $s,t \in [0,T] $: 
\begin{align*}
\text{  i)}& &|\mathcal{M}V (s,x)-\mathcal{M}V (s,y)|&\leq c_1 |x-y|,\\\indent	\text{ii)}&  &|\mathcal{M}V (t,x)-\mathcal{M}V (s,x)|&\leq c_2 |t-s|^{\frac{1}{2}}.	
\end{align*}
\end{lemma}
A proof of the result is reported in \cite{zhang2011stochastic}.

\begin{lemma}\label{Lemma 2.1.6.}
\

\noindent Let $V \in \mathcal{H}$ be a bounded function and $ (\tau,x) \in [0,T]\times S$ where $\tau$ is some $\mathcal{F}-$ measurable stopping time, then the set $\Xi(\tau,x) $ defined by: 
\begin{equation}
\Xi(\tau,x):=\left\{\xi \in \mathcal{Z}:\mathcal{M}V(\tau^-,x)=V(\tau,x+\xi)+c(\tau, \xi) \cdot 1_{\{\tau\leq T\}}  \right\}
\end{equation}
is non-empty.
\end{lemma}

The proof of the lemma is straightforward since we need only prove that the infimum is in fact a minimum. This follows directly from the fact that the cost function is minimally bounded (c.f.  A.\ref{AA.4.}) and that the value functions are also bounded.

We are now in a position to prove the theorem; some ideas for the proof come from \cite{bayraktar2011minimizing,ikeda2014stochastic}:
\begin{refproof}[Proof of Theorem \ref{Verification_theorem_for_Zero-Sum_Stochastic_Differential_Games_of_Control_and_Stopping_with_Impulse_Controls}]
\ 

\noindent In the following, for ease of exposition we introduce the following notation:
\begin{align}
&Y^{y_0,\cdot}(s)\equiv (s,X^{t_0,x_0,\cdot}(t_0+s)), \quad y_0\equiv (t_0,x_0), \; \forall s\in [0,T-t_0],
\\& \hat{Y}^{y_0,\cdot}(\tau)=Y^{y_0,\cdot}(\tau^{-} )+\Delta_N Y^{y_0,\cdot}(\tau), \quad \tau\in\mathcal{T},
\end{align}
where $\Delta_N Y(\tau) $ denotes a jump at time $\tau$ due to $\tilde{N}$.

Correspondingly, we adopt the following impulse response function $\hat{\Gamma}: \mathcal{T}\times S\times \mathcal{Z}\to  \mathcal{T}\times S$ acting on  $y'\equiv (\tau,x')\in \mathcal{T}\times S$ where $x'\equiv X^{t_0,x_0,\cdot}(t_0+\tau^-)$ and $\hat{\Gamma}$ is given by: 
\begin{align}
\hat{\Gamma}(y',\zeta)\equiv (\tau,\Gamma (x',\zeta))=(\tau,X^{t_0,x_0,\cdot} (\tau)),\quad \forall \xi\in\mathcal{Z},\; \forall\tau\in\mathcal{T} .
\end{align}
We begin by fixing the player I control $\hat{u} \in \mathcal{U}$ and let us define $\rho_m:=\rho\wedge m;m=1,2\ldots$. By Dynkin's formula for jump-diffusion processes (see for example Theorem 1.24 in \cite{karatzas2000utility}) we have: \begin{equation}
\mathbb{E}[\phi(Y^{y_0,\hat{u}} (\hat{\tau}_j ))]-\mathbb{E}[\phi(\hat{Y}^{y_0,\hat{u}} (\hat{\tau}_{j+1} ))]=-\mathbb{E}\left[\int_{\hat{\tau}_j}^{\hat{\tau}_{j+1}} \frac{\partial \phi}{\partial s}+\mathcal{L}[\phi(Y^{y_0,\hat{u}} (s)]ds\right]. 	\label{zerosumproofdynkin1}\end{equation}
Summing (\ref{zerosumproofdynkin1}) from $j=0$ to $j=k$ for some $0<k<\mu_{[t_0,\rho_m ]} (\hat{u} )-1$ (recall the definition of $\mu_{[t_0,s]} (u)$ from Definition \ref{Definition 1.2.3.}) and observing that by (\ref{verif_theorem_c-s_zero_item_iv_HJBI_equation}) we have that -$ (\partial_s+\mathcal{L})\phi=f$, we find that: \begin{align}\nonumber
&\phi(y_0)+\sum_{j=1}^{k} \mathbb{E}[\phi(Y^{y_0,\hat{u}} (\hat{\tau}_j ))-\phi(\hat{Y}^{y_0,\hat{u}} (\hat{\tau}_j ))] -\mathbb{E}[\phi(\hat{Y}^{y_0,\hat{u}} (\hat{\tau}_{k+1} ))] \nonumber	
\\&
=-\mathbb{E}\left[\int_{t_0}^{\hat{\tau}_{k+1}} \left(\frac{\partial \phi}{\partial s}+\mathcal{L}\phi(Y^{y_0,\hat{u}_{[t_0,s]}} (s))\right)ds\right]=\mathbb{E}\left[\int_{t_0}^{\hat{\tau}_{k+1}} f(Y^{y_0,\hat{u}_{[t_0,s]}} (s))ds\right]\label{verifproofstep1}. 	
\end{align}
Now by definition of the non-local intervention operator $\mathcal{M}$ and by choice of $\hat{\xi}_j \in \mathcal{Z}$, we have that:
\begin{equation}
\phi(Y^{y_0,\hat{u}} (\hat{\tau}_j ))=\phi(\hat{\Gamma}(\hat{Y}^{y_0,\hat{u}} (\hat{\tau}_j^- ),\hat{\xi}_j ))=\mathcal{M}\phi(\hat{Y}^{y_0,\hat{u}} (\hat{\tau}_j^- ))-c(\hat{\tau}_j,\hat{\xi}_j )\cdot1_{\{\hat{\tau}_j\leq T\}}\label{ch3.verifproofs2},	
\end{equation}
(using the fact that $\inf_{z\in \mathcal{Z}}[\phi(\tau',\Gamma (X(\tau'^-),z))+c(\tau', z)\cdot 1_{\{\tau'\leq T \}}  ]=0$ whenever $\tau'>\tau_S\wedge\rho$). Hence after deducting $\phi(\hat{Y}^{y_0,\hat{u}} (\hat{\tau}_j^- ))$ from both sides we find: 
\begin{align}
\mathcal{M}\phi(\hat{Y}^{y_0,\hat{u}} (\hat{\tau}_j^- ))-\phi(\hat{Y}^{y_0,\hat{u}} (\hat{\tau}_j^- ))-c(\hat{\tau}_j,\hat{\xi}_j )=\phi(Y^{y_0,\hat{u}} (\hat{\tau}_j ))-\phi(\hat{Y}^{y_0,\hat{u}} (\hat{\tau}_j^- )),\label{interventioneq1}	
\end{align}
and by (\ref{verif_theorem_c-s_zero_item_vi_continuity_at_termination}) we readily observe that: $ \phi(Y^{y_0,\hat{u}} (\tau_s ))-\phi(\hat{Y}^{y_0,\hat{u}} (\tau_s ))=0$, hence after plugging (\ref{interventioneq1}) into (\ref{verifproofstep1}) we obtain the following: 
\begin{align}
&\phi(y_0)+\sum_{j=1}^{k} \mathbb{E}[\mathcal{M}\phi(\hat{Y}^{y_0,\hat{u}} (\hat{\tau}_j^- ))-\phi(\hat{Y}^{y_0,\hat{u}} (\hat{\tau}_j^- ))] -\mathbb{E}[\phi(\hat{Y}^{y_0,\hat{u}} (\hat{\tau}_{k+1} ))] 	
\nonumber\\&
=\mathbb{E}\left[\int_{t_0}^{\hat{\tau}_{k+1}} f(Y^{y_0,\hat{u}_{[t_0,s]}} (s))ds+\sum_{j=1}^{k} c(\hat{\tau}_j,\hat{\xi}_j )  \cdot 1_{\{\hat{\tau}_j\leq \tau_S \}}  \right]. \label{zerosumaftersumintervent}	
\end{align}
Note that our choice of $\hat{\xi}_k \in \mathcal{Z}$ induces equality in (\ref{zerosumaftersumintervent}).

Since the number of interventions in (\ref{zerosumaftersumintervent}) is bounded above by $\mu_{[t_0,\rho_m\wedge\tau_S]} (\hat{u})\wedge m$ for some $m< \infty$ and (\ref{zerosumaftersumintervent}) holds for any $k \in \mathbb{N}$, taking the limit as $k\to  \infty$ in (\ref{zerosumaftersumintervent}) gives: 
\begin{align}
&\qquad\qquad\phi(y_0)+\sum_{j=1}^{\mu_{[t_0,\rho_m\vee \tau_S ]} (\hat{u})}\hspace{-3 mm} \mathbb{E}[\mathcal{M}\phi(\hat{Y}^{y_0,\hat{u}} (\hat{\tau}_j^- ))-\phi(\hat{Y}^{y_0,\hat{u}} (\hat{\tau}_j^- ))]  \nonumber	
\\&
= \mathbb{E}\left[\phi (\hat{Y}^{y_0,\hat{u}} (\rho_m  \wedge \tau_S  ))+ \int_{t_0}^{\rho_m  \wedge \tau_S}f(Y^{y_0, \hat{u}_{[t_0,s]}}  )(s ))ds  +\sum_{j=1}^{\mu_{[t_0, \rho_m  \vee  \tau_S  ] } ( \hat{u})}c (\hat{\tau}_j  , \hat{\xi}_j  )  \cdot 1_{\{\hat{\tau}_j  \leq  \rho_m  \wedge \tau_S  \}}   \right]. \label{zerosumafterklim}	
\end{align}

Now $\lim_{m\to  \infty}\sum_{j=1}^{\mu_{[t_0,\rho_m\vee \tau_S ]} (\hat{u})} \mathbb{E}[\mathcal{M}\phi(\hat{Y}^{y_0,\hat{u}} (\hat{\tau}_j^- ))-\phi(\hat{Y}^{y_0,\hat{u}} (\hat{\tau}_j^- ))] =0$ since also by (\ref{verif_theorem_c-s_zero_item_vi_continuity_at_termination}) we have that $\phi(\hat{Y}^{y_0,\cdot} (\hat{\tau}_j ))-\phi(\hat{Y}^{y_0,\cdot} (\hat{\tau}_j^- ))=0$, $\mathbb{P}-$a.s. when $\hat{\tau}_j=\tau_S$. Similarly, we have by (\ref{verif_theorem_c-s_zero_item_vi_continuity_at_termination}) that $\phi(Y^{y_0,\hat{u}} (s))\to G(Y^{y_0,\hat{u}} (\tau_S\wedge\rho)) $  as $s\to \tau_S^-\wedge\rho^-,\;  \mathbb{P}-$a.s.

Now since $\rho_m\wedge\tau_S\to \rho\wedge\tau_s$ as $m\to  \infty$, we can exploit the quasi-left continuity of $X$ (for further details see \cite{protter2005stochastic} (Proposition I.2.26 and Proposition I.3.27)) and the continuity properties of $f$, we find that there exists some $c>0$ s.th. 
\begin{align*}
&\qquad\qquad\qquad\qquad\qquad\left|\lim_{m\to  \infty}\phi(\hat{Y}^{y_0,u} (\rho_m\wedge\tau_s ))+\lim_{m\to  \infty} \int_{t_0}^{\rho_m\wedge\tau_s} f(Y^{y_0,\hat{u}} (s))ds\right| 
    \\&
\leq c \lim_{m\to  \infty} \left(1+\left|\hat{Y}^{y_0,\hat{u}} (\rho_m\wedge\tau_s )\right|+\int_{t_0}^{\rho_m\wedge\tau_s} \left|Y^{y_0,\hat{u}} (s)\right|ds\right) 
\leq c(1+\tau_S )(1+\sup_{s \in [0,T]}|X^{t_0,x_0,\hat{u}} (s)| ) \in \mathbb{L}.
\end{align*}
Hence, taking the limit as $m\to\infty$ and using the Fat\^ou lemma and (\ref{zerosumafterklim}), we find that:\begin{align*}
\phi(y_0)&=\sum_{j=1}^{\mu_{[t_0,\rho_m\vee \tau_S ]} (\hat{u})} \mathbb{E}\Bigg[\mathcal{M}\phi(\hat{Y}^{y_0,\hat{u}} (\hat{\tau}_j^- ))-\phi(\hat{Y}^{y_0,\hat{u}} (\hat{\tau}_j^- ))] \\&\qquad+\mathbb{E}[\phi(\hat{Y}^{y_0,\hat{u}} (\rho_m\wedge\tau_S ))]+\int_{t_0}^{\rho_m\wedge\tau_s} f(Y^{y_0,\hat{u}} (s))ds+\sum_{j\geq  1} c(\hat{\tau}_j,\hat{\xi}_j )  \cdot 1_{\{\hat{\tau}_j\leq {\rho_m\wedge\tau_S }\
}  \Bigg]
\\&=\lim_{m\to  \infty}   \inf\mathbb{E}\Bigg[\sum_{j=1}^{\mu_{[t_0,\rho_m\vee \tau_S ]} (\hat{u})} \mathbb{E}[\mathcal{M}\phi(\hat{Y}^{y_0,\hat{u}} (\hat{\tau}_j^- ))-\phi(\hat{Y}^{y_0,\hat{u}} (\tau_j^- ))] 
\\&\qquad+\phi(\hat{Y}^{y_0,\hat{u}} (\rho_m\wedge\tau_S ))+\int_{t_0}^{\rho_m\wedge\tau_s} f(Y^{y_0,\hat{u}} (s))ds+\sum_{j\geq  1} c(\hat{\tau}_j,\hat{\xi}_j )  \cdot 1_{\{\hat{\tau}_j\leq {\rho_m\wedge\tau_S }\}}  \Bigg]
\\&\geq \mathbb{E}\left[G(\hat{Y}^{y_0,\hat{u}} (\rho\wedge\tau_S ))\cdot1_{\{\rho\wedge \tau_S<\infty\}}+\int_{t_0}^{\rho\wedge\tau_s}f(Y^{y_0,\hat{u}} (s))ds+\sum_{j\geq  1} c(\hat{\tau}_j,\hat{\xi}_j )  \cdot 1_{\{\hat{\tau}_j\leq {\rho\wedge\tau_S }\}}  \right],	\end{align*}
where we have used that $\sum_{j=1}^{\mu_{[t_0,\rho_m\vee \tau_S ]} (\hat{u}) } c(\hat{\tau}_j,\hat{\xi}_j ) =\sum_{j\geq  1} c(\hat{\tau}_j,\hat{\xi}_j )  \cdot 1_{\{\hat{\tau}_j\leq {\rho_m\wedge\tau_S }\}}$.
Since this holds for all $\rho \in \mathcal{T}$ we observe that: 
\begin{align}
\phi(y_0)\geq \sup_{\rho \in \mathcal{T}}\mathbb{E}\left[G (\hat{Y}^{y_0,\hat{u}}  (\rho\wedge \tau_S  ))\cdot1_{\{\rho\wedge \tau_S<\infty\}}+ \int_{t_0}^{\rho\wedge \tau_s } f(Y^{y_0,\hat{u}}(s ))ds+ \sum_{j\geq  1}  c (\hat{\tau}_j  , \hat{\xi}_j  )  \cdot 1_{\{\hat{\tau}_j  \leq  {\rho\wedge \tau_s  }\}}  \right].	\end{align}
After which we easily deduce that: \begin{align}
\phi (y_0)\geq \inf_{u \in \mathcal{U}}\sup_{\rho \in \mathcal{T}}\mathbb{E}\left[G (\hat{Y}^{y_0,u} (\rho\wedge \tau_S  ))\cdot1_{\{\rho\wedge \tau_S<\infty\}}+ \int_{t_0}^{\rho\wedge \tau_s } f (Y^{y_0,u} (s ))ds+ \sum_{j\geq  1}  c (\tau_j  , \xi_j  )  \cdot 1_{\{\tau_j  \leq  {\rho\wedge \tau_s  }\}}  \right] \label{zerosumproofpart1finish} .  	
\end{align}
For the second part of the proof, let us fix $\rho' \in \mathcal{T}_{(0,T)} $ as in (\ref{verif_theorem_c-s_zero_item_v_stopping_time_criterion}) and define: \begin{equation}
\rho_D=\rho_D^u= \inf\{s>t_0;X^{t_0,x_0,u} (s)\notin D_2 \}.\end{equation} 	
Now we choose a sequence $\{D_{2,m}\}_{m=1}^{ \infty}$ of open sets s.th. the set $\bar{D}_{2,m} $ is compact with $\bar{D}_{2,m}\subset D_{2,m+1}$ and  $D_2=\cup_{m=1}^{ \infty} D_{2,m}$  and choose $\rho_D (m)=m\wedge \inf_{s>t_0}
X^{\cdot,u} (s)\notin D_{2,m}$. By Dynkin's formula for jump-diffusion processes and (\ref{verif_theorem_c-s_zero_item_iii_HJBI_stopping_inequality}) we have: \begin{align}
\phi(y_0)+\sum_{j=1}^{k} \mathbb{E}[\phi(Y^{y_0,u} (\tau_j ))-\phi(\hat{Y}^{y_0,u} (\tau_j^- ))] -\mathbb{E}[\phi(\hat{Y}^{y_0,u} (\tau_{k+1}^- ))] 	
\\=-\mathbb{E}\left[\int_{t_0}^{\tau_{k+1}} \frac{\partial \phi}{\partial s}+\mathcal{L}\phi(Y^{y_0,u} (s))ds\right]\leq \mathbb{E}\left[\int_{t_0}^{\tau_{k+1}} f(Y^{y_0,u_{[t_0,s]}} (s))ds\right].  	 \end{align}
Hence,
\begin{align}
\phi (y_0)+ \sum_{j=1}^{k}  \mathbb{E}\left[\phi (Y^{y_0,u} (\tau_j  ))-\phi (\hat{Y}^{y_0,u} (\tau_j^-  ))\right]  \leq  \mathbb{E}\left[\int_{t_0}^{\tau_{k+1}} f (Y^{y_0, u_{[t_0,s]}} (s ))ds  +\phi (\hat{Y}^{y_0,u} (\tau_{k+1}^-  ))\right].  	 \label{zerosumineqproof2}
\end{align}
Now by definition of $\mathcal{M}$ we find that: \begin{equation}
\phi(Y^{y_0,u} (\tau_j ))=\phi(\hat{\Gamma}(\hat{Y}^{y_0,u} (\tau_j^- ),\xi_j ))\geq \mathcal{M}\phi(\hat{Y}^{y_0,u} (\tau_j^- ))-c(\tau_j,\xi_j )\cdot1_{\{\tau_j\leq \tau_S\wedge\rho\}}. \label{zerosuminterventionineq}
\end{equation}
(and again using the fact that $\inf_{z\in \mathcal{Z}}[\phi(\tau',\Gamma (X(\tau'^-),z))+c(\tau', z)\cdot 1_{\{\tau'\leq T \}}  ]=0$ whenever $\tau'>\tau_S\wedge\rho$).

Subtracting $\phi(\hat{Y}^{y_0,u} (\tau_j^- )) $ from both sides of (\ref{zerosuminterventionineq}) and summing and negating, we find that: 
\begin{align}\nonumber
&\sum_{j=1}^{k} \mathbb{E}\left[\phi(Y^{y_0,{u}} (\tau_j ))-\phi(\hat{Y}^{y_0,{u}} (\tau_j^- ))\right]  	
\\&\geq \sum_{j=1}^{k} \mathbb{E}\left[\mathcal{M}\phi(\hat{Y}^{y_0,u} (\tau_j^- ))-\phi(\hat{Y}^{y_0,\hat{u}} (\tau_j^- )) -c(\tau_j,\xi_j )\cdot1_{\{\tau_j\leq \tau_S\wedge\rho\}}\right].  	\label{zerosuminterventionineqsummed} 
\end{align}
Inserting (\ref{zerosuminterventionineqsummed}) into (\ref{zerosumineqproof2}) gives: 
\begin{align}
\phi(y_0)+\sum_{j=1}^{k} \mathbb{E}[\mathcal{M}\phi(\hat{Y}^{y_0,u} (\tau_j^- ))-\phi(\hat{Y}^{y_0,u} (\tau_j^- ))] -\mathbb{E}[\phi(\hat{Y}^{y_0,u} (\tau_{k+1}^- ))]\nonumber
\\\leq \mathbb{E}\left[\int_{t_0}^{\tau_{k+1}} f(Y^{y_0,u_{[t_0,s]}} (s))ds+\sum_{j=1}^{k} c(\tau_j,\xi_j )  \cdot 1_{\{\tau_j\leq {\rho\wedge\tau_s }\}}  \right].  	\label{nonzerosumproof2ineqpostsum}\end{align}
Then letting $k\to  \infty$ in (\ref{nonzerosumproof2ineqpostsum}) gives: 
\begin{equation}\begin{split}
\phi(y_0)\leq &-\sum_{j=1}^{\mu_{[t_0,\rho_D (m)\vee \tau_S]} (u)} \mathbb{E}[\mathcal{M}\phi(\hat{Y}^{y_0,u} (\tau_j^- ))-\phi(\hat{Y}^{y_0,u} (\tau_j^- ))] +\mathbb{E}\Bigg[\phi(\hat{Y}^{y_0,u} (\rho_D (m)\wedge\tau_s ))\\&+\int_{t_0}^{\rho_D (m)\wedge\tau_s} f(Y^{y_0,u} (s))ds+\sum_{j\geq  1} c(\tau_j,\xi_j )  \cdot 1_{\{\tau_j\leq {\rho\wedge\tau_s }\}}  \Bigg]\label{nonzerosumproof2ineqpostklim}. 	 \end{split}
\end{equation}
Again, using the quasi-left continuity of $X$ we find that:\\$
 \lim_{m\to \infty}[\mu_{[t_0,\rho_D (m)\vee \tau_S ]} (u)]\equiv \lim_{m\to  \infty} [\mu_{[t_0,\rho_D (m)] } (u)\vee \mu_{[t_0,\tau_S ]} (u)]=\mu_{[t_0,\rho]}  (u)\vee \mu_{[t_0,\tau_S]} (u)$, hence we have that: $\lim_{m\to  \infty} \sum_{j=1}^{\mu_{[t_0,\rho_D (m)\vee \tau_S]}  (u)} \mathbb{E}[\mathcal{M}\phi(\hat{Y}^{y_0,u} (\tau_j^- ))-\phi(\hat{Y}^{y_0,u} (\tau_j^- ))] =0.$
Moreover, as in the first part of the proof, using the fact that $\rho_D (m)\wedge\tau_S\to \rho_D\wedge\tau_s$ as $m\to  \infty $, we can deduce the existence of a constant $c>0$ s.th. 
\begin{align*}
&\lim_{m\to  \infty}\left(\phi(\hat{Y}^{y_0,u} (\rho_D (m)\wedge\tau_s ))+\int_{t_0}^{\rho_D (m)\wedge\tau_s} f(Y^{y_0,u} (s))ds \right)
\\&\leq c \lim_{m\to  \infty} \left(1+|\hat{Y}^{y_0,u} (\rho_D (m)\wedge\tau_s )|\right)+\lim_{m\to  \infty} \left(\int_{t_0}^{\rho_D (m)\wedge\tau_s} \left|Y^{y_0,u} (s)\right|ds\right) 
\\&\leq c(1+\tau_S )(1+\sup_{s \in [0,T]}|X^{t_0,x_0,u} (s)| ) \in \mathbb{L}. 
\end{align*}
Moreover, using (\ref{verif_theorem_c-s_zero_item_vi_continuity_at_termination}), we observe that $\lim_{m\to\infty}\phi(\hat{Y}^{y_0,u} (\rho_D (m)))=\phi(\hat{Y}^{y_0,u} (\rho_D ))=G(\hat{Y}^{y_0,u} (\rho_D )).$ Hence, by the dominated convergence theorem, after taking the limit $m\to  \infty$ in (\ref{nonzerosumproof2ineqpostklim}) we find that: 
\begin{equation}
\phi(y_0)\leq \mathbb{E}\left[\int_{t_0}^{\rho_D\wedge\tau_s}f(Y^{y_0,u} (s))ds+\sum_{j\geq  1} c(\tau_j,\xi_j ) \cdot 1_{\{\tau_j\leq {\rho\wedge\tau_s }\}}  +G(\hat{Y}^{y_0,u} (\rho_D\wedge\tau_S ))\cdot 1_{\{\rho_D\wedge\tau_S<\infty\}}\right].
\end{equation}
Since this holds for all $u \in \mathcal{U}$ we have that: 
\begin{align}\nonumber
\phi(y_0)\leq    \inf_{u \in \mathcal{U}}\mathbb{E}\Bigg[\int_{t_0}^{\rho_D  \wedge \tau_s} f (Y^{y_0,u} (s ))ds  &+ \sum_{j\geq  1}  c (\tau_j  , \xi_j  ) \cdot 1_{\{\tau_j  \leq  {\rho  \wedge \tau_s  }\}}
\\&+G (\hat{Y}^{y_0,u} (\rho_D  \wedge \tau_S  ))\cdot 1_{\{\rho_D\wedge\tau_S<\infty\}}\Bigg], \end{align}
from which clearly we have that: 
\begin{align}\nonumber
\phi (y_0)\leq    \sup_{\rho \in \mathcal{T}}\inf_{u \in \mathcal{U}}\mathbb{E}\Bigg[\int_{t_0}^{\rho\wedge \tau_s } f (Y^{y_0,u} (s ))ds  &+ \sum_{j\geq  1}  c (\tau_j  , \xi_j  ) \cdot 1_{\{\tau_j  \leq \rho\wedge \tau_s  \}}   \\&+G (\hat{Y}^{y_0,u} (\rho\wedge \tau_s  ))\cdot 1_{\{\rho\wedge\tau_S<\infty\}}\Bigg],\label{nonzerosumproof2mainineq} 	\end{align}
where we observe that by (\ref{nonzerosumproof2mainineq}) and (\ref{zerosumproofpart1finish}) we can conclude that: 
\begin{equation}
 \inf_{u \in \mathcal{U}}\left(\sup_{\rho \in \mathcal{T}}J^{(u,\rho)}  [y_0]\right)\leq \phi(y_0)\leq \sup_{\rho \in \mathcal{T}}\left( \inf_{u \in \mathcal{U}}J^{(u,\rho)}  [y_0]\right). 	\label{nonzerosumproof2doubleineq} \end{equation}

However, for all $u \in \mathcal{U}, \rho \in \mathcal{T} $ and $y \in [0,T]\times S$ we have:\\  $\inf_{u \in \mathcal{U}}(\sup_{\rho \in \mathcal{T}}J^{(u,\rho)}  [y])\geq \sup_{\rho \in \mathcal{T}}( \inf_{u \in \mathcal{U}}J^{(u,\rho)}  [y]) $. Moreover, choosing $u=\hat{u}$ in (\ref{nonzerosumproof2doubleineq}), by (\ref{verif_theorem_c-s_zero_item_iv_HJBI_equation}) we find equality, hence: \begin{equation}
\phi(y_0)=\mathbb{E}\left[ \int_{t_0}^{\hat{\rho}\wedge\tau_s} f(Y^{y_0,\hat{u}} (s))ds+\sum_{j\geq  1} c(\hat{\tau}_j,\hat{\xi}_j ) \cdot 1_{\{\tau_j\leq \hat{\rho}\wedge\tau_s \}}  +G(\hat{Y}^{y_0,\hat{u}} (\hat{\rho}\wedge\tau_s ))\right], 
\end{equation}
from which we find that: \begin{equation}
\phi(y)= \inf_{u \in \mathcal{U}}\left(\sup_{\rho \in \mathcal{T}}J^{(u,\rho)}  [y]\right)=\sup_{\rho \in \mathcal{T}}\left( \inf_{u \in \mathcal{U}}J^{(u,\rho)}  [y]\right),\quad \forall y \in [0,T]\times S, \end{equation} 	
from which we deduce the result.
$\hfill \square$
\end{refproof} 

Theorem \ref{Verification_theorem_for_Zero-Sum_Stochastic_Differential_Games_of_Control_and_Stopping_with_Impulse_Controls} provides a constructive means of finding candidate solutions to the game --- in particular, in order to find the value function we seek a function that satisfies the PDEs (\ref{verif_theorem_c-s_zero_item_iii_HJBI_stopping_inequality}) and (\ref{verif_theorem_c-s_zero_item_iv_HJBI_equation}) and thereafter verify that the conditions of the theorem are satisfied. 

The following result expresses the fact that the state space can be divided into three regions and that the players' actions are governed by which region the state process is within:

\begin{corollary}\label{Corollary 2.3.3.}

\noindent The sample space splits into three regions that represent a region in which, when playing their equilibrium strategies, player I applies impulse interventions $I_1$, a region for player II stops the game $I_2$, and a region $I_3$ in which no action is taken by either player; moreover the three regions are characterised by the following expressions: \begin{align*}
 I_1&=\{(t,x) \in [0,T]\times S: V(t,x)=\mathcal{M}V(t,x),\;\mathcal{L}V(t,x)+f(t,x)\geq 0\},\\
I_2&=\{(t,x) \in [0,T]\times S: V(t,x)=G(t,x),\;\mathcal{L}V(t,x)+f(t,x)\leq 0\},\\   
I_3&=\{(t,x) \in [0,T]\times S: V(t,x)<\mathcal{M}V(t,x),\;V(t,x)>G(t,x);\; \mathcal{L}V(t,x)+f(t,x)=0\}. \end{align*}
\end{corollary}
\section{Stochastic Differential Games of Impulse Control and Stopping with Non-Zero-Sum Payoff}\label{section_non_zero_sum_game}

In this section, we study the game as studied in Section \ref{section_zero_sum_game}, however we now extend the results to a non-zero-sum stochastic differential game. The results of this section are loosely based on \cite{[38]} where we make the necessary adjustments to accommodate both impulse controls and the action of the stopper. We start by proving a non-zero-sum verification theorem for the game in which both players use impulse controls to modify the state process.

Suppose firstly that the uncontrolled passive state process evolves according to a (jump-)diffusion \eqref{uncontrolledstateprocess}.

Since we now wish to study non zero-sum games, we decouple the performance objectives so that we now consider the following pair of payoff functions $J_1$ and $J_2$ for player I and player II respectively:
\begin{align}
&\begin{aligned}
J_1^{(u,\rho)} [t_0,x_0]=\mathbb{E}\Bigg[\int_{t_0}^{\rho \wedge \tau_S} f_1 (s,X^{t_0,x_0,u} (s))ds&-\sum_{j\geq 1} c_1 (\tau_j,\xi_j ) \cdot 1_{\{\tau_j\leq \rho \wedge \tau_S \}}
\\&
+G_1 (\rho \wedge \tau_S,X^{t_0,x_0,u} (\rho \wedge \tau_S ))1_{\{\rho \wedge \tau_S <\infty \}} \Bigg], 
\end{aligned}
 \\ &\begin{aligned} 
 J_2^{(u,\rho)} [t_0,x_0]=\mathbb{E}\Bigg[\int_{t_0}^{\rho \wedge \tau_S} f_2 (s,X^{t_0,x_0,u} (s))ds&-\sum_{j\geq 1} c_2 (\tau_j,\xi_j ) \cdot 1_{\{\tau_j\leq \rho \wedge \tau_S \}}   
 \\&
 +G_2 (\rho \wedge \tau_S,X^{t_0,x_0,u} (\rho \wedge \tau_S ))1_{\{\rho \wedge \tau_S <\infty \}} \Bigg],
 \end{aligned}
     \\&\nonumber\hfill\qquad\qquad\qquad\qquad\qquad\qquad\qquad\qquad\qquad\qquad\qquad\qquad\qquad\qquad\qquad\forall (t_0,x_0)\in [0,T]\times S,
\end{align}
where $\tau_j\in \mathcal{T}; \xi_j\in \mathcal{Z}\;(j\in \mathbb{N})$ are $\mathcal{F}-$measurable intervention times and $\mathcal{F}-$measurable measurable stopping interventions respectively. The control $u\in \mathcal{U}$ is an admissible control for player I. The functions $c_1$ and $c_2$ and $G_1:[0,T]\times S\to\mathbb{R}$ and $G_2:[0,T]\times S\to\mathbb{R}$ are cost functions and bequest functions for player I and player II respectively. 

The function $J_1^{(u,\rho)} [t_0,x_0]$ (resp., $J_2^{(u,\rho)} [t_0,x_0]$) defines the payoff received by the player I (resp., player II) during the game with initial point $(t_0,x_0)\in [0,T]\times S$ when player I uses the control $u\in \mathcal{U}$ and player II decides to stop the game at time $\rho \in \mathcal{T}$.

In order to discuss the notion of equilibrium in a non zero-sum case, we must introduce a relevant equilibrium concept which generalises the minimax (saddle point) equilibrium to the non zero-sum case:

\begin{definition}[Nash Equilibrium]\label{Definition 6.1.}
\

\noindent We say that a pair $(\hat{u},\hat{\rho} )\in \mathcal{U}\times \mathcal{T} $ is a Nash equilibrium of the stochastic differential game if the following statements hold:
\begin{enumerate}
    \item $J_1^{(\hat{u},\hat{\rho} )} [t,x]\geq J_1^{(u,\hat{\rho} )} [t,x], \qquad \forall u\in \mathcal{U}\;,\;\forall (t,x)\in [0,T]\times S$,
\item $ J_2^{(\hat{u},\hat{\rho} )} [t,x]\geq J_2^{(\hat{u},\rho )} [t,x], \qquad \forall \rho \in \mathcal{T}\;,\;\forall (t,x)\in [0,T]\times S.$
\end{enumerate}
\end{definition}

Condition (i) states that given some fixed player II stopping time $\hat{\rho}\in \mathcal{T} $, player I cannot profitably deviate from playing the control policy $\hat{u}\in \mathcal{U}$.

Analogously, condition (ii) is the equivalent statement given the player I's control policy is fixed as $\hat{u}$, player II cannot profitably deviate from $\hat{\rho}\in \mathcal{T} $. We therefore see that $(\hat{u},\hat{\rho} )\in \mathcal{U}\times \mathcal{T} $ is an equilibrium in the sense of a Nash equilibrium since neither player has an incentive to deviate given their opponent plays the equilibrium policy.

Before characterising the value functions in this setting, we give a heuristic motivation of the key features of the verification theorem for the game when the payoff structure is non zero-sum. We perform this task by studying the complete repertoire of tactics that each player can employ throughout the horizon of the game.

\subsection*{Heuristic Analysis of The Value Function}

Suppose firstly that each player's value function is sufficiently smooth on the interior of $S$ to apply Dynkin's formula (i.e. we can take first order temporal derivatives and second order spatial derivatives). Suppose also that the following dynamic programming principle is satisfied for each player's value function:
\begin{align}\nonumber
V_i(t_0,x_0)=& \hspace{-1.5 mm}\sup_{u \in \mathcal{U}}   \sup_{\rho \in \mathcal{T}}  \mathbb{E}\Bigg[\int_{t_0}^{(t_0+h)\wedge\rho} f_i(s,X_s^{t_0,x_0,u} )ds-\sum_{j\geq  1} c_i(\tau_j (\rho),\xi_j (\rho)) \cdot 1_{\{\tau_j (\rho)\leq (t_0+h)\wedge\rho\}}
\\&\hspace{-2 mm}\begin{aligned}
+G_i(\rho\wedge\tau_S,X_{\rho\wedge\tau_S}^{t_0,x_0,u} ) \cdot 1_{\{\rho\wedge\tau_S\leq t_0+h\}} +V_i ((t_0+h)\wedge\rho,X_{t_0+h}^{t_0,x_0,u }) \cdot 1_{\{\rho\wedge\tau_S> t_0+h\}}  \Bigg]&,\label{nonviscDPP_non_zero}
\\
i\in \{1,2,\}.& 
\end{aligned}
\end{align}

We firstly tackle the optimality conditions for player I hence, we focus only on the function $J_1$. Let us therefore fix some player II control $\hat{\rho}\in \mathcal{T}$. 

From (\ref{nonviscDPP_non_zero}) and via a classical limiting procedure, we find that the following condition must hold:
\begin{align}
f_1 (t,x)+\partial_t J_1^{(\hat{\rho},\hat{u})} [t,x]+\mathcal{L}J_1^{(\hat{\rho},\hat{u})} [t,x]
\geq f_1 (t,x)+\partial_t J_1^{(\hat{\rho},u)} [t,x]+\mathcal{L}J_1^{(\hat{\rho},u)} [t,x]
\label{ch3heuristicsdifferencenash}\\ t\leq\hat{\tau}_j<s<\hat{\tau}_{j+1}\leq T,\; \forall u\in \mathcal{U},\;\forall (t,x)\in [0,T]\times S.\nonumber	
\end{align}
Expression (\ref{ch3heuristicsdifferencenash}) is an essential constituent of the verification theorem.   

To deduce (\ref{ch3heuristicsdifferencenash}), firstly we note that:
\begin{align}\nonumber
J_1^{(\hat{\rho},u)} [\tau_j,X^{t_0,x_0,u}_{\tau_j}]=\mathbb{E}\Bigg[\int_{\tau_j}^{t'\wedge \hat{\rho}}f_1 (s,X^{t_0,x_0,u}_s )ds  &+J_1^{(u,\hat{\rho})} (t',X_{t'}^{ t_0,x_0,u } ) \cdot1_{\{t'<\hat{\rho} \}}
\\&+G_1 (\hat{\rho}\wedge\tau_S,X_{\hat{\rho}\wedge\tau_S}^{t_0,x_0,u} )\cdot 1_{\{t'\geq \hat{\rho}\} }  \Bigg],  
\label{ch3heuristicsdifferencenash2firstline}\\&
\qquad\qquad t_0 \leq\tau_j<t'<\tau_{j+1}\wedge \tau_S,\quad \forall u \in\mathcal{U}, \nonumber
\end{align}
where we have used the fact that under the policy $u$ no interventions are executed on the interval $]\tau_j,t']$.

We now apply the Dynkin formula for jump-diffusions by which and, with the smoothing theorem, we find that:
\begin{align}
&\quad J_1^{(\hat{\rho},u)} [\tau_j,X^{t_0,x_0,u}_{\tau_j}]-\mathbb{E}\left[J_1^{(\hat{\rho},u)} [\tau_j,X^{t_0,x_0,u}_{\tau_j}]\cdot 1_{\{t'<\hat{\rho} \}}  \right]\nonumber
\\&
\begin{aligned}
=\mathbb{E}\Bigg[\int_{\tau_j}^{t'\wedge \hat{\rho}}{f_1 (s,X_s^{t_0,x_0,u} )+[\partial_t J_1^{(\hat{\rho},u)} (s,X_s^{t_0,x_0,u} )+\mathcal{L}J_1^{(\hat{\rho},u)} (s,X_s^{t_0,x_0,u} )]\cdot 1_{\{t'< \hat{\rho} \}}  }ds&\\  +G_1 (\hat{\rho},X_{\hat{\rho}}^{t_0,x_0,\hat{u}} )\cdot 1_{\{t'\geq \hat{\rho}\} }  \Bigg]&.	 \label{dynkin_appl_ch3_non_zero_sum_heuristic}
\end{aligned}
\end{align}
Since (\ref{dynkin_appl_ch3_non_zero_sum_heuristic}) holds for all $u\in \mathcal{U}$ and, by \eqref{nonviscDPP_non_zero}, in particular for $u=\hat{u}$, we have that: 
\begin{align} 
&\quad
J_1^{(\hat{\rho},\hat{u})} [\hat{\tau}_j,X^{t_0,x_0,\hat{u}}_{\hat{\tau}_j}]-\mathbb{E}[J_1^{(\hat{\rho},\hat{u})} [\hat{\tau}_j,X^{t_0,x_0,\hat{u}}_{\hat{\tau}_j}]\cdot 1_{\{t'<\hat{\rho} \}}  ]
\\&\begin{aligned}
=\mathbb{E}\Bigg[\int_{\hat{\tau}_j}^{t'\wedge \hat{\rho}}f_1 (s,X_s^{t_0,x_0,\hat{u}} )+\left(\partial_t J_1^{(\hat{\rho},\hat{u})} (s,X_s^{t_0,x_0,\hat{u}} )+\mathcal{L}J_1^{(\hat{\rho},\hat{u})} (s,X_s^{t_0,x_0,\hat{u}) })\right) \cdot1_{\{t'<\hat{\rho} \}}  ds&
\\
+G_1 (\hat{\rho},X_{\hat{\rho}}^{t_0,x_0,\hat{u}} )\cdot 1_{\{t'\geq \hat{\rho} \}}  \Bigg]&, 	\label{ch3heuristicsdifferencenash2penultllim} 
\end{aligned}
\end{align}
where we have used the fact that under the policy $\hat{u}$ no interventions are executed on the interval $]\hat{\tau}_j,t']$.

We now make the following observations:
\begin{align}
J_1^{(\hat{\rho},u)} [\tau_j,X^{t_0,x_0,u}_{\tau_j}]-\mathbb{E}\left[J_1^{(\hat{\rho},u)} [\tau_j,X^{t_0,x_0,u}_{\tau_j}] \cdot1_{\{t'<\hat{\rho} \}}\right]  =
\begin{cases}
\begin{aligned}
&0, &t'<\hat{\rho}
\\&J_1^{(\hat{\rho},u)} [\tau_j,X^{t_0,x_0,u}_{\tau_j}],&t'=\hat{\rho}. 
\end{aligned}
\end{cases}
\forall u\in \mathcal{U}.
\end{align}
Additionally, using the optimality of the policy $\hat{u}$ against $J_1[t,x;\hat{\rho},\cdot]$, we have that $J_1^{(\hat{\rho},u)} [\hat{\tau}_j,X^{t_0,x_0,u}_{\hat{\tau}_j}]\leq J_1^{(\hat{\rho},\hat{u})} [\hat{\tau}_j,X^{t_0,x_0,\hat{u}}_{\hat{\tau}_j}],\; \forall u\in \mathcal{U}$. 

Deducting the terms in expression for $J_1^{(\hat{\rho},u)} [\hat{\tau}_j,X^{t_0,x_0,u}_{\hat{\tau}_j}]-\mathbb{E}[J_1^{(\hat{\rho},u)} [\hat{\tau}_j,X^{t_0,x_0,u}_{\hat{\tau}_j}] 1_{\{t'<\hat{\rho} \}}  ]$ in (\ref{ch3heuristicsdifferencenash2firstline}) from (\ref{ch3heuristicsdifferencenash2penultllim}) we readily deduce that:
\begin{align} 
0\leq\mathbb{E}\Bigg[&\int_{\hat{\tau}_j}^{t'\wedge \hat{\rho}} (f_1 (s,X_s^{t_0,x_0,\hat{u}} )-f_1 (s,X_s^{t_0,x_0,u} ))ds
\\&+\int_{\hat{\tau}_j}^{t'\wedge \hat{\rho}}\left(\partial_t J_1^{(\hat{\rho},\hat{u})} [s,X^{t_0,x_0,\hat{u}}_{s}]-\partial_t J_1^{(\hat{\rho},u)} [s,X^{t_0,x_0,u}_{s}]\right)\cdot 1_{\{t'<\hat{\rho} \}}ds\nonumber
\\&\nonumber
\quad+\int_{\hat{\tau}_j}^{t'\wedge \hat{\rho}}\left(\mathcal{L}J_1^{(\hat{\rho},\hat{u})} [s,X_s^{t_0,x_0,\hat{u}  }]-\mathcal{L}J_1^{(\hat{\rho},u)} [s,X_s^{t_0,x_0,u} ]\right)\cdot 1_{\{t'<\hat{\rho} \}}ds \\&\quad\quad
+\sum_{j\geq 1}c_1 (\tau_j,\xi_j ) \cdot 1_{\{\hat{\tau}_j<\tau_j\leq t'\wedge \hat{\rho}\}} 
+\left(G_1 (\hat{\rho},X_{\hat{\rho}}^{t_0,x_0,\hat{u}} )-G_1 (\hat{\rho},X_{\hat{\rho}}^{t_0,x_0,u} )\right)\cdot 1_{\{t'\geq \hat{\rho} \}}  \Bigg].	 
\label{difference_ineq_ch2_heuristic_NE}
\end{align}
Now since \eqref{difference_ineq_ch2_heuristic_NE} holds for all $j=0,1,2,\ldots$ we have in particular for $j=0$:
\begin{align} \nonumber
0\leq\mathbb{E}\Bigg[&\int_{t_0}^{t'\wedge \hat{\rho}}\left(f_1 (s,X_s^{t_0,x_0,\hat{u}} )-f_1 (s,X_s^{t_0,x_0,u} )\right)ds
\\&+\int_{t_0}^{t'\wedge \hat{\rho}}\left(\partial_t J_1^{(\hat{\rho},\hat{u})} [s,X^{t_0,x_0,\hat{u}}_{s}]-\partial_t J_1^{(\hat{\rho},u)} [s,X^{t_0,x_0,u}_{s}]\right)\cdot 1_{\{t'<\hat{\rho} \}}ds\nonumber
\\&\quad
+\int_{t_0}^{t'\wedge \hat{\rho}}\left(\mathcal{L}J_1^{(\hat{\rho},\hat{u})} [s,X_s^{t_0,x_0,\hat{u}  }]-\mathcal{L}J_1^{(\hat{\rho},u)} [s,X_s^{t_0,x_0,u} ]\right)\cdot 1_{\{t'<\hat{\rho} \}}ds\\&
\quad\quad  +\sum_{j\geq 1}c_1 (\tau_j,\xi_j ) \cdot 1_{\{t_0<\tau_j\leq t'\wedge \hat{\rho}\}} \nonumber
+\left(G_1 (\hat{\rho},X_{\hat{\rho}}^{t_0,x_0,\hat{u}} )-G_1 (\hat{\rho},X_{\hat{\rho}}^{t_0,x_0,u} )\right)\cdot 1_{\{t'\geq \hat{\rho} \}}  \Bigg].	 
\end{align}
After taking the limit $t'\downarrow t_0$ 
we arrive at:
\begin{align}
f_1 (t,x)+\partial_t J_1^{(\hat{\rho},\hat{u})} [t,x]+\mathcal{L}J_1^{(\hat{\rho},\hat{u})} [t,x]
\geq f_1 (t,x)+\partial_t J_1^{(\hat{\rho},u)} [t,x]+\mathcal{L}J_1^{(\hat{\rho},u)} [t,x],&\nonumber
\\ \forall (t,x)\in [0,T]\times S,&	\label{ch3heuristicsdifferencenash2penultllimnonzeros} 
\end{align}
as required.

Note also that by similar reasoning as the zero-sum case, we can also deduce that for the control pair $(\hat{\rho},\hat{u})$, we have that:
\begin{equation}
f_1 (t,x)+\partial_t J_1^{(\hat{\rho},\hat{u})} [t,x]+\mathcal{L}J_1^{(\hat{\rho},\hat{u})} [t,x] \leq 0, \qquad \forall (t,x)\in [0,T]\times S,	\label{ch3heuristicsdifferencenash2conteqdirac} \end{equation}
where we recall that the inequality arises since it may be optimal for player I to execute an impulse intervention at the initial point.

For the player II case, we can straightforwardly adapt the arguments from the zero-sum case. Equations (\ref{ch3heuristicsdifferencenash}) and (\ref{ch3heuristicsdifferencenash2conteqdirac}) are central conditions for equilibrium play and appear in the verification theorem as conditions for equilibrium characterisation. 

Having outlined a heuristic argument for the conditions of the verification theorem, we now give a full statement of the theorem. As for the case in Theorem \ref{Verification_theorem_for_Zero-Sum_Stochastic_Differential_Games_of_Control_and_Stopping_with_Impulse_Controls}, the following theorem says that given some pair of solutions to the pair of non-linear PDEs, $i=\{1,2\}$ in (iii), then these solutions coincide with the functions $J_i$ when player $i$ executes their optimal control policy. 

\begin{theorem}[Verification theorem for non zero-sum controller-stopper games ]\label{Verification theorem for Non-Zero-Sum Stochastic Differential Games of Control and Stopping with Impulse Controls}
\

\noindent Let ${\tau_j },\rho \in \mathcal{T}$ be $\mathcal{F}-$measurable stopping times where $j\in\mathbb{N}$. Suppose that there exist functions $\phi_i\in \mathcal{C}^{1,2} ([0,T],S)\cap\mathcal{C}([0,T],\bar{S}),\; i\in \{1,2\}$ s.th.  conditions (T1) - (T4) hold (see the appendix) and additionally:
 \begin{enumerate}[label={(\roman*\rq)}]
\item \label{verif_theorem_c-s_NON_zero_item_ii'_intervention_stopping_inequalities}	$	\phi_1\geq \mathcal{M} \phi_1$ on $S$ and $\phi_2\geq G_2 $ on $S$ and the regions $D_1 $ and $D_2$ are defined by: $ D_1=\{x\in S;\phi_1 (\cdot,x)>\mathcal{M} \phi_1 (\cdot,x)\}$ and $D_2=\{x\in S;\phi_2 (\cdot,x)>G_2 (\cdot,x)\} $
where we refer to $D_1$ (resp., $D_2$) as the player I (resp., player II) continuation region.
\item \label{verif_theorem_c-s_NON_zero_item_iii'_HJBI_inequalities}	$\frac{\partial \phi_1}{\partial s}+\mathcal{L}\phi_1 (s,X^{\cdot ,u} (s))+f_1 (s,X^{\cdot ,u} (s))\leq \frac{\partial \phi_1}{\partial s}+\mathcal{L}\phi_1 (s,X^{\cdot ,\hat{u}} (s))+f_1 (s,X^{\cdot ,\hat{u}} (s))\leq 0\\   \hspace{0.6 mm}$ on $S\backslash \partial D_1$ and $\forall u \in \mathcal{U}$.
\item\label{verif_theorem_c-s_NON_zero_item_iv'_HJBI_equation}	$ \frac{\partial \phi_i}{\partial s}+\mathcal{L}\phi_i (s,X^{\cdot ,\hat{u}} (s))+f_i (s,X^{\cdot ,\hat{u}} (s))=0$ in $D_1,\;\; i\in \{1,2\}.$
\item\label{verif_theorem_c-s_NON_zero_item_v'_stopping_time_criterion}	For $u\in \mathcal{U}$ define $\rho_D=\rho_D^u=\inf\{s>t_0, X^{\cdot ,u} (s)\notin D_2\}$ and specifically, $\hat{\rho}_D  =\hat{\rho}=\inf\{s>t_0, X^{\cdot ,u} (s)\notin D_2\}$.
\item \label{verif_theorem_c-s_NON_zero_item_vi'_continuity_at_termination}	{$X^{\cdot ,u} (\tau_S )\in \partial S,\; \mathbb{P}-$a.s.} on ${\tau_S<\infty }$ and $\phi_i (s,X^{\cdot ,u} (s))\to G_i (\tau_S\wedge \rho,X^{\cdot ,u} (\tau_S\wedge \rho ))$  as  $s\to \tau_S^-\wedge \rho^-$, $  \mathbb{P}-$a.s.$,i\in \{1,2\},\; \forall u\in \mathcal{U}$.
\end{enumerate}	

Put $\hat{\tau}_0\equiv t_0$ and define $\hat{u}:=[\hat{\tau}_j,\hat{\xi}_j ]_{j\in \mathbb{N}}$ inductively by $\hat{\tau}_{j+1}=\inf\{s>\tau_j;X^{\cdot ,\hat{u}} (s)\notin D_1 \}\wedge \tau_S\wedge \rho$,
then $(\hat{u},\hat{\rho} )\in \mathcal{U}\times \mathcal{T}$ is a Nash equilibrium for the game; that is to say that we have:
\begin{align}
\phi_1 (t,x)&=\sup_{u\in \mathcal{U}}  J_1^{(u,\hat{\rho} )} [t,x]= J_1^{(\hat{u},\hat{\rho} )} [t,x], 	\label{nonzerosumresult1.1} \\
\phi_2 (t,x)&=  \sup_{\rho \in \mathcal{T} }  J_2^{(\hat{u},\rho )} [t,x]= J_2^{(\hat{u},\hat{\rho} )} [t,x],	\qquad {\forall (t,x)\in [0,T]\times S}. \label{nonzerosumresult1.2}
\end{align}
\end{theorem}

In full analogy to Corollary \ref{Corollary 2.3.3.}, we can readily arrive at the following corollary to Theorem \ref{Verification theorem for Non-Zero-Sum Stochastic Differential Games of Control and Stopping with Impulse Controls}:

\begin{corollary}\label{Corollary 6.4.} 
\

\noindent When each player plays their equilibrium control, the sample space splits into three regions that represent a region in which the controller performs impulse interventions in $I_1$, a region in which the stopper stops the process $I_2$ and a region in which no action is taken by either player $I_3$; moreover the three regions are characterised by the following expressions:
\begin{align*}
&I_1=\left\{(t,x)\in [0,T]\times S: V_1 (t,x)= \mathcal{M}  V_1 (t,x),\; \mathcal{L}V_1 (t,x)+f_1 (x)\geq 0\right\},  
\\& I_2=\left\{(t,x)\in [0,T]\times S: V_2 (t,x)=G_2 (t,x),\; \mathcal{L}V_2 (t,x)+f_2 (t,x)\geq 0\right\},  
\\& I_3=\Bigg\{(t,x)\in [0,T]\times S:V_1 (t,x)< \mathcal{M}  V_1 (t,x),
\\&\qquad\qquad\qquad V_2 (t,x)<G_1 (t,x);\;\mathcal{L}V_2 (t,x)+f_2 (t,x)=0,j\in \{1,2\} \Bigg\}.
\end{align*}
\end{corollary}
\begin{refproof}[Proof of Theorem \ref{Verification theorem for Non-Zero-Sum Stochastic Differential Games of Control and Stopping with Impulse Controls}.]
\

\noindent As in the proof of Theorem \ref{Verification_theorem_for_Zero-Sum_Stochastic_Differential_Games_of_Control_and_Stopping_with_Impulse_Controls}, we employ the shorthand:
\begin{align}
    &Y^{y_0,\cdot}(s)\equiv (s,X^{t_0,x_0,\cdot}(t_0+s)), \quad y_0\equiv (t_0,x_0), \; \forall s\in [0,T-t_0], \\& \hat{Y}^{y_0,\cdot}(\tau)=Y^{y_0,\cdot}(\tau^{-} )+\Delta_N Y^{y_0,\cdot}(\tau), \quad \tau\in\mathcal{T},
\end{align}
where $\Delta_N Y(\tau) $ denotes a jump at time $\tau$ due to $\tilde{N}$.

Correspondingly, we adopt the following impulse response function $\hat{\Gamma}: \mathcal{T}\times S\times \mathcal{Z}\to  \mathcal{T}\times S$ acting on  $y'\equiv (\tau,x')\in \mathcal{T}\times S$ where $x'\equiv X^{t_0,x_0,\cdot}(t_0+\tau^-)$ is given by: 
\begin{align}
\hat{\Gamma}(y',\zeta)\equiv (\tau,\Gamma (x',\zeta))=(\tau,X^{t_0,x_0,\cdot} (\tau)),\quad \forall \xi\in\mathcal{Z},\; \forall\tau\in\mathcal{T} .
\end{align}

Let us fix the player II control $\hat{\rho}\in \mathcal{T}$; we firstly appeal to the Dynkin formula for jump-diffusions, hence:
\begin{equation}
\mathbb{E}[\phi_1 (\hat{Y}^{y_0,u} (\tau_{j+1} ))]-\mathbb{E}[\phi_1 (Y^{y_0,u} (\tau_j ))]=\mathbb{E}\left[\int_{\tau_j}^{\tau_{j+1}} \frac{\partial \phi_1}{\partial s}+\mathcal{L}[\phi_1 (Y^{y_0,u} (s))]ds\right]. 	\label{nonzerosumproofdynk} 
\end{equation}
Summing (\ref{nonzerosumproofdynk}) from $j=0$ to $j=k$  implies that:
\begin{align}
&-\phi_1 (y_0)-\sum_{j=1}^k \mathbb{E}\left[\phi_1 (Y^{y_0,u} (\tau_j ))-\phi_1 (\hat{Y}^{y_0,u} (\tau_j^- ))\right] +\mathbb{E}\left[\phi_1 (\hat{Y}^{y_0,u} (\tau_{k+1}^- ))\right] 
\\&=\mathbb{E}\left[\int_{t_0}^{\tau_{k+1}} \frac{\partial \phi_1}{\partial s}+\mathcal{L}[\phi_1 (Y^{y_0,u_{[t_0,s]}} (s))]ds\right]. 	\label{nonzerosumproofsumm} 
\end{align}
Now by \ref{verif_theorem_c-s_NON_zero_item_iii'_HJBI_inequalities} we have that: 
\begin{align}
&\qquad\qquad\qquad\qquad\qquad\qquad\frac{\partial \phi_1}{\partial s}+\mathcal{L}\phi_1 (Y^{y_0,u_{[t_0,s]}} (s))\nonumber
\\&\leq \frac{\partial \phi_1}{\partial s}+\mathcal{L}\phi_1 (Y^{y_0,\hat{u}_{[t_0,s]}} (s))+(f_1 (Y^{y_0,\hat{u}} (s))-f_1 (Y^{y_0,u_{[t_0,s]}} (s)))
\leq -f_1 (Y^{y_0,u_{[t_0,s]}} (s)). \label{nonzerosumproofFeqn}	\end{align}
Hence inserting (\ref{nonzerosumproofFeqn}) into (\ref{nonzerosumproofsumm}) yields
\begin{align}\nonumber
-\phi_1 (y_0)-\sum_{j=1}^k \mathbb{E}[\phi_1 (Y^{y_0,u} (\tau_j ))-\phi_1 (\hat{Y}^{y_0,u} (\tau_j^- ))] +\mathbb{E}[\phi_1 (\hat{Y}^{y_0,u} (\tau_{k+1}^- ))] 	\label{nonzerosumresult2} \\
=\mathbb{E}\left[\int_{t_0}^{\tau_{k+1}} \frac{\partial \phi_1}{\partial s}+\mathcal{L}[\phi_1 (Y^{y_0,u_{[t_0,s]}} (s))]ds\right]\leq -\mathbb{E}\left[\int_{t_0}^{\tau_{k+1}} f_1(Y^{y_0,u_{[t_0,s]}} (s))ds\right]. 	 	
\end{align}
Or equivalently:
\begin{align}\nonumber
&\phi_1 (y_0)+\sum_{j=1}^k \mathbb{E}[\phi_1 (Y^{y_0,u} (\tau_j ))-\phi_1 (\hat{Y}^{y_0,u} (\tau_j^- ))] -\mathbb{E}[\phi_1 (\hat{Y}^{y_0,u} (\tau_{k+1}^- ))] 	
\\&\geq \mathbb{E}\left[\int_{t_0}^{\tau_{k+1}} f_1(Y^{y_0,u_{[t_0,s]}} (s))ds\right]. \label{nonzerosumverifineq}
\end{align}
We now use analogous arguments to (\ref{zerosuminterventionineq}) - (\ref{zerosuminterventionineqsummed}). Indeed, by definition of $\mathcal{M} $ we find that: \begin{equation}
\phi_1(Y^{y_0,u} (\tau_j ))=\phi_1(\hat{\Gamma}(\hat{Y}^{y_0,u} (\tau_j^- ),\xi_j ))\leq \mathcal{M} \phi_1(\hat{Y}^{y_0,u} (\tau_j^- ))+c_1(\tau_j,\xi_j )\cdot 1_{\{\tau_j\leq \tau_S\wedge\rho\}}. \label{nonzerosuminterventionineq}
\end{equation}
(using the fact that $\inf_{z\in \mathcal{Z}}[\phi_1(\tau',\Gamma (X(\tau'^-),z))-c_1(\tau', z)\cdot 1_{\{\tau'\leq T \}}  ]=0$ whenever $\tau'>\tau_S\wedge\rho$).

After  subtracting $\phi_1(\hat{Y}^{y_0,u} (\tau_j^- )) $ from both sides of (\ref{nonzerosuminterventionineq}), summing then negating, we find that: 
\begin{align}
&\qquad\qquad\qquad\sum_{j=1}^{k} \mathbb{E}[\phi_1(Y^{y_0,\hat{u}} (\tau_j ))-\phi_1(\hat{Y}^{y_0,\hat{u}} (\tau_j^- ))]  	
\\&
\leq \sum_{j=1}^{k} \mathbb{E}[\mathcal{M} \phi_1(\hat{Y}^{y_0,u} (\tau_j^- ))-\phi_1(\hat{X}^{t_0,x_0,\hat{u}} (\tau_j^- ))] +\sum_{j=1}^{k} \mathbb{E}[c_1(\tau_j,\xi_j )\cdot 1_{\{\tau_j\leq \tau_S\wedge\rho\}}] . 	 	\label{nonzerosuminterventionineqsummed} \end{align}
After inserting (\ref{nonzerosuminterventionineqsummed}) into (\ref{nonzerosumverifineq}) we find that:
\begin{align}
&\begin{aligned}\phi_1 (y_0)\geq \mathbb{E}\Bigg[\phi_1 (\hat{Y}^{y_0,u} (\tau_{k+1}^- ))-\sum_{j=1}^k [\phi_1 (Y^{y_0,u} (\tau_j ))&-\phi_1 (\hat{Y}^{y_0,u} (\tau_j^- ))]\nonumber
\\&+\int_{t_0}^{\tau_{k+1}} f_1(Y^{y_0,u_{[t_0,s]}} (s))ds\Bigg]\nonumber\end{aligned}
\\&\qquad\quad\begin{aligned}\geq \mathbb{E}\Bigg[\phi_1 (\hat{Y}^{y_0,u} (\tau_{k+1}^- ))&-\sum_{j=1}^k [\mathcal{M} \phi_1 (Y^{y_0,u} (\tau_j^- ))-\phi_1 (\hat{Y}^{y_0,u} (\tau_j^- ))]\\&-\sum_{j=1}^k c_1 (\tau_j,\xi_j )\cdot 1_{\{\tau_j\leq \tau_S\wedge\rho\}} +\int_{t_0}^{\tau_{k+1}} f_1(Y^{y_0,u_{[t_0,s]}} (s))ds\Bigg].\label{nonzerosumsumineq}
\end{aligned}
\end{align}
Define $\hat{\rho}_m\equiv \hat{\beta}_m (u)=\hat{\rho}\wedge m;m=1,2\ldots$.  As in the zero-sum case, since the number of interventions in (\ref{nonzerosumsumineq}) is bounded above by $\mu_{[t_0,\hat{\rho}_m\wedge \tau_S]} (u)\wedge m$ for some $m<\infty$  and (\ref{nonzerosumsumineq}) holds for any $k\in \mathbb{N}$, taking the limit as $k\to \infty$  in (\ref{nonzerosumsumineq}) gives:
\begin{align}\nonumber
\phi_1 (y_0)\geq \mathbb{E}\Bigg[\phi_1 (\hat{Y}^{y_0,u} (\tau_{k+1}^- ))&-\sum_{j=1}^{\mu_{[t_0,\hat{\rho}_m\wedge \tau_S]} (u)\wedge m} [\mathcal{M} \phi_1 (Y^{y_0,u} (\tau_j^- ))-\phi_1 (\hat{Y}^{y_0,u} (\tau_j^- ))] \\&+\int_{t_0}^{\hat{\rho}_m\wedge \tau_s} f_1(Y^{y_0,u_{[t_0,s]}} (s))ds
-\sum_{j\geq 1} c_1 (\tau_j,\xi_j )  \cdot 1_{\{\tau_j\leq {\rho \wedge \tau_s }\}}  \Bigg]. 
\label{nonzerosumsumineqklim}	 
\end{align}
Now, $\lim_{m\to \infty} \sum_{j=1}^{\mu_{[t_0,\hat{\rho}_m\wedge \tau_S]} (u)\wedge m} \mathbb{E}[\mathcal{M} \phi_1 (Y^{y_0,u} (\tau_j ))-\phi_1 (\hat{Y}^{y_0,u} (\tau_j ))] \\
=\sum_{j=1}^{\mu_{[\hat{\rho}\wedge \tau_S ]} (u)} \mathbb{E}[\mathcal{M} \phi_1 (Y^{y_0,u} (\tau_j ))-\phi_1 (\hat{Y}^{y_0,u} (\tau_j ))] =0$ and $\lim_{m\to \infty}[\phi_1 (\hat{Y}^{y_0,u} (\tau_{\mu_{[t_0,\hat{\rho}_m\wedge \tau_S]} } (u)\wedge m) ))]=\phi_1 (\hat{Y}^{y_0,u} (\hat{\rho}\wedge \tau_S ))=G_1 (\hat{Y}^{y_0,u} (\hat{\rho}\wedge \tau_S )).$ Indeed, by \ref{verif_theorem_c-s_NON_zero_item_v'_stopping_time_criterion} we have that\\ $\lim_{m\to \infty }[\tau_{\mu_{[t_0,\hat{\rho}_m\wedge \tau_S]}  (u)}\wedge m) ]=\tau_{\mu_{[t_0,\hat{\rho}\wedge \tau_s ]} (u)}\equiv \hat{\rho}\wedge \tau_s$ . 
Thus, after taking the limit $k,m\to \infty$ in (\ref{nonzerosumsumineqklim})  and noting that by definition, $\lim_{m\to \infty}\hat{\rho}_m=\hat{\rho}$, we have that:
\begin{align}\nonumber
\phi_1 (y_0)\geq \mathbb{E}\Bigg[\int_{t_0}^{\hat{\rho}\wedge \tau_s} f_1(Y^{y_0,u_{[t_0,s]}} (s))ds&-\sum_{j\geq 1} c_1 (\tau_j,\xi_j )  \cdot 1_{\{\tau_j\leq {\hat{\rho}\wedge \tau_s }\}} 
\\&+G_1 (\hat{Y}^{y_0,u} (\hat{\rho}\wedge \tau_s ))\cdot 1_{\{\hat{\rho}\wedge \tau_s<\infty\}}\Bigg].  \end{align}
Since this holds for all $u\in \mathcal{U}$ we find:
\begin{align}\nonumber
\phi_1 (y_0)\geq \sup_{u\in \mathcal{U}}\mathbb{E}\Bigg[G_1 (\hat{Y}^{y_0,u} (\hat{\rho}\wedge \tau_s ))\cdot 1_{\{\hat{\rho}\wedge \tau_s<\infty\}}&+\int_{t_0}^{\hat{\rho}\wedge \tau_s} f_1(Y^{y_0,u_{[t_0,s]}} (s))ds
\\&-\sum_{j\geq 1} c_1 (\tau_j,\xi_j )  \cdot 1_{\{\tau_j\leq {\hat{\rho}\wedge \tau_s }\}}  \Bigg]. 	 
\end{align}
Hence, we find that
\begin{equation}
\phi_1 (y)\geq \sup_{u\in \mathcal{U}}  J_1^{(u,\hat{\rho} )} [y], \qquad \forall y\in [0,T]\times S.
\end{equation}
Now, applying the above arguments with the controls $(\hat{u},\hat{\rho})$ yields the following equality:
\begin{align}
\phi_1 (y)=\sup_{\rho \in \mathcal{T} } J_1^{(\hat{u},\rho )} [y]=J_1^{(\hat{u},\hat{\rho} )} [y],\label{nonzerosumeq1} \qquad\forall y\in [0,T]\times S.
\end{align}
To prove (\ref{nonzerosumresult1.1}) - (\ref{nonzerosumresult1.2}), we firstly fix $\hat{u}\in \mathcal{U}$ as in \ref{verif_theorem_c-s_NON_zero_item_iv'_HJBI_equation}, we again define $\rho_m=\rho \wedge m;m=1,2\ldots$.  Now, by the Dynkin formula for jump-diffusions and by \ref{verif_theorem_c-s_NON_zero_item_iv'_HJBI_equation} and (\ref{nonzerosumproofdynk}) - (\ref{nonzerosumproofsumm}), we have that:
\begin{align*}
\mathbb{E}&[\phi_2 (\hat{Y}^{y_0,\hat{u}} (\rho_m^- ))]-\phi_2 (y_0)-\sum_{j=1}^{\mu_{[t_0,\rho_m]} (\hat{u})} \mathbb{E}[\phi_2 (Y^{y_0,\hat{u}} (\hat{\tau}_j ))-\phi_2 (\hat{Y}^{y_0,\hat{u}} (\hat{\tau}_j^- ))]  	
\\=\mathbb{E}&\left[\int_{t_0}^{\hat{\tau}_{\mu_{[t_0,\rho_m]} (\hat{u})}} \left(\frac{\partial \phi_2}{\partial s}+\mathcal{L}\phi_2 (Y^{y_0,\hat{u}_{[t_0,s]}} (s))ds\right)\right] =-\mathbb{E}\left[\int_{t_0}^{\hat{\tau}_{\mu_{[t_0,\rho_m]} (\hat{u})}} f_2(Y^{y_0,\hat{u}_{[t_0,s]}} (s))ds\right],	 \end{align*}
which (as before, similar to (\ref{interventioneq1})) and by our choice of $\hat{\xi}_j\in \mathcal{Z}$, implies 
\begin{align*}
&\phi_2 (y_0)+\sum_{j=1}^{\mu_{[t_0,\rho_m]} (\hat{u})} \mathbb{E}[\mathcal{M} \phi_2 (Y^{y_0,\hat{u}} (\hat{\tau}_j^- ))-\phi_2 (\hat{Y}^{y_0,\hat{u}} (\hat{\tau}_j^- ))]  
\\&\begin{aligned}	= \mathbb{E}\Bigg[\phi_2  (\hat{Y}^{y_0,\hat{u}} (\hat{\tau}_{\mu_{[t_0, \rho_m]} (\hat{u})} ))&+ \int_{t_0}^{\hat{\tau}_{\mu_{[t_0, \rho_m]} ( \hat{u}) }} f_2  (Y^{y_0, \hat{u}_{[t_0,s ]}} (s ))ds \\&-\sum_{j=1}^{\mu_{[t_0,\rho_m]} ( \hat{u}) } c_2  (\hat{\tau}_j  , \hat{\xi}_j  )  \cdot 1_{\{\hat{\tau}_j  \leq  \rho_m  \}}  \Bigg],\end{aligned} 	 \end{align*}
which we may rewrite as
\begin{align}
\nonumber \phi_2 (y_0)=\mathbb{E}\Bigg[\phi_2 (\hat{Y}^{y_0,\hat{u}} (\hat{\tau}_{\mu_{[t_0,\rho_m]} (\hat{u} ) } ))&+\int_{t_0}^{\hat{\tau}_{\mu_{[t_0,\rho_m]} (\hat{u} ) }} f_2(Y^{y_0,\hat{u}_{[t_0,s]}} (s))ds
\\&\nonumber-\sum_{j=1}^{\mu_{[t_0,\rho_m]} (\hat{u} )} c_2 (\hat{\tau}_j,\hat{\xi}_j )  \cdot 1_{\{\hat{\tau}_j\leq \tau_{\mu_{[t_0,\rho_m]} (\hat{u} ) } \} }  \Bigg]
\\&-\sum_{j=1}^{\mu_{[t_0,\rho_m]} (\hat{u})} \mathbb{E}\left[\mathcal{M} \phi_2 (Y^{y_0,\hat{u}} (\tau_j^- ))-\phi_2 (\hat{Y}^{y_0,\hat{u}} (\tau_j^- ))\right] . 	\label{verfsummedeq2} 
\end{align}
Now, since $\mu_{[t_0,\rho_m]} (\hat{u})\to \mu_{[t_0,\rho \wedge \tau_S]} (\hat{u})$ as $m\to \infty$  and $\lim_{s\to \tau_S }\phi_i (\hat{Y}^{y_0,\cdot} (s))\to G_i (\hat{Y}^{y_0,\cdot} (\tau_S )),\\i\in \{1,2\}$ using (v) and $\hat{\tau}_{\rho \wedge \tau_S}\equiv \rho \wedge \tau_S$ then using (\ref{verfsummedeq2}) and by the Fat\^{o}u lemma we find that:
\begin{align}\nonumber
\phi_2&(y_0)\geq \lim \inf_{m\to \infty}\mathbb{E}\Bigg[\phi_2 (\hat{Y}^{y_0,\hat{u}} (\hat{\tau}_{\mu_{[t_0,\rho_m]} (\hat{u}})  ))+\int_{t_0}^{\hat{\tau}_{\mu_{[t_0,\rho_m]} (\hat{u}})} f_2(Y^{y_0,\hat{u}_{[t_0,s]}} (s))ds
\\\nonumber&-\sum_{j\geq 1} c_2 (\hat{\tau}_j,\hat{\xi}_j )  \cdot 1_{\{\tau_j\leq \hat{\tau}_{\mu_{[t_0,\rho_m]} (\hat{u}})  \}} -\sum_{j=1}^{\mu_{[t_0,\rho_m]} (\hat{u})} \mathbb{E}[\mathcal{M} \phi_2 (Y^{y_0,\hat{u}} (\tau_j ))-\phi_2 (\hat{Y}^{y_0,\hat{u}} (\tau_j ))] ]   	
\\&=\mathbb{E}\left[G_2 (Y^{y_0,u} (\rho \wedge \tau_S ))\cdot 1_{\{{\rho}\wedge \tau_s<\infty\}}+\int_{t_0}^{\rho \wedge \tau_S} f_2 (Y^{y_0,u} (s))ds-\sum_{j\geq 1} c_2 (\hat{\tau}_j,\hat{\xi}_j )  \cdot 1_{\{\hat{\tau}_j\leq \rho \wedge \tau_S \}} \right]. 	\end{align}
Since this holds for all $\rho \in \mathcal{T}$ we find that:
\begin{align}
\phi_2 (y_0)\geq \sup_{\rho \in \mathcal{T}  }\mathbb{E}\Bigg[G_2 (Y^{y_0,u} (\rho \wedge \tau_S ))\cdot 1_{\{{\rho}\wedge \tau_s<\infty\}}+&\int_{t_0}^{\rho \wedge \tau_S} f_2 (Y^{y_0,\hat{u}} (s))ds
\\&-\sum_{j\geq 1} c_2 (\hat{\tau}_j,\hat{\xi}_j )  \cdot 1_{\{\hat{\tau}_j\leq \rho \wedge \tau_S \}}  \Bigg]. 	\end{align}
Hence, we find that
\begin{equation}  
\phi_2 (y)\geq \sup_{\rho \in \mathcal{T}  } J_2^{(\hat{u},\rho )} [y], \qquad \forall y\in [0,T]\times S.
\label{nonzerosineq2}
\end{equation}

Now, applying the above arguments with the controls $(\hat{u},\hat{\rho})$ yields the following equality:
\begin{equation}
\phi_2 (y)=\sup_{u\in \mathcal{U} }  J_2^{(u,\hat{\rho} )} [y]=J_2^{(\hat{u},\hat{\rho} )} [y], \qquad \forall y\in [0,T]\times S.
\label{nonzerosequality} 
\end{equation}

We therefore observe using (\ref{nonzerosequality}) in conjunction with (\ref{nonzerosumeq1}) and that $(\hat{u},\hat{\rho})$ is a Nash equilibrium and hence the thesis is proven.
$\hfill \square$\end{refproof}

\section{Examples}\label{section_examples}
In order to demonstrate the workings of the theorems, we give some example calculations. 

\begin{example}\label{ch_2_zero_sum_example}
Consider a system with passive dynamics that are described by a stochastic process $X_s=X(s,\omega ):[0,T]\times \Omega\to\mathbb{R}$ which obeys the following SDE:
\begin{align}
dX(r)=X(r^-)(\alpha dr+\beta dB(r)),\qquad \forall r\in ]0,T], \label{example_calc_ch2_zs_state_process}
\end{align}
where $\alpha,\beta\in\mathbb{R}_{>0}$ are fixed constants, $B(r)$ is a 1-dimensional Brownian motion and $T\in\mathbb{R}_{>0}$ is some finite time horizon. The state process (\ref{example_calc_ch2_zs_state_process}) is  \textit{geometric Brownian motion}. Geometric Brownian motion is widely used to model various financial processes \cite{hull2016options} and is a particular case of geometric L\'evy process (c.f. the optimal liquidity control and lifetime ruin problem) that is restricted to have continuous sample paths.

The state process $X$ is modified by a controller, player I that exercises an impulse control policy $u=[\tau_j,\xi_j]\in\mathcal{U}$. Additionally, at any point $\rho<T$ a second player, player II can choose to stop the process where $\rho\in\mathcal{T}$ is an $\mathcal{F}-$measurable stopping time. The controlled state process therefore evolves according to the following expression:
\begin{align}
X(r)=x_0+ \alpha\int_0^{r\wedge\rho} X(s) ds+\beta\int_0^{r\wedge\rho} X(s)dB(s)-\sum_{j\geq 1}(\kappa_1+(1+\lambda)\xi_j)\cdot 1_{\{\tau_j\leq\rho\wedge {r}\}},&\label{example_calc_ch2_zs_control_state_process}
\\\forall r\in [0,T],\mathbb{P}-{\rm a.s.},&\nonumber
\end{align}
where $\kappa_1>0$ and $\lambda>0$ are the fixed part and the proportional part of the transaction cost incurred by player I for each intervention (resp.).

Player I seeks to choose an admissible impulse control $u=[\tau_j,\xi_j]$ that maximises its reward $J$ where $\{\tau_j\}_{\{j\geq 1\}}$ are intervention times and each $\xi_{j\geq 1}\in\mathcal{Z}$ is an impulse intervention. Player II seeks to choose an $\mathcal{F}-$ measurable stopping time $\rho\in \mathcal{T}$ that minimises the same quantity $J$ which is given by the following expression:
\begin{align}
J^{\rho,u}[s,x]=\mathbb{E}\left[e^{-\delta(s+\rho)}(X(\rho)-\kappa_2)+\sum_{j\geq 1}e^{-\delta(s+\tau_j)}\xi_j\cdot 1_{\{\tau_j\leq\rho\wedge T\}} \right], \qquad \forall (s,x)\in[0,T]\times\mathbb{R}.
\end{align}

An example of a setting for this game is an interaction between a project manager (player I) that seeks to maximise project investments $\{\xi_j\}_{\{j\geq 1\}}$ over some time horizon $T$, and a second interested party (player II), e.g. a firm owner, that can choose to terminate the project at any point $\rho\leq T$. Whenever the firm owner chooses to terminate the project, they receive a discounted payment of $\kappa_2$. The owner however, seeks to terminate the project when the unspent cash flow $X$ is minimal.

The problem is to find a function $\phi\in\mathcal{C}^{1,2}([0,T],\mathbb{R})$ s.th.
\begin{align}
\inf_{\rho}\;\sup_{u}J^{\rho,u}[s,x]=\sup_{u}\;\inf_{\rho}J^{\rho,u}[s,x]=\phi(s,x),\quad \forall (s,x)\in[0,T]\times\mathbb{R},
\end{align}
We recognise this as a zero-sum stochastic game of impulse control and stopping. We therefore seek to apply Theorem \ref{Verification_theorem_for_Zero-Sum_Stochastic_Differential_Games_of_Control_and_Stopping_with_Impulse_Controls} to compute the equilibrium controls. 

Firstly, we observe that by (\ref{example_calc_ch2_zs_state_process}) and using (\ref{generator}), the generator $\mathcal{L}$ for the process $X$ is given by: 
\begin{align}
\mathcal{L}\psi(s,x)=\frac{\partial \psi}{\partial s}(s,x)+\alpha x \frac{\partial \psi}{\partial x}(s,x)+\frac{1}{2}\beta^2x^2\frac{\partial^2 \psi}{\partial x^2}(s,x),    \label{example_ch2_generator}
\end{align}
for some test function $\psi\in\mathcal{C}^{1,2}([0,T],\mathbb{R})$.

We now wish to derive the functional form of the function $\phi$. Applying (\ref{verif_theorem_c-s_zero_item_iv_HJBI_equation}) of Theorem \ref{Verification_theorem_for_Zero-Sum_Stochastic_Differential_Games_of_Control_and_Stopping_with_Impulse_Controls} leads to the \textit{Cauchy-Euler equation} $\mathcal{L}\phi=0$, (here, $f\equiv 0$ in Theorem \ref{Verification_theorem_for_Zero-Sum_Stochastic_Differential_Games_of_Control_and_Stopping_with_Impulse_Controls}). Following this, we make the following ansatz: $\phi(s,x)=e^{-\delta s}\psi(x)$ where $\psi(x):=ax^c$ for some as yet, undetermined constants $a,c\in\mathbb{R}$.

Plugging the ansatz for the function $\phi$ and using (\ref{verif_theorem_c-s_zero_item_iv_HJBI_equation}) of Theorem \ref{Verification_theorem_for_Zero-Sum_Stochastic_Differential_Games_of_Control_and_Stopping_with_Impulse_Controls} into (\ref{example_ch2_generator}) immediately gives:
\begin{align}
-\delta+\alpha c +\frac{1}{2}\beta^2(c-1)c =0.    
\end{align}
After some manipulation, we deduce that there exist two solutions for $c$ which we denote by $c_+$ and $c_-$ s.th. $c_+>c_-$ with $c_+>0$ and $|c_-|>0$ which are given by the following:
\begin{align}
c_{\pm}=-\frac{\alpha-\frac{1}{2}\beta^2}{\beta^2}\pm\frac{1}{\beta^2}\sqrt{(\alpha-\frac{1}{2}\beta^2)^2+2\beta^2\delta}.\label{c_values_ch2_example}
\end{align}
We now apply the HJBI equation (\ref{verif_theorem_c-s_zero_item_iv_HJBI_equation}) of Theorem \ref{Verification_theorem_for_Zero-Sum_Stochastic_Differential_Games_of_Control_and_Stopping_with_Impulse_Controls} to characterise the function $\phi$ on the region $D_1\cap D_2$. Following our ansatz, we  observe that by (\ref{verif_theorem_c-s_zero_item_iv_HJBI_equation}) the following expression for the function $\phi$: 
\begin{align}
\phi(s,x)&=e^{-\delta s}\psi(x),\qquad &&\forall (s,x)\in[0,T]\times D_1\cap D_2,\\
\psi(x)&=(a_1x^{c_+}+a_2x^{c_-}), \qquad &&\forall x\in D_1\cap D_2, \label{ch_2_phi_a}
\end{align}
where $a_1$ and $a_2$ are constants that are yet to be determined and $D_1$ and $D_2$ are the continuation regions for player I and player II respectively.

In order to determine the constants $a_1$ and $a_2$, we firstly observe that $\phi(\cdot,0)=0$. This then implies that $a_1=-a_2:=a$. We now deduce that the function $\psi$ is given by the following expression:
\begin{align}
\psi(x)=a(x^{c_+}-x^{c_-}),\qquad \forall x\in D_1\cap D_2.      
\end{align}
In order to characterise the function over the entire state space and find the value $a$, using conditions (\ref{verif_theorem_c-s_zero_item_ii_intervention_stopping_inequalities}) - (\ref{verif_theorem_c-s_zero_item_vi_continuity_at_termination}) of Theorem \ref{Verification_theorem_for_Zero-Sum_Stochastic_Differential_Games_of_Control_and_Stopping_with_Impulse_Controls}, we study the behaviour of the function $\phi$ given each player's control.

Firstly, we consider the player I impulse control problem. In particular, we seek conditions on the impulse intervention applied when $\mathcal{M}\phi=\phi$. To this end, let us firstly conjecture that the player I continuation region $D_1$ takes the following form:
\begin{equation}
    D_1=\{x\in\mathbb{R}; 0<x<\tilde{x}\},\label{ch2_example_cont_region_d1}
\end{equation}
for some constant $\tilde{x}$ which we shall later determine.

Our first task is to determine the optimal value of the impulse intervention. We now define the following two functions which will be of immediate relevance:
\begin{align}
&\psi_0(x):=a(x^{c_+}-x^{c_-}),\\
&h(\xi):=\psi(x-\kappa_1-(1+\lambda)\xi)+\xi,\\&\qquad\qquad\qquad\qquad\qquad\qquad\qquad\qquad\qquad\qquad\qquad\qquad\qquad\qquad
\forall x\in \mathbb{R},\forall \xi\in\mathcal{Z}.\nonumber
\end{align}
In order to determine the value $\hat{\xi}$ that maximises $\Gamma(x({\tau}^- ),\xi)$ at the point of intervention, we investigate  the first order condition on $h$ i.e. $0=h'(\xi)$. This implies the following:
\begin{align}
\psi'(\tilde{x}-\kappa_1-(1+\lambda)\xi)&=\frac{1}{1+\lambda}.\label{impulse_1st_condition}
\end{align}
Using the expression for $\psi$ (\ref{ch_2_phi_a}) we also observe the following:
\begin{align}
\psi'_0(x)= c_+{x}^{c_+-1}&-c_-{x}^{c_--1}>0,\qquad&&\forall x\in\mathbb{R},  \label{ch_2_example_FOC}
\\
\psi''_0(x)=c_+(c+-1){x}^{c_+-2}&-c_-(c_--1){x}^{c_--1}<0,&&\forall x<x^{\#}:=\Big|\frac{c_-(c_--1)}{c_+(c_+-1)}\Big|^{\frac{1}{c_+-c_-}},
\label{ch_2_example_SOC}
\end{align}
from which we deduce the existence of two points $x^{\star},x_{\star}$ for which the condition $\psi'_0(\cdot)=(1+\lambda)^{-1}$ holds. W.l.o.g. we assume $x^{\star}>x_{\star}$. Now by (\ref{verif_theorem_c-s_zero_item_ii_intervention_stopping_inequalities}) of Theorem \ref{Verification_theorem_for_Zero-Sum_Stochastic_Differential_Games_of_Control_and_Stopping_with_Impulse_Controls} we require that $e^{-\delta s}\psi(x)=\mathcal{M}e^{-\delta s}\psi_0(x)$ for any $s\in [0,T]$ whenever $x\geq \tilde{x}$ (c.f. $D_1$ in equation (\ref{ch2_example_cont_region_d1})), hence we find that:
\begin{equation}
    \psi(x)=\psi_0(x_{\star})+\hat{\xi}(x), \qquad \forall x\geq \tilde{x},
\end{equation}
where $x-\kappa_1-(1+\lambda)\hat{\xi}(x)=x_{\star}$ from which we readily find that the optimal impulse intervention value is given by:
\begin{equation}
    \hat{\xi}(x)=\frac{x-x_{\star}-\kappa_1}{1+\lambda},\qquad \forall x\geq \tilde{x}.\label{ch_2_example_optim_intervention}
\end{equation}
Having determined the optimal impulse intervention and constructed the form of the continuation region for Player I, we now turn to the optimal stopping criterion for Player II. 

We conjecture that the continuation region for player II, $D_2$,  takes the following form:
\begin{equation}
    D_2=\{x\in\mathbb{R};\; x>\hat{x}\}.\label{ch2_example_cont_region_d2}
\end{equation}
Now using condition (\ref{verif_theorem_c-s_zero_item_vi_continuity_at_termination}) of Theorem \ref{Verification_theorem_for_Zero-Sum_Stochastic_Differential_Games_of_Control_and_Stopping_with_Impulse_Controls}
we observe the following:
\begin{align}
    \psi(x) =(x-\kappa_2),\qquad \forall x\notin  D_2.
\end{align}
Additionally, we recall that by (\ref{ch_2_phi_a}) we have the following
\begin{align}
\psi(x) =
a(x^{c_+}-x^{c_-}),\qquad \forall x\in D_1\cup D_2,
\end{align}
where the constant $a$ is to be determined.

Putting the above facts together we can give a characterisation for the function $\psi$:
\begin{align}
\psi(x) =
\begin{cases}
\begin{aligned}
&a(x^{c_+}-x^{c_-}), &\forall x\in D_1\cap D_2,
\\&(x-\kappa_2), &\forall x\notin  D_2,
\\&a(x^{c_+}_{\star}-x^{c_-}_{\star})+\frac{x-x_{\star}-\kappa_1}{1+\lambda}, &\forall x\notin D_1,
\end{aligned}
\end{cases}
\end{align}
where the constants $c_+$ and $c_-$ are specified in equation (\ref{c_values_ch2_example}).

Using the facts above, we are now in a position to determine the value of the constants $a,\hat{x}$ and $\tilde{x}$. To do this, we assume the \textit{high contact principle} --- a condition that asserts the continuity of the value function at the boundary of the continuation region (for exhaustive discussions on the condition, see \cite{oksendalapplied2007, oksendal1989high}). 

For player II, using (\ref{verif_theorem_c-s_zero_item_vi_continuity_at_termination}) it then follows that the following conditions must hold:
\begin{align}
&\phi(\cdot,\hat{x})=G(\cdot,\hat{x})&&\implies  a(\hat{x}^{c_+}-\hat{x}^{c_-})&&\hspace{-13 mm}=\hat{x}-\kappa_2, \label{high_contact_1st_order}
\\&\phi'(\cdot,\hat{x})=G'(\cdot,\hat{x})&&\implies a(c_+\hat{x}^{c_+-1}-c_-\hat{x}^{c_--1})&&\hspace{-13 mm}=1, \label{high_contact_2nd_order}
\end{align}
Using (\ref{high_contact_1st_order}) - (\ref{high_contact_2nd_order}) we deduce that the value $a$ is given by:
\begin{align}
a=\kappa_2[(1-c_-)\hat{x}^{c_-}-(1-c_+)\hat{x}^{c_+})]^{-1}.\label{ch2_example_a}
\end{align}
Additionally, by (\ref{high_contact_1st_order}) - (\ref{high_contact_2nd_order}) we find that the value of $\hat{x}$ is the solution to the equation:
\begin{equation}
    p(\hat{x})=0, \label{hat_x_ch2_example_solution}
\end{equation}
where $p(x)={\kappa_2}^{-1}[(1-c_-){x}^{c_--1}-(1-c_+){x}^{c_+-1})](c_+{x}^{c_+-1}-c_-{x}^{c_--1})^{-1}$.

Lastly, we apply the high contact principle to find the boundary of the continuation region $D_1$. Indeed, continuity at $\tilde{x}$ leads to the following:
\begin{align}
\psi(\tilde{x})=\psi_0(x_{\star})+\hat{\xi}(\tilde{x}),  \implies a(\tilde{x}^{c_+}-\tilde{x}^{c_-})=a(x_{\star}^{c_+}-x_{\star}^{c_-})+\frac{\tilde{x}-x_{\star}-\kappa_1}{1+\lambda},
\end{align}
from which we find that $\tilde{x}$ is the solution to the following equation:
\begin{align}
    m(\tilde{x})&=0,\label{ch_2_example_p1_m_function}
    \\
        m(x)&=x-a(1+\lambda)[x^{c_+}-x^{c_-}+x^{c_-}_{\star}-x^{c_+}_{\star}]-x_{\star}+\kappa_1. \label{ch_2_example_p1_cont_region_value}
\end{align}

Equations (\ref{hat_x_ch2_example_solution}) and (\ref{ch_2_example_p1_m_function}) are difficult to solve analytically for the general case but can however, be straightforwardly solved numerically using a root-finding algorithm.

To summarise, the solution is as follows: whenever $X\in D_1\cap D_2$ neither player intervenes. Player I performs an impulse intervention of size $\hat{\xi}$ given by (\ref{ch_2_example_optim_intervention}) whenever the process reaches the value $\tilde{x}$ and player II terminates the game if the process hits the value $\hat{x}$. The value function for the problem is $\phi(s,x)\equiv e^{-\delta s}\psi(x),\; \forall (s,x)\in[0,T]\in\mathbb{R}$, where is $\psi$ given by:
\begin{align}
\psi(x) =
\begin{cases}
\begin{aligned}
&a(x^{c_+}-x^{c_-}), &\forall x\in D_1\cap D_2,
\\&(x-\kappa_2), &\forall x\notin  D_2,
\\&a(x^{c_+}_{\star}-x^{c_-}_{\star})+\frac{x-x_{\star}-\kappa_1}{1+\lambda}, &\forall x\notin D_1,
\end{aligned}
\end{cases}
\end{align}
and where the player I and player II continuation regions are given by:
\begin{align}
    &D_1=\{x\in\mathbb{R}; 0<x<\tilde{x}\},
    \\ &D_2=\{x\in\mathbb{R}; x>\hat{x}\},
\end{align}
where the constants $a,\hat{x}$ and $\tilde{x}$ are determined by (\ref{ch2_example_a}), (\ref{hat_x_ch2_example_solution}) and (\ref{ch_2_example_p1_cont_region_value}) respectively and the constants $c_{\pm}$ are given by (\ref{c_values_ch2_example}).
\end{example}
\begin{example}\subsubsection*{The Optimal Liquidity Control and Lifetime Ruin Problem }\label{8}

We now revisit the optimal liquidity control and lifetime ruin problem in Section \ref{section_investment_problem} and solve the model. In the following analysis, we use the results of the stochastic differential game of impulse control and stopping to solve our model. 

Before stating results, using (\ref{generator}), (\ref{ch2firmliquidityprocess}) and (\ref{ch2invliquidityprocess}), we firstly make the following observation on the stochastic generator $\mathcal{L}^\theta$ which is given by the following expression:
\begin{align}
&\mathcal{L}^\theta\Phi(s,\cdot)
\\&\begin{aligned}= &erx\frac{\partial \Phi}{\partial x}(s,\cdot)+\Gamma  y \frac{\partial \Phi}{\partial y}(s,\cdot)+\frac{1}{2}\sigma^2_fx^2\frac{\partial^2 \Phi}{{\partial x}^2}(s,\cdot)+\frac{1}{2}\pi^2\sigma^2_Iy^2\frac{\partial^2 \Phi}{{\partial y}^2}(s,\cdot)
+\frac{1}{2}q^2\frac{\partial^2 \Phi}{\partial q^2}(s,\cdot)
\\&\begin{aligned}+\int_{\mathbb{R}}\Big\{\Phi(s,x+x\gamma_f(z),y, q-q\theta_1(z))-\Phi(s,\cdot)-x\gamma_f(z)\frac{\partial \Phi}{\partial x}(s,\cdot)+q\theta_1(z)\frac{\partial \Phi}{\partial q}(s,\cdot)\Big\}\nu(dz)&
\\-\sigma_fxq \frac{{\partial}^2 \Phi}{\partial x\partial q}(s,\cdot).&
\label{stochgeneratoropinvA}
\end{aligned}
\end{aligned}
\end{align}

The following result provides a complete characterisation of the investor's value function:
\begin{theorem}\label{capital_injection_value_function_double_obstacle} 
\

\noindent The investor's problem reduces to the following double obstacle variational inequality:
\begin{equation}
\inf\left\{\sup\left[\psi(s,\cdot)-(\kappa_I+\alpha_I(\hat{y}-y)),-\left(\frac{\partial}{\partial_s}+\mathcal{L}^{\hat{\theta}}\right)\psi(s,\cdot)\right],\psi(s,\cdot)-G(s,\cdot)\right\}=0, \label{chliqcontroldoubleobs}
\end{equation}
where $G(s,x,y,q)=e^{-\delta s}(g_1xq+\lambda_T+g_2y)$; the constants $g_1\in ]0,1]$ and $\lambda_T\geq0$ represent the fraction of the firm's liquidity process and some fixed amount each received by the investor upon exit respectively, $\hat{y}\in\mathbb{R}$ is an endogenous constant and lastly $q\in \mathbb{R}$ is the value of a stochastic process (later described in Lemma \ref{Lemma 7.4.}).
\end{theorem}
Theorem \ref{capital_injection_value_function_double_obstacle} establishes that the complete problem facing the investor can be written as a double obstacle problem from which the value function for the investment can be computed. Explicit solutions for the problem can be derived in cases in which the investor's wealth and firm liquidity processes do not contain jumps see Sec. \ref{sec_anal_solvability}.  

We now state the main theorem for this example:

\begin{theorem}\label{Theorem 7.2.} 
\

\noindent For the investor's optimal liquidity control and exit problem, the sequence of optimal capital injections $(\hat{\tau}, \hat{Z})\equiv [\hat{\tau}_j,\hat{z}_j ]_{j\in \mathbb{N}}\equiv\sum_{j\geq 1}\hat{z}_j  \cdot 1_{\{\hat{\tau}_j\leq \hat{\rho}\wedge T \}}  (s)$ is characterised by
the investment times $\{\hat{\tau}_j\}_{j\in\mathbb{N}}$ and magnitudes $\{\hat{z}_j\}_{j\in\mathbb{N}}$ where $[\hat{\tau}_j,\hat{z}_j ]_{j\in \mathbb{N}}$ are constructed recursively via the following expressions: 
\renewcommand{\theenumi}{\roman{enumi}}
\begin{enumerate}
\item$\hat{\tau}_0\equiv t_0\text{ and  }\hat{\tau}_{j+1}=\inf\{s>\tau_j;Y^{(\hat{\tau}, \hat{Z})_{[t_0,s]} } (s)\geq \tilde{y}|s\in\mathcal{T}\}\wedge \hat{\rho}$,
\item $\hat{z}_j=\hat{y}-y(\hat{\tau}_j)$,
\end{enumerate}
where the duplet $(\tilde{y},\hat{y})\in\mathbb{R}\times\mathbb{R}$ consists of endogenous constants.

The investor's non-investment region is given by:
\begin{equation}
D_2=\left\{y<\tilde{y}\right|y\in\mathbb{R}_{>0}\}.
\end{equation}

The optimal exit time $\hat{\rho}\in\mathcal{T}$ for the investor is given by:
\begin{equation}
\hat{\rho}=\inf\left\{s\geq t_0;\left(X(s),Q(s)\right)\notin D_1|s\in\mathcal{T}\right\}\wedge \tau_S, \end{equation}
where $Q$ is a stochastic process (c.f. Lemma \ref{Lemma 7.4.}) 
and the set $D_1$ (non-stopping region) is defined by:
\begin{equation}
D_1=\{xq> \omega^{\star}|x\in\mathbb{R}_{>0},q\in\mathbb{R}\},
\label{non_stopping_region_D1}
\end{equation} 
where $\omega^\star\in\mathbb{R}$ is an endogenous constant.
\end{theorem}

Theorem \ref{Theorem 7.2.} says that the investor performs discrete capital injections over a sequence of intervention times $\{\hat{\tau}_k\}_{k\in \mathbb{N}}$ over the time horizon of the problem. The decision to invest is determined by the investor's wealth process --- in particular, whenever the investor's wealth process reaches $\tilde{y}$, then the investor performs capital injections of magnitudes $\{\hat{z}_k\}_{k\in\mathbb{N}}$ to increase the firm's liquidity levels in order to provide the firm with maximal liquidity to perform market operations. Therefore, the value $\tilde{y}$ represents an investment threshold. This in turn maximises the liquidity that the investor makes available to the firm whilst the investor remains in the market after which the investor liquidates all investment holdings. However, if the firm's liquidity process exits the region $D_1$, in order to avoid the prospect of loss on investment, the investor immediately exits the market by liquidating all market holdings in the firm.

The non-stopping region $D_1$ is defined by (\ref{non_stopping_region_D1}) and the function $\psi$ is the investor's value function.  We later (Proposition \ref{Proposition 7.5}) provide a full characterisation of the investor's value function and the set of endogenous constants.

From Theorem \ref{Theorem 7.2.} we also arrive at the following result that enables us to state the exact points at which the investor performs an injection, when the investor exits the market and when the investor does nothing:

\begin{corollary}\label{Corollary 7.3.}
\

\noindent For the optimal liquidity control and lifetime ruin problem, the investor's wealth process $X$ lies within a space that splits into three regions: a region in which the investor performs a capital injection --- $I_1$, a region in which no action is taken --- $I_2$ and lastly a region in which the investor exits the market by selling all firm holdings --- $I_3$. Moreover, the three regions are characterised by the following expressions:
\begin{align*}
I_1&=\{y\geq \tilde{y},qx> \omega^{\star}|x,y\in\mathbb{R}_{>0},q\in\mathbb{R}\},\\
I_2&=\{qx> \omega^{\star}, y< \tilde{y}|x,y\in\mathbb{R}_{>0},q\in\mathbb{R}\},\\
I_3&=\{qx\leq \omega^{\star}|x\in\mathbb{R}_{>0},q\in\mathbb{R}\},
\end{align*}
where $\tilde{y},\omega^{\star}$ are fixed endogenous constants.
\end{corollary}
Lemma \ref{Lemma 7.4.} provides an expression for the process $Q$:
\begin{lemma}\label{Lemma 7.4.} 
\

\noindent The process $Q$ is determined by the expression:
\begin{align}
Q(s)= Q(0)\exp{\Big\{\frac{1}{2}\sigma^2_fs-\sigma_fB_f(s)+\int^s_0\int_\mathbb{R}(\ln (1+\hat{\theta}_1(r,z))-\hat{\theta}_1(r,z))\tilde{N}_f(dr,dz)\Big\}}&,\nonumber\\
 \label{ch2qdefn}\quad \forall s\in [0,T],&
\end{align}
where $\hat{\theta}_1$ is a solution to the equation $H(\psi)=0$ where $H$ is given by:
\begin{equation}
H(\psi)=\int_\mathbb{R}([\Xi(\psi(z))]^k-1)\nu(dz), \label{foctheta1final}
\end{equation}
where $\Xi(\psi(z)):=(1-\hat{\theta}_1(z))(1+\gamma_f(z)))$ and $k$ is an endogenous constant. 
\end{lemma}

The following result provides a complete characterisation of the investor's value function and the set of endogenous constants within the problem:
\begin{proposition}\label{Proposition 7.5}  
\

\noindent The value function $\psi:[0,T]\times \mathbb{R}_{>0}\times \mathbb{R}_{>0} \times \mathbb{R}\to \mathbb{R}$ for the investor's joint problem (\ref{investor_stopping_problem}) - (\ref{investor_control_problem}) is given by:
  \begin{align} 
\psi(s,x,y,q)= \begin{cases}
\begin{aligned}
&A_1(s,x,y,q),    &\left(\mathbb{R}\backslash{\partial  
 D_2}\right)\cap D_1 \\
  &A_2(s,x,y,q),&\mathbb{R}\backslash{\partial D_1}\\
  &A_3(s,x,y,q),&D_1\cap D_2
  \end{aligned}
  \end{cases}
  \label{valuefunctionopinvA}
  \end{align} 
where $A_1,A_2,A_3$ are given by:
\begin{align*}
A_1(s,x,y,q)&:= e^{-\delta s}q\left\{c\left(y^{d_1}-y^{d_2}\right)-q^{-1}\left(\kappa_I+\alpha_I\left(\hat{y}-y\right)\right)+ax^kq^k\right\},\\
A_2(s,x,y,q)&:=e^{-\delta (T\wedge\hat{\rho})}\left(g_1xq+\lambda_T+g_2y\right),\\
A_3(s,x,y,q)&:=qe^{-\delta s}\left(c(y^{d_1}-y^{d_2})+ax^kq^k\right),\end{align*}
where the constants $\delta, \kappa_I, \alpha_I$ are the investor's discount factor, the fixed part of the transaction cost and the proportional part of the transaction cost respectively and the constants $a,d_1, d_2$ and $\omega^{\star}$ are given by:
\begin{align}
  \omega^{\star}&=\frac{\lambda_Tk}{g_1(1-k)} 
  \\
  a&=\left(\frac{g_1}{k}\right)^k\left(\frac{\lambda_Tk}{1-k}\right)^{1-k} 
\\d_1&=\frac{1}{2}-\frac{1}{\pi^2\sigma_I^2}\left(\sqrt{(\Gamma-\frac{1}{2}\pi^2\sigma_I^2)^2+2\pi^2\sigma_I^2\delta}+\Gamma\right)\\
d_2&=\frac{1}{2}+\frac{1}{\pi^2\sigma_I^2}\left(\sqrt{(\Gamma-\frac{1}{2}\pi^2\sigma_I^2)^2+2\pi^2\sigma_I^2\delta}-\Gamma\right). 
\end{align}
The constant $k$ is a solution to the equation $p(k)=0$ where the function $p$ is given by:
\begin{equation}
p(k):= -\delta+(er-\sigma_f^2)k+k\int_{\mathbb{R}}(\hat{\theta}_1(z)-\gamma_f(z))\nu(dz),
\end{equation}
and lastly the constants $c,\hat{y}_2,\tilde{y}_2$ are determined by the set of equations:
\begin{align}
\tilde{y}_2^{d_1}-\hat{y}_2^{d_1}+\tilde{y}_2^{d_2}-\hat{y}_2^{d_2}&=c^{-1}(\alpha_I(\tilde{y}_2-\hat{y}_2)-\kappa_I)\\
d_1\hat{y}_2^{d_1-1}-d_2\hat{y}_2^{d_2-1}&=\alpha_Ic^{-1}\\
d_1\tilde{y}_2^{d_1-1}-d_2\tilde{y}_2^{d_2-1}&=\alpha_Ic^{-1}.
\end{align}
\end{proposition}
Proposition \ref{Proposition 7.5} therefore provides a complete characterisation of the value function for the investor's problem and the endogenous constants appearing in Theorem \ref{Theorem 7.2.} - Lemma \ref{Lemma 7.4.}.

\subsection*{Analytic Solvability of the Investment Problem} 
For the case without jumps, we can compute an exact closed analytic expression for the value function which is presented in Lemma \ref{Lemma 7.6.} (see Sec. \ref{sec_anal_solvability}). For the general case in which jumps are included, an analytic solution is not available. However, numerical approximations to solutions of quasi-variational HJBI equations are accessible through finite difference approximation schemes. In particular, under certain stability conditions, the Howard policy iteration algorithm can be shown to converge to the optimal strategy for the impulse control problem. 

Such matters are discussed in \cite{lapeyre2003understanding,kushner1990} and in particular, the numerical approximation of solutions to the quasi-variational inequality arising from impulse control problems is discussed extensively in \cite{azimzadeh2017impulse}. A proof of the convergence the Howard policy iteration algorithm for a general class of problems under which the current impulse control problem falls is discussed is given in \cite{chancelier2007policy}.  

\begin{refproof}[Proof of Theorem \ref{Theorem 7.2.}.]
\

\noindent We prove the theorem by applying Theorem \ref{Verification theorem for Non-Zero-Sum Stochastic Differential Games of Control and Stopping with Impulse Controls} to the model. We separate the proof into components that study the investor's capital injections problem, the investor's optimal stopping individually before combing the calculations to construct the full solution to the problem. 

The scheme of the proof is as follows:
\begin{enumerate}[label=\;\;\;\;\;\uline{\textit{Step \arabic*:}}]
    \item derive the functional form for the value function.
    \item characterise the non-investment region and optimal capital injection.
    \item characterise the continuation (non-exit) region and exit criterion.
    \item compute value function for the complete problem and show that the value function is a solution to a double obstacle problem.
\end{enumerate}
We then finalise with some remarks on the solution and discuss the cases where the underlying processes contain no jumps and the corresponding closed analytic solutions.

We wish to fully characterise the optimal investment strategies for the investor. To put problem (\ref{ch2divinvestorproblemcontrol}) - (\ref{ch2divinvestorproblemtime})  in terms of the framework of Theorem \ref{Verification_theorem_for_Zero-Sum_Stochastic_Differential_Games_of_Control_and_Stopping_with_Impulse_Controls}, we firstly note that we now seek the triplet $(\hat{\theta},(\tau,Z),\hat{\rho})$  s.th.
\begin{equation}
J^{\hat{\rho},\hat{u},\hat{\theta}}(t,y_1,y_2,y_3)=\sup_{\rho\in\mathcal{T}}\Big(\inf_{(\tau,Z)\in\Phi}\Big(\inf_{\theta}J^{\rho,u,\theta}(t,y_1,y_2,y_3)\Big)\Big),
\end{equation}
where
\begin{align}\nonumber
 J^{\rho,u,\theta}(t,y_1,y_2,y_3)=  \mathbb{E}\Bigg[-\sum_{m\geq 1}&e^{-\delta \tau_m} [\kappa_I+\alpha_I (\tau_m ) z_m ] \cdot 1_{\{\tau_m\leq\tau_S \}} \\&+e^{-\delta (\tau_S\wedge\rho)}\left(g_1{Y_1}_{\tau_S \wedge \rho }^{t,y_1,(\tau, Z) }Y_3^{t,y_3,(\theta_0,\theta_1)}+\lambda_T+g_2{Y_2}_{\tau_S \wedge \hat{\rho} }^{t,y_2,(\tau, Z) }\right)\Bigg],\label{ch2divinvestorpayoffch3}	
\end{align}
and $\theta\equiv(\theta_0,\theta_1):[0,T]\times\Omega\times[0,T]\times\Omega\to\Theta\subset\mathbb{R}^2$ and the dynamics of the state processes $Y:=(Y_0,Y_1,Y_2,Y_3)$ are expressed via the following:
\begin{align}
&dY_0(s)\equiv dt,\; \forall s\in [0,T],\label{chliqcontrolstateprocess1}
\\&
dY_1(s)\equiv dX_s^{t_0,x_0,(\tau,Z)};\; X_{t_0}^{t_0,x_0,\cdot}=y_1,\; \forall s\in [0,T], \label{chliqcontrolstateprocess2}
\\&
dY_2(s)\equiv dY_s^{t_0,y_0,(\tau,Z)};\; Y_{t_0}^{t_0,y_0,\cdot}=y_2,\; \forall s\in [0,T],\label{chliqcontrolstateprocess3}
\\&
dY_3^{t_0,y_3,(\theta_0,\theta_1)}(s)=-Y^{t_0,y_3,(\theta_0,\theta_1)}_3(s)\left[\theta_0(s)dB_f(s)+\int_{\mathbb{R}}\theta_1(s,z)\tilde{N}_f(ds,dz)\right]\label{chliqcontrolstateprocess4},
\end{align}
so that $Y_1, Y_2$ are processes which represent the firm liquidity processes and the investor's wealth process respectively. The processes $Y_0$ and $Y_3^{t_0,y_3,(\theta_0,\theta_1)}$ represent time and market adjustments to the investor's wealth process respectively.

In this section, we suppress the indices on the process and write $Y_3\equiv Y_3^{t_0,y_3,(\theta_0,\theta_1)}$. We also occasionally employ the following shorthands:
$\frac{\partial \phi}{\partial y_i}\equiv \partial_{y_i}\phi, \frac{\partial^2 \phi}{\partial y_i\partial y_j}\equiv \partial_{y_i,y_j}\phi$ for $i\in\{0,1,2,3\}$.

Lastly, we have the following relations for the state process coefficients:
 \begin{equation}
\mu(\cdot,y_2)=\Gamma y_2,\hspace{3 mm} \mu(\cdot,y_1)=er y_1.
\end{equation} 

We restrict ourselves to the case when:
\begin{equation}
\gamma_I(\cdot,y_2)=0. \label{ch3investorconditions}
\end{equation}
For the case that includes jumps (i.e $\gamma_I\not\equiv 0,\theta_1\not\equiv 0$) we impose a condition on the firm's discounted rate of return and the discount rate, that is we assume the following condition holds:
\begin{equation}
e >\Big(\frac{r}{\sigma^2_f}\Big)^{-1}. \label{ch3investorcondition1}
\end{equation}
The continuation regions $D_2$ and $D_1$ for the controller and the stopper respectively now take the form:
\begin{align}
D_2&=\left\{y\in[0,T]\times\mathbb{R}_{>0}\times\mathbb{R}_{>0}\times \mathbb{R};\psi(y)<\mathcal{M}  \psi(y)\right\},\\
D_1&=\left\{(y_0,y_1,y_2)\in[0,T]\times\mathbb{R}_{>0}\times\mathbb{R}_{>0}; \psi(y_0,y_1,y_2,\cdot)-G(y_0,y_1,y_2)>0\right\},
\end{align}
where given some $\phi\in\mathcal{H}$ the intervention operator $\mathcal{M}:\mathcal{H}\to\mathcal{H}$ is given by:
\begin{equation}
\mathcal{M}  \phi(y_0,y_1,y_2,y_3)=\inf_{\zeta\in\mathcal{Z}}\left\{\phi(y_0,y_1,y_2-\zeta,y_3)-(\kappa_I+\alpha_I\zeta), \zeta>0\right\},
\end{equation}
for all $y\in[0,T]\times \mathbb{R}_{>0}\times\mathbb{R}_{>0}\times\mathbb{R}$ and the stopping time $\hat{\rho}$ is defined by:
\begin{equation}
\hat{\rho}=\inf\{y_0>t_0; (y_0,y_1,y_2)\notin D_1|y_0\in \mathcal{T}\}\wedge\tau_S.
\end{equation}
\underline{\textit{Step 1.}}
Our first task is to characterise the value of the game. Now by conditions (ii) - (vi) of Theorem \ref{Verification_theorem_for_Zero-Sum_Stochastic_Differential_Games_of_Control_and_Stopping_with_Impulse_Controls}, we observe that the following expressions must hold:
\begin{align}
&\psi(y_0,y_1,y_2,y_3)=e^{-\delta y_0}\left(g_1{y_1} y_3+\lambda_T+g_2y_2\right), &&\hspace{-2.5 mm}\forall (y_0, y_3) \in [0,T]\times  \mathbb{R};\forall (y_1,y_2)\in \mathbb{R}_{>0}^2\backslash{D_1}     &&\hspace{-3 mm}\text{(condition (ii))}\label{psiterminalcondition}
\\&
\psi(y_0,y_1,y_2,y_3)\geq e^{-\delta y_0}\left(g_1{y_1} y_3+\lambda_T+g_2y_2\right),   &&\hspace{-3.5 mm}\forall y\in [0,T]\times \mathbb{R}_{>0}\times\mathbb{R}_{>0}\times\mathbb{R},     &&\hspace{-3 mm}\text{(condition (v))}
\\
&\frac{\partial \psi}{\partial y_0}+\mathcal{L}^\theta\psi(y)\geq 0, &&\hspace{-3.5 mm}\forall y\in [0,T]\times \mathbb{R}_{>0}\times\mathbb{R}_{>0}\times\mathbb{R},      &&\hspace{-3 mm}\text{(condition (vi))}
\\&
\inf_{\theta}\left\{\frac{\partial \psi}{\partial y_0}+\mathcal{L}^\theta\psi(y)\right\}=0, &&\hspace{-3.5 mm}\forall y\in   [0,T]\times D_1\times D_2\times\mathbb{R}, &&
\nonumber\\& && &&\hspace{-3 mm}\text{(condition (vi))}\label{chliqcontrolconditioneq}
\end{align}
where the condition labels refer to the conditions of Theorem \ref{Verification_theorem_for_Zero-Sum_Stochastic_Differential_Games_of_Control_and_Stopping_with_Impulse_Controls}.
Now using (\ref{chliqcontrolstateprocess1}) - (\ref{chliqcontrolstateprocess4}) we find that the generator is given by the following expression:
\begin{align}\nonumber
\mathcal{L}^\theta\psi(y)=& ery_1\frac{\partial \psi}{\partial y_1}(y)+\Gamma  y_2 \frac{\partial \psi}{\partial y_2}(y)+\frac{1}{2}\sigma^2_fy_1^2\frac{\partial^2 \psi}{\partial y_1^2}(y)+\frac{1}{2}\pi^2\sigma^2_Iy_2^2\frac{\partial^2 \psi}{\partial y_2^2}(y)
\\&\begin{aligned}
+\frac{1}{2}\theta^2_0y^2_3\frac{\partial^2 \psi}{\partial y_3^2}(y)
-\theta_0y_1y_3\sigma_f \frac{{\partial}^2 \psi}{\partial y_1\partial y_3}(y)+\int_{\mathbb{R}}\Big\{\psi(s,y_1+y_1\gamma_f(z),y_2, y_3-y_3\theta_1(z))&
\\-\psi(y)-y_1\gamma_f(z)\frac{\partial \psi}{\partial y_1}(y)+y_3\theta_1(z)\frac{\partial \psi}{\partial y_3}(y)\Big\}\nu(dz)&,
\end{aligned} \label{chliqcontrolgeneratror1}
\end{align}
\begin{equation}
\sup_{\rho\in\mathcal{T}}\left[\inf_{(\tau,Z)\in\Phi}\left(\inf_{\theta}J^{\rho,(\tau,Z),\theta}(y)\right)\right]=0.
\end{equation}
By (\ref{chliqcontrolconditioneq}) and (\ref{chliqcontrolgeneratror1}) we readily deduce that the first order condition on $\hat{\theta}_0$ is given by the following expression:
\begin{equation}
\hat{\theta}_0y_3^2\frac{\partial^2 \psi}{\partial y_3^2}-y_1y_3\sigma_f\frac{\partial^2 \psi}{\partial y_1\partial y_3} =0,
\end{equation}
which after some simple manipulation we find that:
\begin{equation}
\hat{\theta}_0 = y_1y_3^{-1}\sigma_f\partial^2_{y_1,y_3}\psi(\partial_{y_3}^2\psi)^{-1}. \label{chliqcontrolthetanought1}
\end{equation}
Now by (vi) of Theorem \ref{Verification_theorem_for_Zero-Sum_Stochastic_Differential_Games_of_Control_and_Stopping_with_Impulse_Controls} we have that on $D_1$:
\begin{equation}
\frac{\partial \psi}{\partial s}+\mathcal{L}^\theta\psi=0,
\end{equation}
(here $f=0$) which implies that:
\begin{align}
0=& \frac{\partial \psi}{\partial y_0}(y)+ery_1\frac{\partial \psi}{\partial y_1}(y)+\Gamma  y_2 \frac{\partial \psi}{\partial y_2}(y)+\frac{1}{2}\sigma^2_fy_1^2\frac{\partial^2 \psi}{{\partial y_1}^2}(y)+\frac{1}{2}\pi^2\sigma^2_Iy_2^2\frac{\partial^2 \psi}{{\partial y_2}^2}(y)
\\&
\begin{aligned}-\frac{1}{2}\sigma^2_f y_1^2\Bigg[\frac{\partial^2 \psi}{\partial y_1\partial y_3}(y)\Bigg]^2\Big(\frac{\partial^2 \psi}{\partial y_3^2}(y)\Big)^{-1}
+\int_{\mathbb{R}}\Big\{\psi(y_0,y_1+y_1\gamma_f(z),y_2, y_3-y_3\theta_1(z))&
\\-\psi(y)-y_1\gamma_f(z)\frac{\partial \psi}{\partial y_1}(y)+y_3\theta_1(z)\frac{\partial \psi}{\partial y_3}(y)\Big\}\nu(dz)&.\label{chliqcontrolgeneratror2}
\end{aligned}
\end{align}
Let us try as our candidate function $\psi(y)=e^{-\delta y_0}y_3(\phi_2(y_2)+\phi_\omega(\omega))$, where $\omega:=y_1y_3$. Then after plugging our expression for $\psi$ into (\ref{chliqcontrolgeneratror2}) we find that:
\begin{align}\nonumber
 0&=-\delta [\phi_2(y_2)+\phi_\omega(\omega)]+er\omega\phi_\omega'(\omega)+\Gamma \phi_2(y_2)+\frac{1}{2}\sigma_f^2\omega^2\phi_\omega''(\omega)+
 \\&\begin{aligned}
 +\frac{1}{2}\pi^2\sigma_I^2y_2^2\phi_2''(y_2)+\frac{1}{2}\hat{\theta}_0^2\omega(2\phi_\omega'(\omega)+\omega\phi_\omega''(\omega))-\hat{\theta}_0\sigma_f\omega(2\phi_\omega'(\omega)+\omega\phi_\omega''(\omega))&
 \\+\int_{\mathbb{R}}\Big\{(1-\theta_1(z))\big[\phi_2(y_2)+\phi_\omega(\omega(1+\gamma_f(z))(1-\theta_1(z)))\big]&
 \\-(1-\theta_1(z))[\phi_2(y_2)+\phi_\omega(\omega)]+\omega\phi_\omega(\omega)'(\theta_1(z)-\gamma_I(z))\Big\}\nu(dz)&, \label{chliqcontrolgeneratror3}
\end{aligned}
\end{align}
and (\ref{chliqcontrolthetanought1}) now becomes:
\begin{equation}
\hat{\theta}_0=\sigma_f\frac{y_1(2y_3\phi'_\omega(\omega)+y_3\omega\phi''_\omega(\omega)}{y_3(2y_1\phi'_\omega(\omega)+y_1\omega\phi''_\omega(\omega))}=\sigma_f. \label{chliqcontrolthetanought2}
\end{equation}
Hence, substituting (\ref{chliqcontrolthetanought2}) into (\ref{chliqcontrolgeneratror3}) we find that:
\begin{align}\nonumber
0=&-\delta[\phi_\omega(\omega)+\phi_2(y_2)]+er+\sigma^2_f)\omega\phi'_\omega(\omega)+\Gamma y_2\phi'_2(y_2)+\frac{1}{2}\pi^2\sigma_I^2y_2^2\phi_2''(y_2)
\\&\begin{aligned}
&+\int_{\mathbb{R}}\Big\{(1-\theta_1(z))\big[\phi_2(y_2)+\phi_\omega(\omega(1+\gamma_f(z))(1-\theta_1(z)))\big]
\\&-(1-\theta_1(z))[\phi_2(y_2)+\phi_\omega(\omega)]+\omega\phi_\omega'(\omega)(\theta_1(z)-\gamma_I(z))\Big\}\nu(dz). \label{chliqcontrolcontequationfullcomp}
\end{aligned}
\end{align}
Additionally, our first order condition on $\hat{\theta}_1$ becomes:
\begin{equation}
\int_{\mathbb{R}}\Big\{\phi_\omega(\omega\Xi(\hat{\theta}_1(z)))+\omega\Xi(\hat{\theta}_1(z))\phi_\omega'(\omega\Xi(\hat{\theta}_1(z)))-\phi_\omega-\omega\phi'_\omega\Big\}\nu(dz)=0, \label{ch3thetahateqnaftersubs}
\end{equation} where $\Xi(\hat{\theta}_1(z)):=(1-\hat{\theta}_1(z))(1+\gamma_f(z)))$.

We can decouple (\ref{chliqcontrolcontequationfullcomp}) after which we find that when $(y_1,y_2)\in D_1\times D_2$ we have that:
\begin{enumerate}[label={\roman*)}]
\item{\begin{align}\nonumber \hspace{-4 mm}\int_{\mathbb{R}}\Big\{(1&-\theta_1(z))\big[\phi_\omega(\omega(1+\gamma_f(z))(1-\theta_1(z)))\big]-(1-\theta_1(z))\phi_\omega(\omega)\\&+\omega\phi_\omega'(\omega)(\theta_1(z)-\gamma_I(z))\Big\}\nu(dz)
 -\delta\phi_\omega(\omega)+(er-\sigma^2_f)\omega\phi'_\omega(\omega)=0,\label{chliqcontrolphi1solequation}
 \end{align}}
 \item{$\qquad\qquad 
 -\delta\phi_2(y_2)+\Gamma y_2\phi'_2(y_2)+\frac{1}{2}\pi^2\sigma_I^2y_2^2\phi_2''(y_2)=0.
 \label{chliqcontrolphi2unsolved.}$}
\end{enumerate}
We can solve the Cauchy-Euler equation \ref{chliqcontrolphi2unsolved.} --- after performing some straightforward calculations we find that:
\begin{equation}
\phi_2(y_2)=c_1y_2^{d_1}+c_2y_2^{d_2}, \label{phi2ansatz.0.0}
\end{equation}
for some (as yet undetermined) constants $c_1$ and $c_2$. The constants $d_1$ and $d_2$ are given by:
\begin{align}
d_1&=\frac{1}{2}-\frac{1}{\pi^2\sigma_I^2}\left(\sqrt{(\Gamma-\frac{1}{2}\pi^2\sigma_I^2)^2+2\pi^2\sigma_I^2\delta}+\Gamma\right)\label{d1phi2}\\
d_2&=\frac{1}{2}+\frac{1}{\pi^2\sigma_I^2}\left(\sqrt{(\Gamma-\frac{1}{2}\pi^2\sigma_I^2)^2+2\pi^2\sigma_I^2\delta}-\Gamma\right). \label{d2phi2}
\end{align}
Since $\psi(0)=Y_2(0)=0$, we easily deduce that $c_2=-c_1$, after which we deduce that $\phi_2$ is given by the following expression:
\begin{equation}
\phi_2(y_2)=c(y_2^{d_1}-y_2^{d_2}), \label{phi2ansatz.0}    
\end{equation}
where $c:=c_1=-c_2$ is some as of yet undetermined constant.

To obtain an expression for the function $\phi_\omega$, in light of (\ref{chliqcontrolphi1solequation}) we conjecture that $\phi_\omega$ takes the form:
\begin{equation}
\phi_\omega=a\omega^k, \label{phiomegaansatz}
\end{equation} where $a$ and $k$ are some constants. Using (\ref{phiomegaansatz}) and (\ref{chliqcontrolphi1solequation}), we find the following:
\begin{equation}
\mathcal{L}^{q}\phi_\omega(\omega)=a\omega^kp(k),
\end{equation}
where the operator $\mathcal{L}^q$ is defined by the following expression for some function $\phi\in\mathcal{C}^{1,2}$:
\begin{align}\nonumber
&\mathcal{L}^q\phi(w):= -\delta\phi_w(w)+(er-\sigma^2_f)w\phi'_w(w)
\\&\hspace{-1 mm}+\int_{\mathbb{R}}\Big\{(1-\hat{\theta}_1(z))\big[\phi(w(1+\gamma_f(z))(1-\hat{\theta}_1(z)))-\phi(w)\big]+w\phi'(w)(\hat{\theta}_1(z)-\gamma_f(z))\Big\}\nu(dz),
\end{align} and $p(k)$ is defined by:
\begin{align}
p(k):= -\delta+(er-\sigma^2_f)k+\int_{\mathbb{R}}\Big\{(1-\hat{\theta}_1(z))[\Xi(c\hat{\theta}_1(z))^k-1]+k(\hat{\theta}_1(z)-\gamma_f(z))\Big\}\nu(dz),\label{pkequation1}
\end{align} where  $\Xi(\hat{\theta}_1(z)):=(1-\hat{\theta}_1(z))(1+\gamma_f(z))).$\

Note that using (\ref{phiomegaansatz}), the first order condition on $\hat{\theta}_1$ (c.f. (\ref{ch3thetahateqnaftersubs})) now becomes:
\begin{equation}
\int_\mathbb{R}(\Xi^k-1)\nu(dz)=0. \label{foctheta1final2}
\end{equation}
Hence using (\ref{foctheta1final2}), (\ref{pkequation1}) becomes:
\begin{equation}
p(k)= -\delta+(er-\sigma^2_f)k+k\int_{\mathbb{R}}(\hat{\theta}_1(z)-\gamma_f(z))\nu(dz).\label{pkequation2}
\end{equation}
We now make the following observations:
\begin{equation}
p(0)=- \delta <0, \hspace{3 mm} p(b)> 1-\delta+b\int_\mathbb{R}(\theta_1(z)-\gamma_f(z))\nu(dz)\geq 0,   \mathbb{P}-\text{a.s.},
\end{equation}
for any $b> (er-\sigma^2_f)^{-1}$ using condition (\ref{ch3investorcondition1}) and the fact that $\delta\in ]0,1]$.
We therefore deduce the existence of a value $z\in ]0,1[$ s.th. $p(z)=0$. We now conclude that:
\begin{equation}
\phi_\omega(\omega)=a\omega^{k}, \label{phiomegaansatz2}
\end{equation}
where $a$ is an arbitrary constant and where $k$ is a solution to the equation:
\begin{equation}
p(k)=0.
\end{equation}
We now split the analysis into two parts in which we study the investor's capital injections (impulse control) problem and the investor's optimal stopping problem separately. We then later recombine the two problems to construct our solution to the problem.\\

\noindent\underline{\textit{Step 2: The Investor's Capital Injections Problem}}\\

We firstly tackle the investor's capital injections problem, in particular we wish to ascertain the form of the function $\phi_2$ and describe the intervention region and the optimal size of the investor's capital injections.

Our ansatz for the continuation region $D_2$ is that it takes the form:
\begin{equation}
D_2=\{y_2< \tilde{y}_2,|y_2,\tilde{y}_2\in \mathbb{R}\}. \label{contregionansatzopinv}
\end{equation}
Therefore by (ii) of Theorem \ref{Verification_theorem_for_Zero-Sum_Stochastic_Differential_Games_of_Control_and_Stopping_with_Impulse_Controls} for $y_2 \notin D_2$ we have that:
\begin{align}
\psi(s,y)=\mathcal{M}\psi(s,y)=\inf\left\{\psi(s,y_1,y_2-\zeta,y_3)+(\kappa_I+\alpha_I\zeta), \zeta>0\right\}\nonumber\\
\iff \hspace{3 mm} \phi_2(y_2)=\inf\left\{\phi_2(y_2-\zeta)+(\kappa_I+\alpha_I\zeta), \zeta>0\right\}.\label{chliqcontrolheqn1}\end{align}

Let us define the function $h$ by the following expression:
\begin{equation}
h(\zeta)=\phi_2(y_2-\zeta)-(\kappa_I+\alpha_I\zeta).
\end{equation}
We therefore see that the first order condition for the minima $\hat{\zeta}(y_2)\in \mathcal{Z}$ of the function $h$ is:
\begin{equation}
h'(\zeta)=\phi'_2(y_2-\hat{\zeta})=\alpha_I.
\end{equation}

Let us now consider a unique point $\hat{y}_2\in]0,\tilde{y}_2[$ s.th.
\begin{equation}
\phi_2'(\hat{y}_2)=\alpha_I, \label{chliqcontrolphi2deriveqn1}
\end{equation}
and
\begin{equation}
\hat{y}_2=y_2-\hat{\zeta}(y_2) \text{  or  } \zeta(\hat{y}_2)=\hat{y}_2-y_2.
\end{equation}
After imposing a continuity condition at $y_2=\tilde{y}_2$, by (\ref{chliqcontrolheqn1}) we have that:
\begin{equation}
  \phi_2(\tilde{y}_2)=\phi_{2,0}(\hat{y}_2)-(\kappa_I+\alpha_I(\hat{y}_2-\tilde{y}_2)), \label{chliqcontrolphi1hatstareqn}
\end{equation}
where $\phi_{2,0}(y_2)=\phi_2(y_2)$ on $D_2$ and where $\phi_2$ is given by (\ref{phi2ansatz.0}). 
Additionally, by construction of $\tilde{y}_2$ we have that:
\begin{equation}
\phi_2'(\tilde{y}_2)=\alpha_I. \label{chliqcontrolphi2deriveqn2}
\end{equation}
Hence we deduce that the function $\phi_2$ is given by the following expression:
 \begin{equation}
\phi_2(y_2)= \begin{cases}
c\left(y^{d_1}_2-y^{d_2}_2\right)-\left(\kappa_I+\alpha_I\left(\hat{y}_2-y_2\right)\right), \hspace{5 mm}     y_2\geq \tilde{y}_2 \\
  c\left(y^{d_1}_2-y^{d_2}_2\right),\hspace{5 mm}y_2< \tilde{y}_2,\label{ch3phi2tab} 
  \end{cases}
  \end{equation}
where $d_1$ and $d_2$ are given by (\ref{d1phi2}) - (\ref{d2phi2}).

In order to compute the constants $a$, $\hat{y}_2$ and $\tilde{y}_2$, we use the system of equations (\ref{chliqcontrolphi2deriveqn1}), (\ref{chliqcontrolphi1hatstareqn}) and (\ref{chliqcontrolphi2deriveqn2}).\\

\noindent\underline{\textit{Step 3: The Investor's Optimal Stopping Problem}}\\

Our ansatz for the continuation region $D_1$ is that it takes the form:
\begin{equation}
D_1=\{\omega=y_1y_3> y_1^{\star}y_3^{\star}=\omega^{\star}|y_1,y_1^{\star}\in \mathbb{R};y_3,y_3^{\star}\in\mathbb{R}\}.
\end{equation}

If we assume that the \textit{high contact principle}\footnote{Recall that the high contact principle is a condition that asserts the continuity of the value function at the boundary of the continuation region. In the current case this implies that $\phi_\omega(\omega)|^{\omega=\omega^{\star}}=G(\omega)|^{\omega=\omega^{\star}},\; \frac{\partial }{\partial \omega}\phi_\omega(\omega)|^{\omega=\omega^{\star}}=\frac{\partial }{\partial \omega}G(\omega)|^{\omega=\omega^{\star}}$} holds, in particular if we have differentiability at $\omega^{\star}$ then, using (\ref{phiomegaansatz}) we obtain the following equations:
\renewcommand{\theenumi}{\roman{enumi}}
 \begin{enumerate}[leftmargin= 6.5 mm]
\item \hspace{50 mm}$ a\omega^{*k}=g_1\omega^{\star}+\lambda_T \label{phi2hcpr1case1}$
\item \hspace{50 mm}$ak\omega^{*k-1}=g_1,$  \label{phi2hcpr2case2}
\end{enumerate}
by continuity and differentiability at $\omega^{\star}$. Since the system of equations (\ref{phi2hcpr1case1}) - (\ref{phi2hcpr2case2}) completely determine the constants $a$ and $\omega^{\star}$, we can compute the values of $\omega^{\star}$ and $a$ in (\ref{phiomegaansatz}), after which we find:
\begin{equation}
\omega^{\star}=\frac{\lambda_Tk}{g_1(1-k)}, \hspace{3 mm} a=\Big(\frac{g_1}{k}\Big)^k\Big(\frac{\lambda_Tk}{1-k}\Big)^{1-k}.
\end{equation}
  
We are now in a position to prove Lemma \ref{Lemma 7.4.}; in particular, using (\ref{chliqcontrolstateprocess4}) and (\ref{chliqcontrolthetanought2}) we now see that the process $Y_3$ is determined by the expression:
\begin{equation}
dY_3(s)=-\left[\sigma_fY_3dB_f(s)+Y_3(s)\int_{\mathbb{R}}\hat{\theta}_1(s,z)\tilde{N}_f(ds,dz)\right],\quad \mathbb{P}-a.s.,\label{ch3y3solved}
\end{equation}
where $\hat{\theta}_1$ is determined by the equation (c.f. (\ref{foctheta1final2})):
\begin{equation}
\int_\mathbb{R}(\Xi^k(\hat{\theta}_1(z))-1)\nu(dz)=0, \label{foctheta1final3}
\end{equation}
where $\Xi(\hat{\theta}_1(z)):=(1-\hat{\theta}_1(z))(1+\gamma_f(z)))$. 

Using It\={o}'s formula for It\={o}-L\'evy processes, we can solve (\ref{foctheta1final3}), moreover since\\ $\mathbb{E}_{\mathbb{Q}}\Bigg[X+\lambda_T\Bigg]=\mathbb{E}_{\mathbb{P}}\Bigg[\Big(X+\lambda_T\Big)Y_3\Bigg]$, (c.f. (\ref{ch2divinvestorproblemtime})), the process $Y_3$ represents the Radon-Nikodym derivative of the measure $\mathbb{Q}$ with respect to the measure $\mathbb{P}$ (i.e. $Y_3(s)=\frac{d(\mathbb{Q}|\mathcal{F}_s)}{d(\mathbb{P}|\mathcal{F}_s)}$). Combining these two statements and denoting $Y_3$ by $Q$ immediately gives the result stated in Lemma \ref{Lemma 7.4.}. $\hfill \square$\\

\noindent\underline{\textit{Step 4: The Investor's Value Function and Joint Problem}}\\

Our last task is to combine the results together and fully characterise the investor's value function. We firstly note that putting the above results together yields the following double obstacle variational inequality:
\begin{align}
\sup\left\{\inf\left[\psi(y)-(\kappa_I+\alpha_I(\hat{y}_2-y_2)),-\left(\frac{\partial}{\partial_{y_0}}+\mathcal{L}^{\hat{\theta}}\right)\psi(y)\right],\psi(y)-G(y)\right\}=0, 
\end{align}
where $y=(y_0,y_1,y_2,y_3)$ and $G(y)=e^{-\delta y_0}(g_1y_1y_3+g_2y_2+\lambda_T)$ and the investor's stopping time is given by:
\begin{align}
\hat{\rho}=\inf\left\{s\geq t_0; Y_1(s)Y_3(s)\notin D_1|s\in \mathcal{T}\right\},
\end{align}
where the stochastic generator $\mathcal{L}^{\hat{\theta}}$ acting on a test function $\psi\in\mathcal{C}^{1,2}$ is defined via the following expression:
\begin{align}\nonumber
\mathcal{L}^\theta\psi(y)&= ery_1\frac{\partial \psi}{\partial y_1}(y)+\Gamma  y_2 \frac{\partial \psi}{\partial y_2}(y)+\frac{1}{2}\sigma^2_fy_1^2\frac{\partial^2 \psi}{{\partial y_1}^2}(y)+\frac{1}{2}\pi^2\sigma^2_Iy_2^2\frac{\partial^2 \psi}{{\partial y_2}^2}(y)
\\&\nonumber
\begin{aligned}
+\frac{1}{2}\theta^2_0y^2_3\frac{\partial^2 \psi}{\partial y_3^2}(y)
-\theta_0y_1y_3\sigma_f \frac{{\partial}^2 \psi}{\partial y_1\partial y_3}(y)+\int_{\mathbb{R}}\Big\{\psi(y_0,y_1+y_1\gamma_f(z),y_2, y_3-y_3\theta_1(z))\\\qquad\qquad\qquad-\psi(y)-y_1\gamma_f(z)\frac{\partial \psi}{\partial y_1}(y)+y_3\theta_1(z)\frac{\partial \psi}{\partial y_3}(y)\Big\}\nu(dz),
\end{aligned}
\end{align}
and where the function $\hat{\theta}_1$ satisfies the first order condition (\ref{foctheta1final3}).

The double obstacle problem in (\ref{chliqcontroldoubleobs}) characterises the value for the game, this proves Theorem \ref{capital_injection_value_function_double_obstacle}. 

Proposition \ref{Proposition 7.5} provides a full expression of the value function for the investor's problem. To prove Proposition \ref{Proposition 7.5}, we need to collect the results on the constituent functions of the value function and assemble the complete function. Combining (\ref{phiomegaansatz}) and (\ref{ch3phi2tab}) and using (\ref{psiterminalcondition}) shows that the value function $\psi$ is given by:
 \begin{equation}
\psi(y)= 
\begin{cases}
\begin{aligned}
&A_1(y), \;&     (\mathbb{R}\backslash{\partial  
 D_2}) \cap D_1\\
  &A_2(y),& \mathbb{R}\backslash{\partial D_1}
  \\&A_3(y),& D_1\cap D_2,
  \end{aligned}
  \end{cases}
\end{equation} 
where
\begin{align*}
    A_1(y)&:=e^{-\delta y_0}y_3\left(c\left(y^{d_1}_2-y^{d_2}_2\right)-y^{-1}_3\left[\kappa_I+\alpha_I(\hat{y}_2-y_2)\right]+ay_1^ky_3^k\right),\\
    A_2(y)&:=e^{-\delta (T\wedge\hat{\rho})}\left(g_1y_1y_3+\lambda_T+g_2y_2\right),\\
    A_3(y)&:=y_3e^{-\delta y_0}\left(ay_1^ky_3^k+c(y^{d_1}_2-y^{d_2}_2)\right).
\end{align*}
The constants $a,\omega^{\star}$ are given by:
\begin{equation}
  \omega^{\star}=\frac{\lambda_Tk}{g_1(1-k)}, \hspace{3 mm} a=\left(\frac{g_1}{k}\right)^k\left(\frac{\lambda_Tk}{1-k}\right)^{1-k}, \label{omegastarastarfinal}
\end{equation}
and the constants $d_1$ and $d_2$ are given by:
\begin{align}
d_1&=\frac{1}{2}-\frac{1}{\pi^2\sigma_I^2}\left(\sqrt{(\Gamma-\frac{1}{2}\pi^2\sigma_I^2)^2+2\pi^2\sigma_I^2\delta}+\Gamma\right)\\
d_2&=\frac{1}{2}+\frac{1}{\pi^2\sigma_I^2}\left(\sqrt{(\Gamma-\frac{1}{2}\pi^2\sigma_I^2)^2+2\pi^2\sigma_I^2\delta}-\Gamma\right). 
\label{d1-d2}
\end{align}
The constants $c,\hat{y}_2,\tilde{y}_2$ are determined by the set of equations:
\begin{align}
\tilde{y}_2^{d_1}-\hat{y}_2^{d_1}+\tilde{y}_2^{d_2}-\hat{y}_2^{d_2}&=c^{-1}(\alpha_I(\tilde{y}_2-\hat{y}_2)-\kappa_I)\label{aimpulseeq1}\\
d_1\hat{y}_2^{d_1-1}-d_2\hat{y}_2^{d_2-1}&=\alpha_Ic^{-1}\\
d_1\tilde{y}_2^{d_1-1}-d_2\tilde{y}_2^{d_2-1}&=\alpha_Ic^{-1},\label{aimpulseeq3}
\end{align}
and the constant $k$ is a solution to the equation $p(k)=0$ where the function $p$ is given by\footnote{For the case that includes jumps in the firm liquidity process, we assume that the firm's discounted rate of return is greater than $1$ and the discount rate is relatively small compared to the volatility parameter $\sigma_f$ as given in condition (\ref{ch3investorcondition1}).}:
\begin{equation}
p(k):= -\delta+(er-\sigma_f^2)k+k\int_{\mathbb{R}}(\hat{\theta}_1(z)-\gamma_f(z))\nu(dz),
\label{ch3pkequation}\end{equation}
where $\hat{\theta}_1$ is a solution to (\ref{foctheta1final}). This proves Proposition \ref{Proposition 7.5}.
$\hfill \square$\end{refproof}

\subsubsection*{The Case $\mathbf{\gamma_I= 0,\theta_1= 0}$}
\label{sec_anal_solvability}

If the investor's liquidity process contains no jumps (i.e. $\gamma_f\equiv 0$ and  $\theta_1\equiv 0$ in (\ref{ch2firmliquidityprocess}) and (\ref{chliqcontrolstateprocess4}) (resp.)) then we can obtain closed analytic solutions for the parameters of the function $\phi_\omega$. Indeed, when $\gamma_I\equiv 0$ and $\theta_1\equiv 0$ using (\ref{ch3pkequation}), we see that the expression for  $p(k)$ reduces to:
\begin{equation}
p(k):= -\delta+(er-\sigma_f^2)k.
\end{equation}
We can therefore solve for $k$ after which we find that the function $\phi_\omega$ is given by: 
\begin{align}
\phi_\omega(\omega)= a\omega^{k},
\end{align}
where $k=\delta(er-\sigma^2_f)^{-1}$. The constants $a,\omega^{\star},d_1,d_2$ are determined by (\ref{omegastarastarfinal}) - (\ref{d1-d2}) and the constants $c,\hat{y},\tilde{y}$ are determined by the set of equations:
\begin{align}
\tilde{y}^{d_1}-\hat{y}^{d_1}+\tilde{y}^{d_2}-\hat{y}^{d_2}&=c^{-1}(\alpha_I(\tilde{y}-\hat{y})-\kappa_I)\\
d_1\hat{y}^{d_1-1}-d_2\hat{y}^{d_2-1}&=\alpha_Ic^{-1}\\
d_1\tilde{y}^{d_1-1}-d_2\tilde{y}^{d_2-1}&=\alpha_Ic^{-1}.
\end{align}

We therefore arrive at the following result which provides a complete characterisation of the value function in terms of a closed analytic solution for the investor's problem when the liquidity process contains no jumps: 
\begin{lemma}\label{Lemma 7.6.}
\

\noindent For the case in which the investor's liquidity process contains no jumps (i.e. $\gamma_f\equiv 0$ in (\ref{ch2firmliquidityprocess})), the function $\psi$ is given by the following:
\begin{align}
\psi(y)=
\begin{cases}
\begin{aligned}
&A_1'(y), &(\mathbb{R}\backslash{\partial  
 D_2}) \cap D_1\\
  &A_2'(y), &\mathbb{R}\backslash{\partial D_1}\\
  &A_3'(y),  &D_1\cap D_2,
  \end{aligned}
 \end{cases}
 \end{align}  
where 
\begin{align*}
A_1'(y)&:=e^{-\delta y_0}y_3(c(y^{d_1}_2-y^{d_2}_2)-y_3^{-1}(\kappa_I+\alpha_I(\hat{y}-y_2))+ay_1^ky_3^k)),
\\
A_2'(y)&:= e^{-\delta (T\wedge\hat{\rho})}(g_1y_1y_3+\lambda_T+g_2y_2),
\\
A_3'(y)&:=y_3e^{-\delta y_0}(c(y^{d_1}_2-y^{d_2}_2)+ay_1^ky_3^k),
\end{align*}
where $k=\delta(er-\sigma_f^2)^{-1},\;
\omega^{\star}=\lambda_Tkg_1^{-1}(1-k)^{-1}, \; a=g_1^{k}k^{-k}(\lambda_Tk)^{1-k}(1-k)^{-(1-k)}.\label{omegaequationch2}
$
\end{lemma}
Lastly, the process $Q$ is determined by the expression:
\begin{equation}
Y_3(t)= Y_3(0)\exp{\Big\{\frac{1}{2}\sigma^2_ft-\sigma_fB_f(t)\Big\}}, \qquad \forall t\in [0,T].
\end{equation}
$\hfill \square$
\subsubsection*{The General Case ($\mathbf{\gamma_I\neq 0,\theta_1\neq 0}$)}

Though obtaining a closed analytic solution to $p(k)=0$ represents a difficult task, the solution may be approximated using numerical methods. 

The following results follow directly from the degeneracy of the game:
\begin{corollary}\label{Corollary 2.5.5.}
\

\noindent Consider the above problem when investment set $\Phi$ is a singleton, the game collapses to an optimal stopping problem when the investor seeks to minimise risk of ruin.\footnote{ This is a specific (Markovian and for which the set of probability measures is dominated) case of the game considered in for example \cite{nutz2015optimal}.} In this case the double-obstacle problem reduces to:
\begin{equation}
\inf\Big\{-\Big(\frac{\partial}{\partial_s}+\mathcal{L}^{\hat{\theta}}\psi(y)\Big),\psi(y)-G(y)\Big\}=0. \label{chliqcontroldoubleobsnocontrol}
\end{equation}
\end{corollary}
\begin{corollary}\label{Corollary 2.5.6.}
Consider the above problem, when $\theta\equiv 0$, the game collapses to a problem of impulse control with discretionary stopping.
\end{corollary}
\end{example}

\section*{Appendix}\label{appendix}
\noindent \textbf{Technical Conditions for (T1) - (T4)}.\

\renewcommand{\theenumi}{\roman{enumi}}
 \begin{enumerate}[leftmargin= 6 mm]

\item[(T1)] Assume that $\mathbb{E}[\int_0^T 1_{\partial D} (X^{\cdot,u} (s))ds]=0\;\; \forall X \in S,\forall u \in \mathcal{U}$ where $D\equiv D_1\cup D_2$.

\item[(T2)] $\partial D$ is a Lispchitz surface --- that is to say that $\partial D$ is locally the graph of a Lipschitz continuous function: $\phi \in \mathcal{C}^2 (S\backslash \partial D) $ with locally bounded derivatives.\

\item[(T3)] The sets $\{\phi^- (X^{\cdot,u} (\tau_m));$ $ \tau_m \in \mathcal{T},\forall  m \in \mathbb{N}\}$ and $\{\phi^- (X^{\cdot,u} (\rho));\rho \in \mathcal{T}\}$ are uniformly integrable $\forall  x \in S$ and $\forall u \in \mathcal{U}$.

	\item[(T4)] $\mathbb{E}[|\phi(X^{\cdot,u} (\tau ))|+|\phi(X^{\cdot,u} (\rho))|+\int_{0}^{T} |\mathcal{L}\phi(X^{\cdot,u} (s))|ds]< \infty,$ $ \\   \text{for all intervention times } \tau_,\rho \in \mathcal{T}$ and $\forall u \in \mathcal{U}.$
\end{enumerate}

\medskip

\printbibliography[
heading=bibintoc,
title={References}]

\end{document}